\documentclass[11pt]{article}
\usepackage[pagebackref,colorlinks=true,pdfpagemode=none,urlcolor=blue,linkcolor=blue,citecolor=blue]{hyperref}
\usepackage{subfigure}
\usepackage{amsthm, amssymb, bm}
\usepackage{soul, color}
\usepackage[]{amsmath}
\usepackage[]{amsfonts}
\usepackage[]{fancyhdr}
\usepackage[]{graphicx}
\usepackage[]{empheq}
\usepackage[toc,page]{appendix}
\newtheorem{theorem}{Theorem}[section]
\newtheorem{lemma}[theorem]{Lemma}

\newtheorem{remark}[theorem]{Remark}

\makeatletter

\newcommand{\Rmnum}[1]{\expandafter\@slowromancap\romannumeral #1@}
\makeatother

\def \Cm {\mathbb{C}}
\def \Imm {\mathbb{I}}

\def \Rm {\mathbb{R}}
\def \Sm {\mathbb{S}}

\newcommand{\bk}{\mathbf k}

\newcommand{\bfe}{ {\bf e}}

\newcommand{\tr}{ {\text{tr }}}

\newcommand{\cout}[1]{}

\def \bmrho {{\boldsymbol\rho}}
\def \bmphi {{\boldsymbol\varphi}}

\title{Imaging of anisotropic conductivities from current densities in two dimensions}

\author{Guillaume Bal\thanks{Department of Applied Physics and Applied Mathematics, Columbia University,  New York NY, 10027; gb2030@columbia.edu} \and Chenxi Guo\thanks{Department of Applied Physics and Applied Mathematics, Columbia University,  New York NY, 10027; cg2597@columbia.edu} \and Fran\c cois Monard\thanks{Department of Mathematics, University of Washington, Seattle WA, 98195; fmonard@uw.edu}}

\begin{document}
\maketitle
\begin{abstract}
We consider the imaging of anisotropic conductivity tensors $\gamma=(\gamma_{ij})_{1\leq i,j\leq 2}$ from knowledge of several internal current densities $\mathcal{J}=\gamma\nabla u$ where $u$ satisfies a second order elliptic equation $\nabla\cdot(\gamma\nabla u)=0$ on a bounded domain $X\subset \Rm^2$ with prescribed boundary conditions on $\partial X$. We show that $\gamma$ can be uniquely reconstructed from four {\em well-chosen} functionals $\mathcal{J}$ and that noise in the data is differentiated once during the reconstruction. The inversion procedure is local in the sense that (most of) the tensor $\gamma(x)$ can be reconstructed from knowledge of the functionals $\mathcal{J}$ in the vicinity of $x$. We obtain the existence of an open set of boundary conditions on $\partial X$ that guaranty stable reconstructions by using the technique of complex geometric optics (CGO) solutions. The explicit inversion procedure is presented in several numerical simulations, which demonstrate the influence of the choice boundary conditions on the stability of the reconstruction. This problem finds applications in the medical imaging modality called Current Density Imaging or Magnetic Resonance Electrical Impedance Tomography.
%
%
\end{abstract}
\section{Introduction}
Current Density Impedance Imaging (CDII), also called Magnetic Resonance Electrical Impedance Tomography (MREIT) is a medical imaging technique that belongs to the class of coupled-physics imaging modalities. Such modalities aim to combine a high-contrast modality, such as Electrical Impedance Tomography (EIT), with a high-resolution modality, such as Magnetic Resonance Imaging (MRI) or ultrasound. EIT, which aims to reconstruct the electrical properties of tissues, leads to a nonlinear inverse boundary problem known as the Calder{\'o}n problem, which has been extensively studied (see \cite{Uhl2009} for a review). The Calder\'on problem consists in understanding what can be reconstructed in $\gamma$ from knowledge of all possible measurements performed at the boundary $\partial X$ of a domain $X\in\Rm^n$, in other words from knowledge of all pairs $(u,\gamma\nu\cdot\nabla u)$ at the boundary $\partial X$ with outward unit normal $\nu(x)$ for $x\in\partial X$, where $u$ is an arbitrary solution to the equation $\nabla\cdot \gamma\nabla u=0$ in $X$. It is known that (i) anisotropic tensors $\gamma$ cannot uniquely be reconstructed from such boundary data; and (ii) that when $\gamma=\beta\tilde\gamma$ with ${\rm det}\tilde\gamma=1$ and $\tilde\gamma$ known, then $\beta$ can uniquely be reconstructed with a stability estimate that is logarithmic; see \cite{Uhl2009}. This stability estimate, which intuitively corresponds to differentiating noise in the data an infinite number of times, results in typically low-resolution reconstructions.


In a new class of inverse problems (called hybrid or coupled-physics inverse problems), it is possible to overcome the limitations of reconstructions from classical boundary data by acquiring and using internal functionals of the coefficients of interest. These high-resolution internal functionals of the unknown conductivity allow for the high-resolution reconstruction of a fully anisotropic structure.  For a very incomplete list of works on these problems in the mathematical literature, we refer the reader to \cite{Ammari2008,AS-IP-12,Ren2011,Bal2010,Stefanov2012}. Different types of internal functionals, such as \emph{current densities} and \emph{power densities}, corresponding to different physical couplings have been analyzed to recover the unknown conductivity. In the case of power densities, which share some similarities with the problem of interest here, we refer the reader to, e.g.,\cite{Bal2012e,Capdeboscq2009,Kuchment2011a,Kuchment2011,Monard2012a,Monard2011}. 

In this paper,  we consider the problem of reconstructing an anisotropic conductivity $\gamma$ in a domain $X$ from measurement of internal current densities $H$. Internal current densities can be obtained by the technique of Current Density Imaging (CDI). The idea is to use Magnetic Resonance Imaging (MRI) to determine the magnetic field $B$ induced by an input current $I$ (see \cite{Ider1997}). The current density is then defined by $H=\nabla \times B$.  The explicit reconstructions we propose require that all components of $H$ be measured. This may be challenging in practice as it requires a rotation of the domain being imaged or of the MRI scanner. The reconstruction of $\gamma$ from knowledge of only some components of $H$, ideally only one component for the simplest practical experimental setup, is open at present.

In two dimensions, a numerical reconstruction algorithm based on the construction of equipotential lines was given in \cite{Kwon2002}. An iterative algorithm known as \emph{J}-substitution algorithm was proposed by Kwon \emph{et al} in \cite{Kwon2002a}. With knowledge of the magnitude of only one current density $|H|=|\gamma\nabla u|$, the problem was studied in \cite{Nachman2007,Nachman2009} in the isotropic case. The anisotropic case in a known conformal class was studied in \cite{HMN2013}. The present authors recently derived a local reconstruction procedure for fully anisotropic tensors in \cite{Guo2013a} and showed that the tensor can be uniquely and stably reconstructed with a loss of one derivative from the measurements to the reconstructed image. The result was also extended by the first two authors to the full Maxwell's system in \cite{Guo2013b}.

The explicit reconstruction method provided in \cite{Guo2013a} requires that some matrices constructed from available data satisfy appropriate conditions of linear independence. In the present work, we show that in $\Rm^2$, such assumptions can be globally guaranteed with a set of well-chosen illuminations based on the construction of Complex Geometrical Optics (CGO) solutions, provided that one can prescribe Dirichlet (or other) conditions over the full boundary. Several numerical experiments presented in Section \ref{numerical} confirm the theoretical predictions. The numerical simulations show that the reconstruction procedure works well for different types of tensors containing both smooth and discontinuous coefficients. Using the decomposition $\gamma=\beta\tilde\gamma$ with $\beta=(\det\gamma)^{\frac{1}{2}}$, the simulation results also show that both the isotropic and the anisotropic parts of the tensor can be stably reconstructed, with a better robustness to noise for the scalar $\beta$. This is consistent with theoretical results in \cite{Guo2013a}, where the stability of the inversion on $\beta$ is better than on the anisotropy $\tilde\gamma$. 

Our CGO-based theoretical results exhibit a specific class of boundary conditions that ensure stable reconstructions. In practice, a much larger class of boundary conditions than those that can be analyzed mathematically still provide stable reconstructions. Yet, when only a part of the boundary conditions is accessible for current injection, the linear independence of specific matrices needed in the reconstruction deteriorates. The reconstructions then become unstable in some parts of the domain. This phenomenon is demonstrated in several numerical simulations. All simulations are performed in two dimensions of space, although we expect the conclusions to still hold qualitatively in higher dimensions as well.

 
 The rest of the paper is structured as follows. The main results are presented in section \ref{sec:res}. The reconstruction procedure is detailed in section \ref{sec:rec}. The numerical implementation of the algorithm and the effects of the choice of boundary conditions are shown in section \ref{numerical}. Section \ref{sec:conclu} offers some concluding remarks.

 
\section{Statements of the main results}\label{sec:res}
Let $X\subset \Rm^2$ be a bounded domain with a $C^{2,\alpha}$ boundary $\partial X$. Although most of the following results generalize to arbitrary spatial dimensions, we restrict the setting to $\Rm^2$; see \cite{Guo2013a} for results in higher dimensions. We consider the inverse problem of reconstructing an anisotropic conductivity tensor in the second-order elliptic equation,
\begin{align}
    \nabla\cdot(\gamma\nabla u) =0 \quad (X), \qquad u|_{\partial X} = g,
    \label{eq:conductivity}
\end{align}
from knowledge of internal current densities of the form $H=\gamma\nabla u$, where $u$ solves \eqref{eq:conductivity}. The above equation has real-valued coefficients and $\gamma=(\gamma_{ij})_{1\leq i,j\leq 2}$ is a symmetric (real-valued) tensor satisfying the uniform ellipticity condition 
\begin{align}
    \kappa^{-1}\|\xi\|^2\le \xi\cdot\gamma\xi \le \kappa\|\xi\|^2, \quad \xi\in \Rm^2, \quad \text{for some } \kappa\ge 1,
    \label{positive definite}
\end{align}
so that \eqref{eq:conductivity} admits a unique solution in $H^1(X)$ for $g\in H^{\frac{1}{2}}(\partial X)$. 
\subsection{Global reconstructibility condition}
We start by selecting $4$ boundary conditions $(g_1,g_2,g_3,g_4)$ and the corresponding current densities 
\begin{align}\label{measurement}
H_i = \gamma\nabla u_i, \quad 1\leq i\leq 4
\end{align}
where the function $u_i$ solves \eqref{eq:conductivity}.  Assuming that over $X$, the two solutions $u_1,u_2$ satisfy the following positivity condition
\begin{align}\label{posi condition}
\inf_{x\in X}|\det(\nabla u_1,\nabla u_2)|\geq c_1 > 0
\end{align}
then the gradients of additional solutions $\nabla u_3,\nabla u_4$ can be decomposed as linear combinations in the basis $(\nabla u_1,\nabla u_2)$,
\begin{align}\label{lnc}
\left\{\begin{array}{lll}
\nabla u_3 = \mu_1\nabla u_1+\mu_2\nabla u_2\\
\nabla u_4 = \lambda_1\nabla u_1+\lambda_2\nabla u_2
\end{array}\right.
\end{align}
where the coefficients $\{\mu_i\}_{1\leq i\leq 2}$ can be computed by Cramer's rule as
\begin{align}\label{lncoefficient}
(\mu_1,\mu_2)=(\frac{\det(\nabla u_3,\nabla u_2)}{\det(\nabla u_1,\nabla u_2)}, \frac{\det(\nabla u_1,\nabla u_3)}{\det(\nabla u_1,\nabla u_2)})=(\frac{\det(H_3,H_2)}{\det(H_1,H_2)},\frac{\det(H_1,H_3)}{\det(H_1,H_2)}).
\end{align}
The same expression holds for $\{\lambda_i\}_{1\leq i\leq 2}$ by replacing $u_3$ by $u_4$ in the above equation. Therefore these coefficients are computable from the available current densities. The reconstruction procedures will make use of the matrices $Z_k$ defined by 
\begin{align}
    Z_k=\left[Z_{k,1} ,Z_{k,2}\right],\quad \text{where }\quad Z_{1,i}=\nabla \mu_i \quad Z_{2,i}=\nabla \lambda_i, \quad 1\leq i, k\leq 2.
\label{Y Z}
\end{align}
These matrices are also uniquely determined by the known current densities. Denoting the matrix $H=[H_1,H_2]$ and the skew-symmetric matrix $J=\bfe_2\otimes\bfe_1-\bfe_1\otimes\bfe_2$, we construct two matrices as follows,
\begin{align}\label{constraint matrice}
M_k = (Z_kH^TJ)^{sym}, \quad \text{for}\quad k = 1,2.
\end{align}
The calculations in the following section show that condition \eqref{posi condition} and the independence of $M_1,M_2\in S_2(\Rm)$ give a sufficient condition for a global reconstruction of $\gamma$. Condition \eqref{posi condition} may be fulfilled using \cite[Theorem 4]{Alessandrini} which guarantees that \eqref{posi condition} holds if the map $\partial X\ni x\rightarrow(g_1(x),g_2(x))$ is a homeomorphism onto its image. That all required conditions are met for some boundary conditions is provided in the following lemma.
\begin{lemma}\label{cond cgo}
Let $\gamma(x)\in H^{5+\epsilon}(X)$ for some $\epsilon >0$ and satisfy the uniform elliptic condition \eqref{positive definite}. Then there exists an open set of illuminations $\{g_i\}_{1\leq i\leq 4}$, such that the solutions $\{u_i\}_{1\leq i\leq 4}$ satisfy the following conditions:
\begin{itemize}
	\item[A.] $\inf\limits_{x\in X}|\det(H_1,H_2)|\geq c_0 > 0$ holds on $X$. 
	\item[B.] The two matrices $M_1,M_2$ constructed by \eqref{constraint matrice} are independent in $S_2(\Rm)$ throughout $X$.
    \end{itemize}  
\end{lemma}
Since $\gamma$ is uniformly elliptic on $X$, condition {\em A} is equivalent to equation \eqref{posi condition}. Note that  {\em A} and {\em B} are expressed in terms of the measured quantities $\{H_j\}_j$, and as such can be checked directly during experiments. When the above constant $c_0$ is deemed too small, or the matrices $M_j$ are not sufficiently independent, then acquiring additional measurements might be considered.

The proof of Lemma \ref{cond cgo} is based on the construction of Complex Geometrical Optics(CGO) solutions and will be given in Section \ref{plemma}.
\begin{remark}
Conditions {\em A} and {\em B} are all that is required from the available internal functional $\{H_j\}$. The above lemma shows the existence of boundary conditions such that they hold. In practice, these conditions are met for a large class of boundary conditions not covered by the above lemma; see section \ref{numerical}.
\end{remark}
\begin{remark}
For the general $n$ dimensional case, Lemma \ref{cond cgo} does not necessarily hold globally. However, it holds locally with $4$ well-chosen illuminations. The proof is based on the Runge approximation; see \cite{Guo2013a} for details.
\end{remark}
\subsection{Uniqueness and stability results}
We denote by $M_2(\Rm)$ the space of $2\times 2$ matrices with inner product $\langle A,B\rangle:=\tr(A^TB)$. Assuming that there exist $4$ illuminations $\{g_i\}_{1\leq i\leq 4}$ with their corresponding solutions $(u_i)_{1\leq i\leq 4}$ satisfying the conditions in Lemma \ref{cond cgo}. Then the isotropic part $\beta$ can be reconstructed via a redundant elliptic system with a prior knowledge of the anisotropic part $\tilde\gamma$. In particular, the matrices $M_1,M_2$ constructed by \eqref{constraint matrice} are independent and of codimensiton $1$ in $S_2(\Rm)$. We will see that $\tilde\gamma$ is orthogonal to  $M_1,M_2$ which can be calculated from knowledge of $\{H_i\}_{1\leq i\leq 4}$. Together with the fact that $\det\tilde\gamma=1$ and $\tilde\gamma$ is positive, $\tilde\gamma$ can be completely determined by $(H_i)_{1\leq i\leq 4}$. The algorithm is based on an appropriate generalization of the cross-product. The reconstruction formulas can be found in Section \ref{retilde} and \ref{rebeta}. This algorithm leads to a unique and stable reconstruction in the sense of the following theorem.
\begin{theorem}\label{stability}
Suppose that Lemma \ref{cond cgo}.$A$ holds over $X$ for two couples $\{u_i\}_{i=1}^{2}$ and $\{u'_i\}_{i=1}^{2}$, solutions of the conductivity equation \eqref{eq:conductivity} with the tensors $\gamma=\beta\tilde\gamma$ and $\gamma'=\beta'\tilde\gamma'$ satisfying the uniform ellipticity condition \eqref{positive definite}, where $\tilde\gamma,\tilde\gamma'\in W^{1,\infty}(X)$ are known. Then $\beta$ can be uniquely reconstructed in $X$ with the following stability estimate,
\begin{align}\label{stability beta}
	\|\log\beta-\log\beta'\|_{W^{p,\infty}(X)}\le \epsilon_0+C \left( \sum_{i=1,2} \|H_i-H'_i\|_{W^{p,\infty}(X)}+\|\tilde\gamma-\tilde\gamma'\|_{W^{p,\infty}(X)} \right).
    \end{align}
    Here, $\epsilon_0=|\log\beta(x_0)-\log\beta'(x_0)|$ is the error committed at some fixed $x_0\in X$. If in addition Lemma \ref{cond cgo}.$B$ holds for the two sets $\{u_i\}_{i=1}^{4}$ and $\{u'_i\}_{i=1}^{4}$ as above, then $\tilde\gamma$ can be reconstructed with the stability as follows,
    \begin{align}
	\|\tilde\gamma - \tilde\gamma'\|_{W^{p,\infty}(X)}\le C \sum_{i=1}^{4} \|H_i-H'_i\|_{W^{p+1,\infty}(X)}.
	\label{eq:stabgammatilde}
    \end{align}
\end{theorem}
\begin{remark} 
From Theorem \ref{stability}, with a prior knowledge of the anisotropic part $\tilde\gamma$, the reconstruction of the scalar $\beta$ has a better stability estimate than $\tilde\gamma$. This will be demonstrated by the numerical experiments in Section \ref{NEX}.
\end{remark}

\section{Reconstruction approaches}\label{sec:rec}
The reconstruction approaches were presented in \cite{Guo2013a} for a general $n$ dimensional case. To make this paper self-contained, we briefly list the algorithm for the $2$ dimensional case and prove the \emph{global} reconstructibility condition in Lemma \ref{cond cgo}. We first present the reconstruction formula for $\beta$, assuming that the anisotropic part $\tilde\gamma$ is known from prior informations or reconstructed by current densities.  
\subsection{Reconstruction of $\beta$}\label{retilde}
Denoting the curl operator in $\Rm^2$ by $J\nabla\cdot$, where $J=\bfe_2\otimes\bfe_1-\bfe_1\otimes\bfe_2$. We rewrite \eqref{measurement} as $\frac{1}{\beta}\tilde{\gamma}^{-1}H_i=\nabla u_i$ for $i = 1,2$ and apply the curl operator to both sides. Using the fact that $\nabla u_i$ is curl free, we get the following equation,
\begin{align*}
\nabla\log\beta\cdot(J\tilde\gamma^{-1}H_i)=-J\nabla\cdot(\tilde\gamma^{-1}H_i).
\end{align*}
Considering both $j=1,2$, simple calculations lead to 
\begin{align}\label{nabla beta}
\nabla\log\beta=-J\tilde\gamma H^{-T}\left(\begin{array}{c} J\nabla\cdot(\tilde\gamma^{-1}H_1)\\  J\nabla\cdot(\tilde\gamma^{-1}H_2)\end{array}\right).
\end{align}
Since both first order derivatives of $\log\beta$ can be reconstructed by \eqref{nabla beta}, together with the boundary condition, the above equation leads to an over-determined elliptic system for $\beta$.
\subsection{Reconstruction of $\tilde\gamma$}\label{rebeta}
We now develop the reconstruction algorithm for $\tilde\gamma$. This reconstruction is algebraic and \emph{local} in nature: the reconstruction of $\gamma$ at $x_0\in X$ requires the knowledge of current densities for $x$ only in the vicinity of $x_0$. In addition to $H_1, H_2$, we pick 2 more measurements $H_3, H_4$ satisfying Lemma \ref{cond cgo}.$B$. We apply the curl operator $J\nabla\cdot$ to the linear combinations in \eqref{lnc}. Again, using the fact that $\nabla u_i = \gamma^{-1}H_i$ is curl free, we obtain the following equation, 
\begin{align*}
\sum_{i=1,2}Z_{k,i}\cdot(J\tilde\gamma^{-1}H_i)=0  \quad \text{where} \quad k=1,2.
\end{align*}
Using the fact that $\tr(A) = \tr(S^{-1}AS)$ and $\gamma$ is symmetric, the above equation amounts to
\begin{align*}
0= \tilde\gamma:Z_kH^TJ=\tilde\gamma:(Z_kH^TJ)^{sym}=\tilde\gamma:M_k.
\end{align*}
Since $\{M_1,M_2\}$ are of codimension $1$ in $S_2(\Rm)$, the above equation leads to the fact that $\tilde\gamma$ must be parallel to the following matrix constructed with $M_1,M_2$,
\begin{align}\label{B}
B=\left(\begin{array}{cc}
2M_1^{22}M_2^{12}-2M_1^{12}M_2^{22}& M_1^{11}M_2^{22}-M_1^{22}M_2^{11}\\
M_1^{11}M_2^{22}-M_1^{22}M_2^{11}    & 2M_1^{12}M_2^{11}-2M_1^{11}M_2^{12}
\end{array}\right).
\end{align}
Here, $M_k^{ij}$ denotes the $ij$ element of the symmetric matrix $M_k$. Notice that $B$ vanishes only if $M_1$ and $M_2$ are linearly dependent. Together with the fact that $\det\tilde\gamma=1$ and $\tilde\gamma$ is positive, we obtain the following explicit reconstruction, 
\begin{align}\label{re tilde}
\tilde\gamma= \text{sign}(B^{11})|B|^{-\frac{1}{2}}B.
\end{align}
\paragraph{Proof of Theorem \ref{stability}:}The proof is straightforward by noticing that one derivative is taken in the reconstruction procedure for $\tilde\gamma$. The stability for $\beta$ is a direct result from the standard regularity estimate for elliptic operators. See \cite{Guo2013a} for details.

\subsection{Proof of Lemma \ref{cond cgo}}\label{plemma} 
\paragraph{Isotropic tensors $\gamma=\beta\Imm_2$.}The proof is based on the construction of complex geometrical optics (CGO) solutions.  As shown in \cite{Bal2010}, letting $\beta\in H^{5+\varepsilon}(X)$, one is able to construct a complex-valued solution of \eqref{eq:conductivity} of the form 
    \begin{align}
	u_\bmrho = \frac{1}{\sqrt{\beta}} e^{\bmrho\cdot x} (1+ \psi_\bmrho),
	\label{eq:urho}
    \end{align}
where $\bmrho\in \Cm^2$ is of form $\bmrho = \rho (\bk + i\bk^\perp)$ with $\bk\in\Sm^1$ and $\bk\cdot \bk^\perp=0$. Thus $e^{\bmrho\cdot x}$ is a harmonic complex plane wave with $\bmrho\cdot\bmrho=0$. With the assumed regularity, we have the following estimate (see \cite[Proposition 3,3]{Bal2010}),
\begin{align*}
\lim_{\rho\rightarrow \infty}\|\psi_\bmrho\|_{\mathcal{C}^2(\bar X)}=0.
\end{align*}
Computing the gradient of $u_{\bmrho}$ and rearranging terms, we obtain that
 \begin{align*}
	\nabla u_{\bmrho} = e^{\bmrho\cdot x} (\bmrho + \bmphi_{\bmrho}), \quad\text{with} \quad \bmphi_\bmrho := \nabla\psi_\bmrho +\psi_\bmrho\bmrho- (1+\psi_{\bmrho} )\nabla \log\sqrt{\beta},
    \end{align*}
where $\|\bmphi_\bmrho\|_{\mathcal{C}^1(\bar X)}$ is uniformly bounded independent of $\bmrho$. Since $\beta$ is real-valued, both the real and imaginary parts of $u_{\bmrho}$ are real-valued solutions of \eqref{eq:conductivity} and we obtain the following expression
    \begin{align*}
	\nabla u_\bmrho^\Re &= \frac{\rho e^{\rho\bk\cdot x}}{\sqrt{\beta}} \left( (\bk + \rho^{-1} \bmphi_\bmrho^\Re) \cos(\rho \bk^\perp\cdot x) - (\bk^\perp + \rho^{-1} \bmphi_\bmrho^\Im) \sin (\rho\bk^\perp\cdot x) \right), \\
	\nabla u_\bmrho^\Im &= \frac{\rho e^{\rho\bk\cdot x}}{\sqrt{\beta}} \left(  (\bk^\perp + \rho^{-1} \bmphi_\bmrho^\Im) \cos(\rho \bk^\perp\cdot x) + (\bk + \rho^{-1} \bmphi_\bmrho^\Re) \sin (\rho\bk^\perp\cdot x)  \right).
    \end{align*}
Straightforward computations lead to 
    \begin{align*}
	\det(\nabla u_\bmrho^\Re,\nabla u_\bmrho^\Im )= \frac{\rho^2 e^{2\rho\bk\cdot x}}{\beta} ( 1 + f_{\bmrho} ), \quad \text{where} \quad \lim_{\rho\rightarrow \infty}\|f_\bmrho\|_{\mathcal{C}^1(\bar X)}=0.
    \end{align*}
Now we identify $\bk=\bfe_1$ and define $\bk_1=\bk$, $\bk_2=-\bk$.  For $j=1,2$, define $\bmrho_j:=\rho(\bk_j+i\bk_j^\perp)$. Considering the solutions $(u_{\bmrho_1}^\Re,u_{\bmrho_1}^\Im,u_{\bmrho_2}^\Re,u_{\bmrho_2}^\Im)$, the previous calculations show that 
\begin{align}
\inf_{x\in X}|\det(\nabla u_{\bmrho_1}^\Re,\nabla u_{\bmrho_1}^\Im )| \geq c_0>0,
\end{align}
Together with the uniform ellipticity of $\gamma$, the above inequality implies condition $A$  . Then using Cramer's rule in \eqref{lncoefficient}, simple algebra shows that 
\begin{align*}
\mu_1 = \frac{ \sin (2\rho\bk^\perp\cdot x)+g_{\mu_1}}{e^{2\rho\bk\cdot x}(1+f_{\bmrho_1})},\quad \mu_2 = \frac{ -\cos (2\rho\bk^\perp\cdot x)+g_{\mu_2}}{e^{2\rho\bk\cdot x}(1+f_{\bmrho_1})}
\end{align*}
and similarly,
\begin{align*}
\lambda_1 = \frac{ -\cos (2\rho\bk^\perp\cdot x)+g_{\lambda_1}}{e^{2\rho\bk\cdot x}(1+f_{\bmrho_1})},\quad \lambda_2 = \frac{ -\sin (2\rho\bk^\perp\cdot x)+g_{\lambda_2}}{e^{2\rho\bk\cdot x}(1+f_{\bmrho_1})}
\end{align*}
where $\|g_{\mu_i}\|_{\mathcal{C}^1(\bar X)},\|g_{\lambda_i}\|_{\mathcal{C}^1(\bar X)}$ are bounded for $i=1,2$. Then by the definition of $Z_k$ in \eqref{Y Z},  we obtain the following expression,
\begin{align*}
Z_1 &=2\rho e^{-2\rho\bk\cdot x}[(-\bk\sin (2\rho\bk^\perp\cdot x) +\bk^\perp\cos (2\rho\bk^\perp\cdot x),\bk\cos (2\rho\bk^\perp\cdot x) +\bk^\perp\sin (2\rho\bk^\perp\cdot x))+o (\rho^{-1})]\\
Z_2 &=2\rho e^{-2\rho\bk\cdot x}[(\bk\cos (2\rho\bk^\perp\cdot x) +\bk^\perp\sin (2\rho\bk^\perp\cdot x),\bk\sin (2\rho\bk^\perp\cdot x) -\bk^\perp\cos (2\rho\bk^\perp\cdot x))+o (\rho^{-1})].
\end{align*}
Together with $\bk=\bfe_1$ and $H=\beta(\nabla u_{\bmrho_1}^\Re,\nabla u_{\bmrho_1}^\Im )$, we obtain that,
\begin{align*}
(Z_1H^TJ)^{sym} &= 2\rho^2 \sqrt{\beta}e^{-\rho\bk\cdot x}\left[\left(\begin{array}{cc}
\cos (\rho\bk^\perp\cdot x) &\sin (\rho\bk^\perp\cdot x) \\
\sin (\rho\bk^\perp\cdot x)   & -\cos (\rho\bk^\perp\cdot x) 
\end{array}\right)+o (\rho^{-1})\right]\\
(Z_2H^TJ)^{sym} &= 2\rho^2 \sqrt{\beta}e^{-\rho\bk\cdot x}\left[\left(\begin{array}{cc}
\sin (\rho\bk^\perp\cdot x) &-\cos (\rho\bk^\perp\cdot x) \\
-\cos (\rho\bk^\perp\cdot x)   & -\sin (\rho\bk^\perp\cdot x) 
\end{array}\right)+o (\rho^{-1})\right].
\end{align*}
Since
\begin{align*}
\frac{M_1:M_2}{\|M_1\|\|M_2\|}=o (\rho^{-1})
\end{align*}
$M_1$, $M_2$ are almost orthogonal as $\rho$ is large enough, which implies the independence. This proves condition $B$.
\paragraph{General case:}Following the idea in \cite[Theorem 4.4]{BU-CPAM-12}, we extend $\gamma$ to a smooth tensor on $\Rm^2\simeq\mathbb{C}$, which remains uniformly positive definite and equal to $\Imm_2$ outside of a compact domain. For $\varphi:\Rm^2\ni x \mapsto \varphi(x) =y\in \Rm^2$ a diffeomorphism, we denote the push-forward of $\gamma$ by the $\varphi$ as follows,
\begin{align}
\varphi_*\gamma(y) = \frac{D\varphi(x)\gamma(x)D\varphi^t(x)}{|\det(D\varphi)|}|_{x=\varphi^{-1}(y)}
\end{align}
The theory of quasi-conformal mappings \cite{Astala} implies that there exists a unique such diffeomorphism $\varphi$ satisfying the \emph{Beltrami} system,
\begin{align*}
\varphi_*\gamma(y)=|\gamma|^{\frac{1}{2}}\circ\varphi^{-1}(y), \quad \varphi(z) =z+\mathcal{O}(z^{-1}) \quad \text{as}\quad z\rightarrow \infty,
\end{align*}
which means that the conductivity $\gamma$ is conformal to the Euclidean conductivity $\Imm_2$. As in the isotropic case, we can construct CGOs of the form,
  \begin{align}
	v_\bmrho = \frac{1}{\sqrt{\varphi_*\gamma(y)}} e^{\bmrho\cdot y} (1+ \psi_\bmrho(y)),
	\label{eq:urho}
    \end{align}
where $\lim_{\rho\rightarrow \infty}\|\psi_\bmrho\|_{\mathcal{C}^2(\varphi(X))}=0$. Using the change of variables, we construct $u=v\circ \varphi$,
 \begin{align}
	u_\bmrho = \frac{1}{\sqrt{\varphi_*\gamma\circ\varphi(x)}} e^{\bmrho\cdot \varphi(x)} (1+ \phi_\bmrho(x)),
	\label{eq:urho}
    \end{align}
where $\lim_{\rho\rightarrow \infty}\|\phi_\bmrho\|_{\mathcal{C}^2(X)}=0$. By the method in the isotropic case, we construct the solutions $(v_1,v_2,v_3,v_4)=(v_{\bmrho_1}^\Re,v_{\bmrho_1}^\Im,v_{\bmrho_2}^\Re,v_{\bmrho_2}^\Im)$, with ${\bmrho_{1},\bmrho_2}$ defined as before. Then for $1\leq i\leq 4$, the functions $u_i=v_i\circ\varphi$ satisfy the conductivity equation \eqref{eq:conductivity}. Using the chain rule $\nabla u_i = \nabla(v_i\circ\varphi)=D\varphi^t\nabla v_i\circ \varphi$, condition $A$ is satisfied with $\rho$ sufficiently large since $\nabla v_1,\nabla v_2$ are linearly independent as indicated in the isotropic case. Denote the skew-symmetric matrice $J'=D\varphi^tJ D\varphi$ and $Z'_k(y)=Z_k(x)|_{x=\varphi^{-1}(y)}$. Again by the chain rule, the following relation holds for every $x\in X$, 
\begin{align*}
(Z_kH^TJ')^{sym} &= (D\varphi^t Z'_k(\nabla v_1,\nabla v_2 )^t D\varphi\gamma D\varphi^tJ D\varphi)^{sym}\\
&= \det(D\varphi)D\varphi^t ((Z'_k\beta(\nabla v_1,\nabla v_2 )^t J)^{sym}\circ\varphi )D\varphi
\end{align*}
where $J'=D\varphi^tJ D\varphi$ is skew-symmetric. As in the proof in the isotropic case, we see that $(Z'_k\beta(\nabla v_1,\nabla v_2 )^t J)^{sym}$ are linearly independent over $\varphi(X)$ for $k= 1,2$. Thus, $M_1,M_2$ are linearly independent thoughout $X$, which proves condition $B$.

\section{Numerical experiments}\label{numerical}
To demonstrate the computational feasibility of the reconstruction algorithm, we performed some numerical experiments to validate the reconstruction algorithms from the previous section, assess their robustness to noisy measurements and determine how reconstructions are affected by boundary conditions limited to a part of the domain.


\subsection{Preliminary facts on the numerical implementation}  \label{NEX}

Recall that we decompose $\gamma$ into the following form with three unknown coefficients $\{\xi, \zeta, \beta\}$,
\begin{align}\label{3coef}
\gamma= \beta\tilde\gamma = \beta\left[\begin{array}{cc} \xi &\zeta\\ \zeta & \frac{1+\zeta^2}{\xi}\end{array}\right], \quad \xi>0
\end{align}
where $\beta=|\gamma|^{\frac{1}{2}}$ and $|\tilde\gamma|=1$. The full reconstruction is a two-step procedure, starting  with the reconstruction of the anisotropy $\tilde\gamma(\xi,\zeta)$ via formula \eqref{re tilde}. This requires implementing the formula
\begin{align}
\tilde\gamma= \frac{\sum_{i=1}^m\text{sign}(B_i^{11})B_i}{\sum_{i=1}^m{|\det B_i|^{\frac{1}{2}}}}
\end{align}
where each $B_i$ is constructed via \eqref{B} by choosing two additional current densities. Once $\tilde\gamma$ is reconstructed, $\beta$ is in turn reconstructed via the redundant elliptic system \eqref{nabla beta}. 

\paragraph{Regularized inversion.} Since we have explicit reconstruction formulas for $\gamma$, we use a total variation method as the denoising procedure by minimizing the following functional,
\begin{align}
f=\arg \min_g \frac{1}{2}\|g-f_{\text{rc}}\|^2_2 +\rho \|Mg\|_{\text{\scriptsize TV}}
\end{align}
where $f_{\text{rc}}$ denotes the explicit reconstructions of the coefficients of $\gamma$ and $M$ is the discretized version of the gradient operator. We choose the $l^1$-norm as the regularization TV norm for discontinuous, piecewise constant, coefficients. In this case, the minimization problem can be solved using the split Bregman method presented in \cite{Goldstein2009}. To recover smooth coefficients, we minimize the following least square problem,
\begin{align*}
    f=\arg \min_g \frac{1}{2}\|g-f_{\text{rc}}\|^2_2 +\rho \|Mg\|^2_2
\end{align*}
where the Tikhonov regularization functional admits an explicit solution $f=(\Imm+\rho M^*M)^{-1}M^*f_{\text{rc}}$. 

\subsection{Experiment with control over the full boundary}

In the numerical experiments below, we take the domain of interest to be the square $X=[-1,1]^2$ and use the notation $\mathbf{x}= (x,y)$. We use a $\mathsf{N+1\times N+1}$ square grid with $\mathsf{N}=80$, the tensor product of the equi-spaced subdivision $\mathsf{x =-1:h:1}$ with $\mathsf{h=2/N}$. The internal current densities $H(x)$ used are synthetic data that are constructed by solving the conductivity equation \eqref{eq:conductivity} using a finite difference method implemented with {\tt MatLab}. Although the data constructed this way may contain some noise, we refer to these data as the ``noise-free" or ``clean" data. 

We also perform the reconstructions with noisy data by perturbing the internal functionals $H(\mathbf{x})$ so that,
\begin{align*}
\widetilde {H}(\mathbf{x}) = H(\mathbf{x}).*(1+\alpha* \text{random}(\mathbf{x})),
\end{align*}
where $\text{random}(\mathbf{x})$ is a $\mathsf{N+1\times N+1}$ random matrix taking uniformly distributed values in $[-1,1]$ and  $\alpha$ is the noise level. We then run a de-noising process on the random matrix, which we chose as a low-pass filter constructed
by a 5-point sliding averaging process.

We use the relative $L^2$ error between reconstructed and true coefficients to measure the quality of the reconstructions. $\mathcal{E}^C_{\xi},\mathcal{E}^N_{\xi},\mathcal{E}^C_{\zeta},\mathcal{E}^N_{\zeta},\mathcal{E}^C_{\beta},\mathcal{E}^N_{\beta}$ denote the relative $L^2$ error in the reconstructions from clean and noisy data for $\xi$, $\zeta$ and $\beta$, respectively.

\begin{figure}[htp]
  \centering
  \subfigure[true $\xi$]{
      \includegraphics[width=37mm,height=35mm]{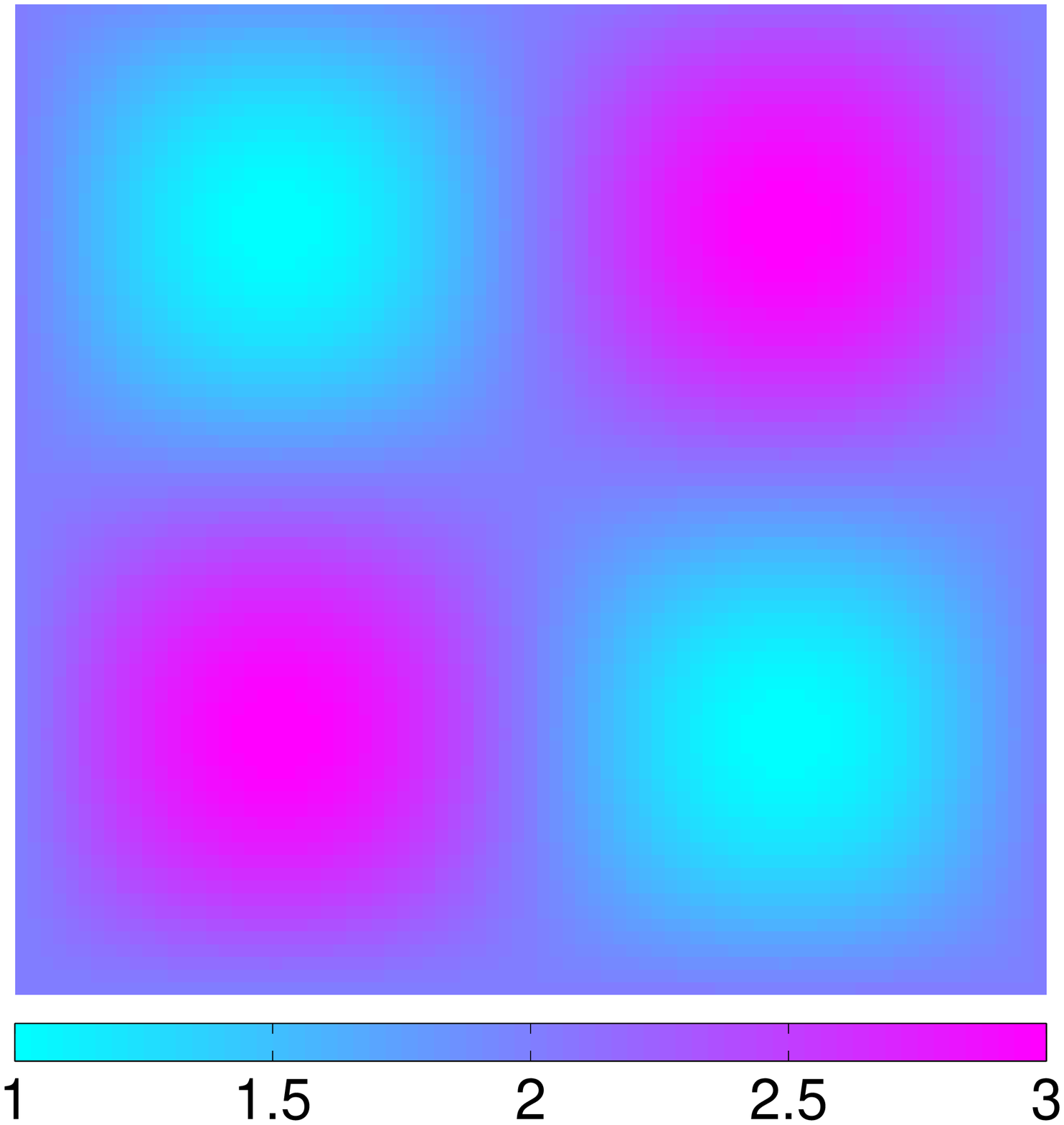}
      \label{ex1txi}
      } 
  \subfigure[$\xi$ ($\alpha=0\%$)]{    
     \includegraphics[width=37mm,height=35mm]{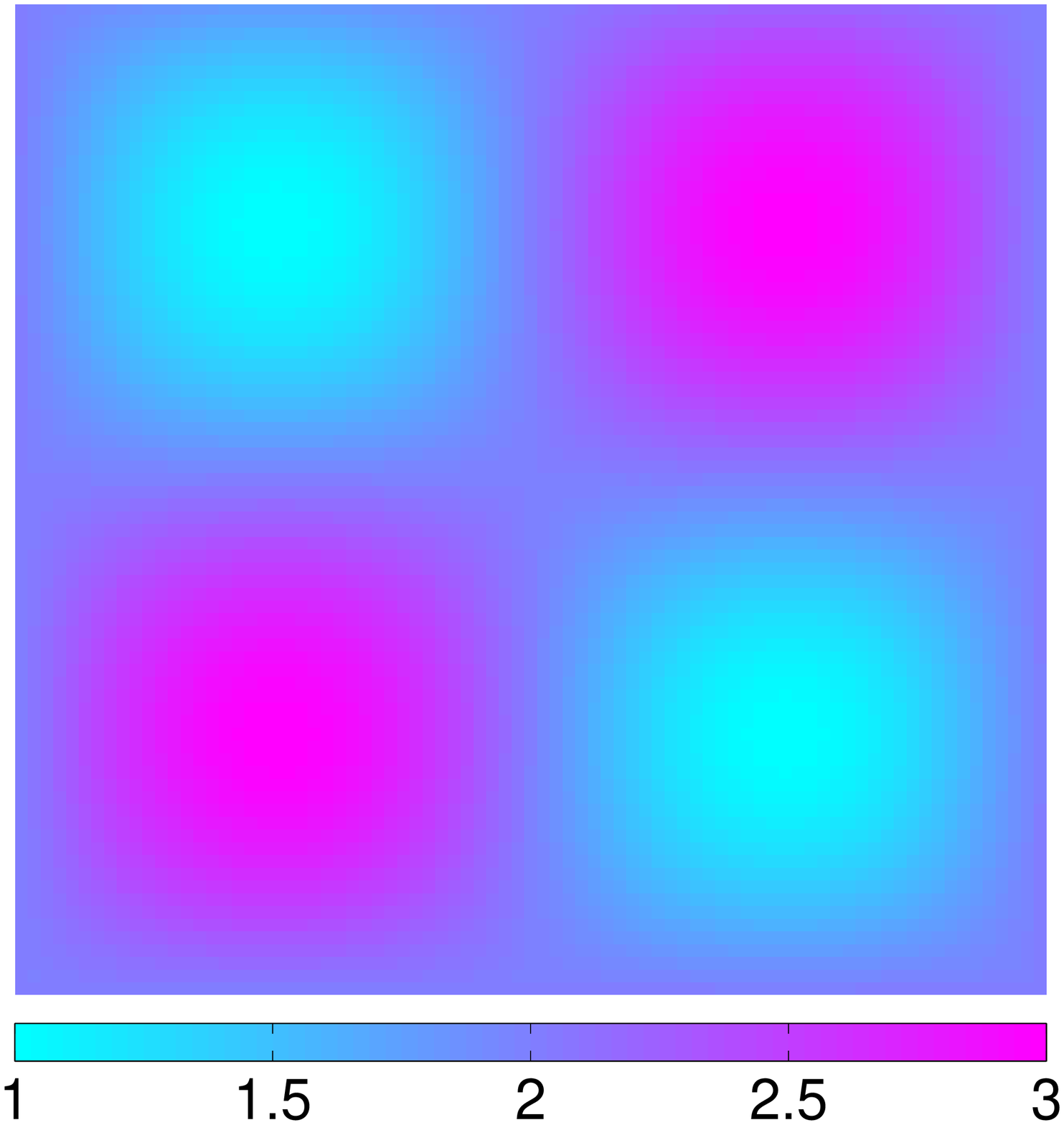}
     \label{ex1cxi}
     }
   \subfigure[$\xi$ ($\alpha=4\%$)]{
    \includegraphics[width=37mm,height=35mm]{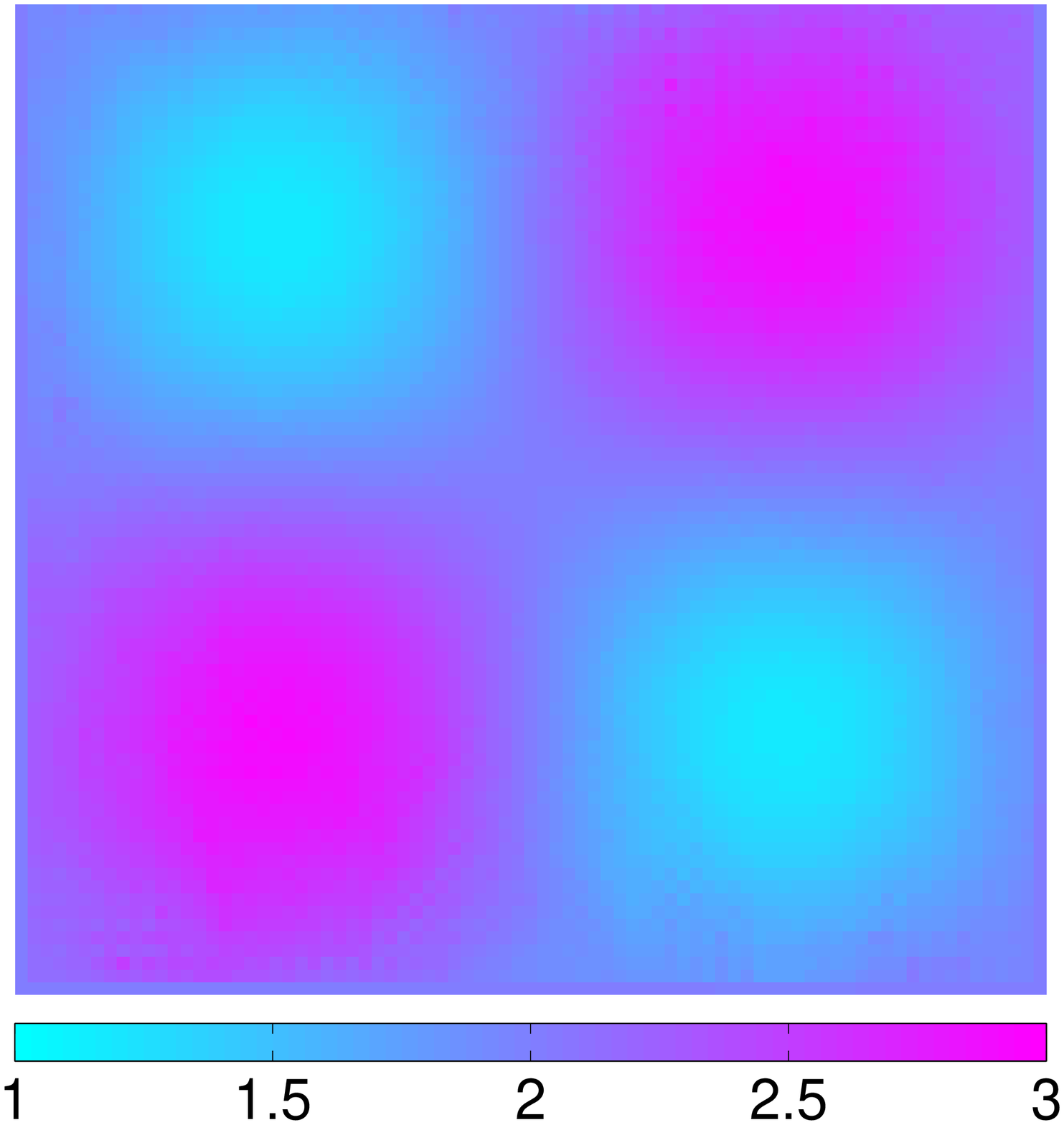}
    \label{ex1nxi}
    }
    \subfigure[$\xi$ at \{$y=0$\}]{
     \includegraphics[trim=10mm 5mm 10mm 0mm,clip,width=35mm,height=38mm]{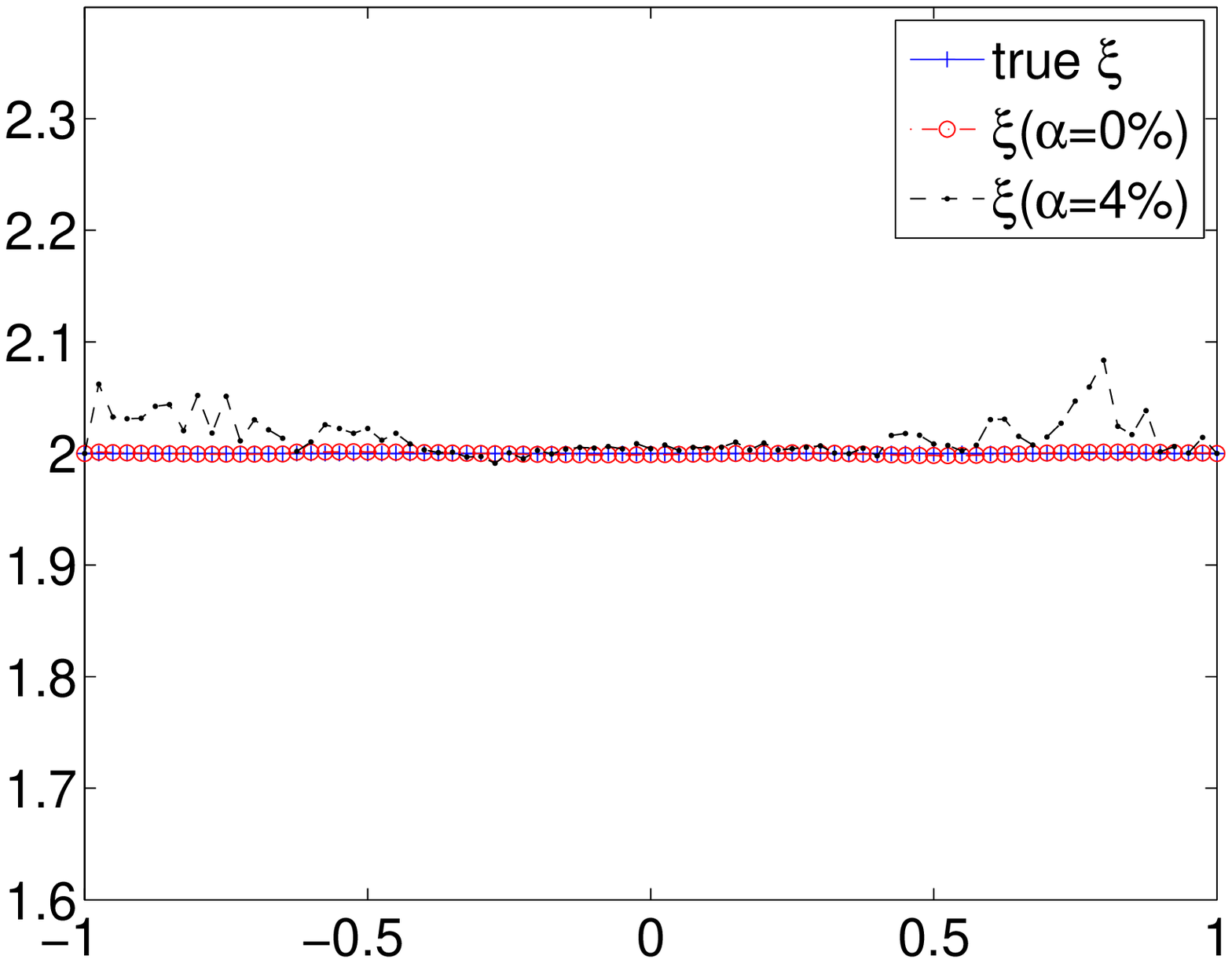} 
     \label{ex1rxi}
     }
     \subfigure[true $\zeta$]{
      \includegraphics[width=37mm,height=35mm]{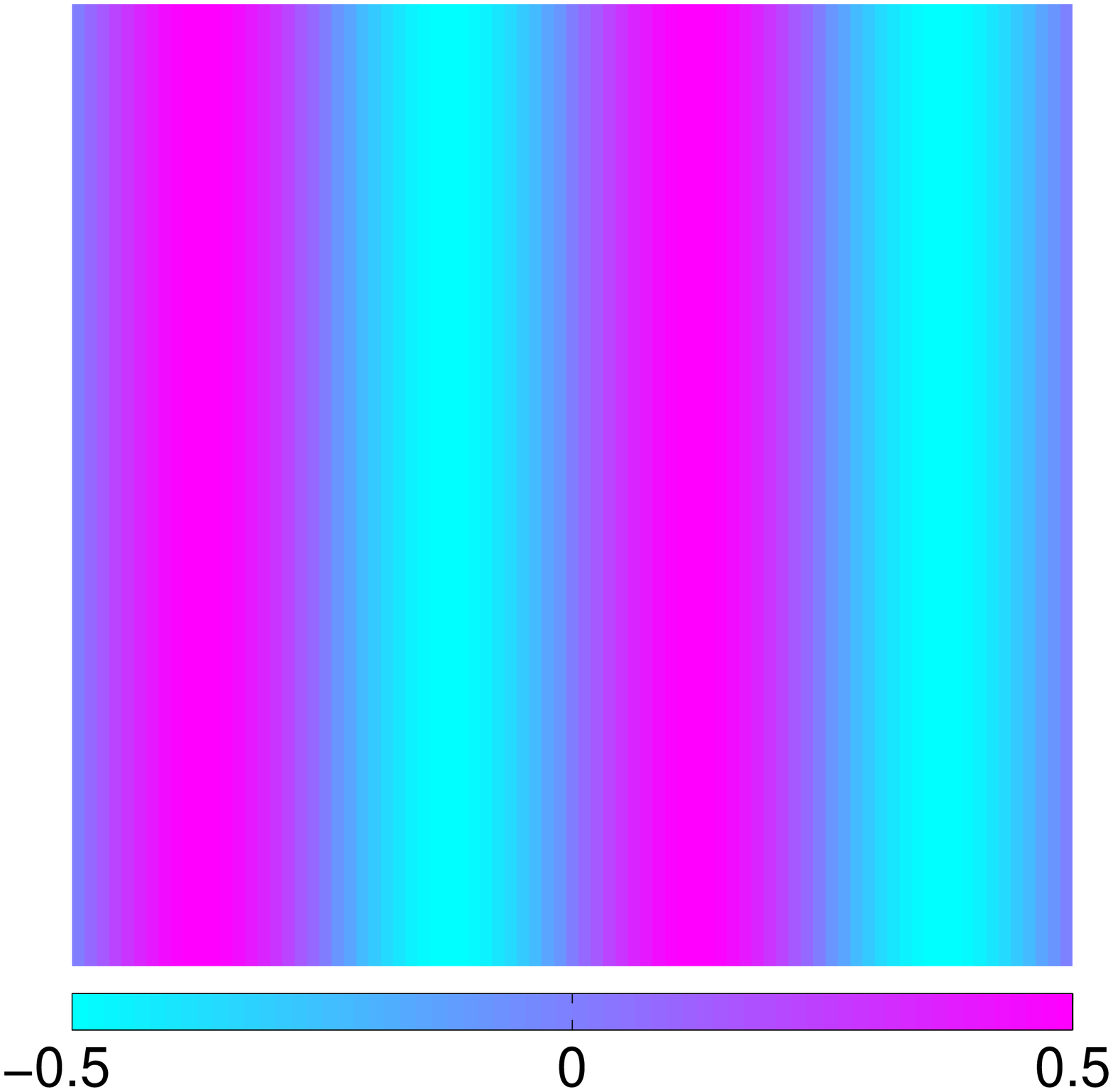}
      \label{ex1ttau}
      } 
    \subfigure[$\zeta$ ($\alpha=0\%$)]{ 
     \includegraphics[width=37mm,height=35mm]{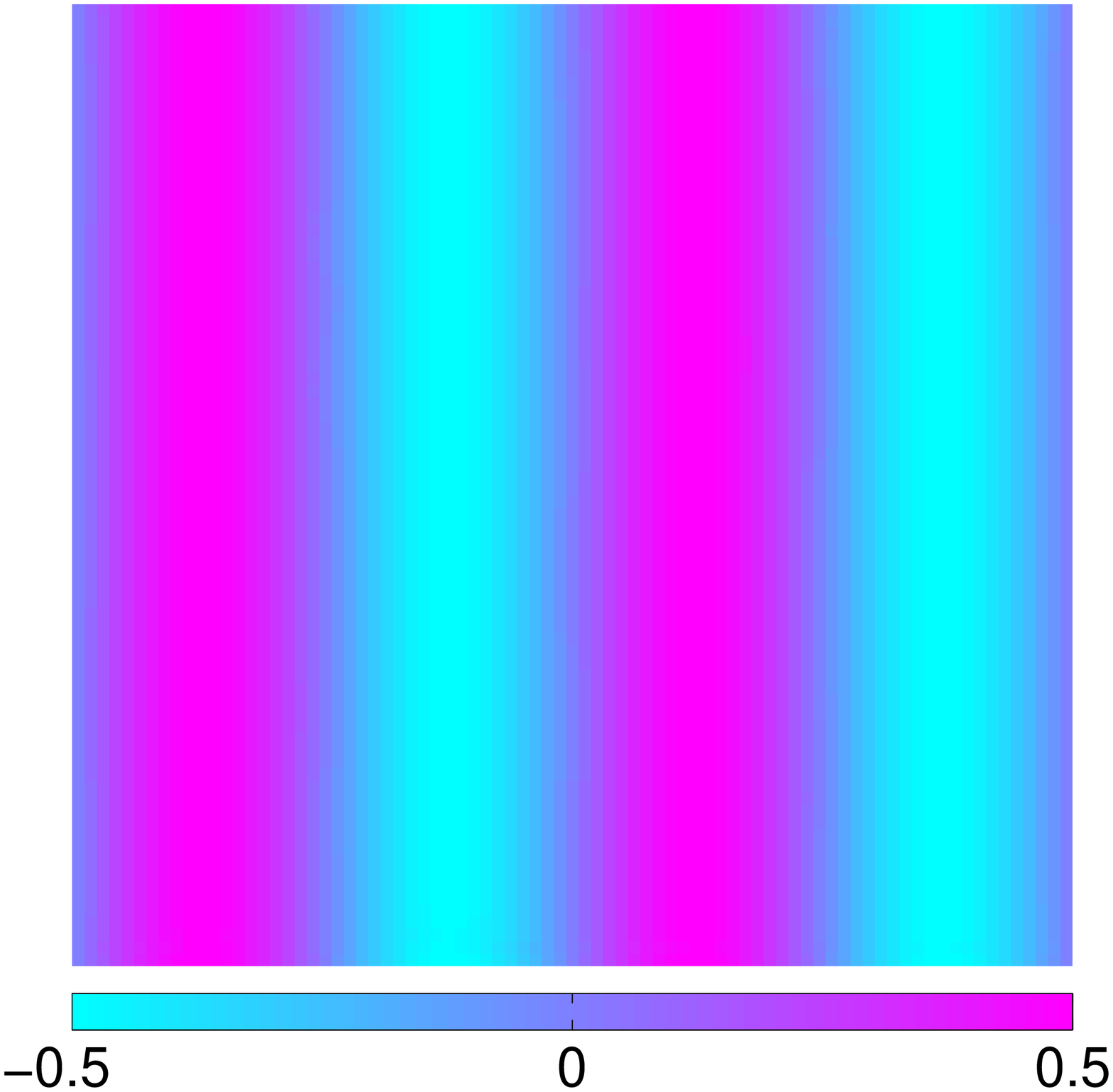}
     \label{ex1ctau}
     }
   \subfigure[$\zeta$ ($\alpha=4\%$)]{ 
    \includegraphics[width=37mm,height=35mm]{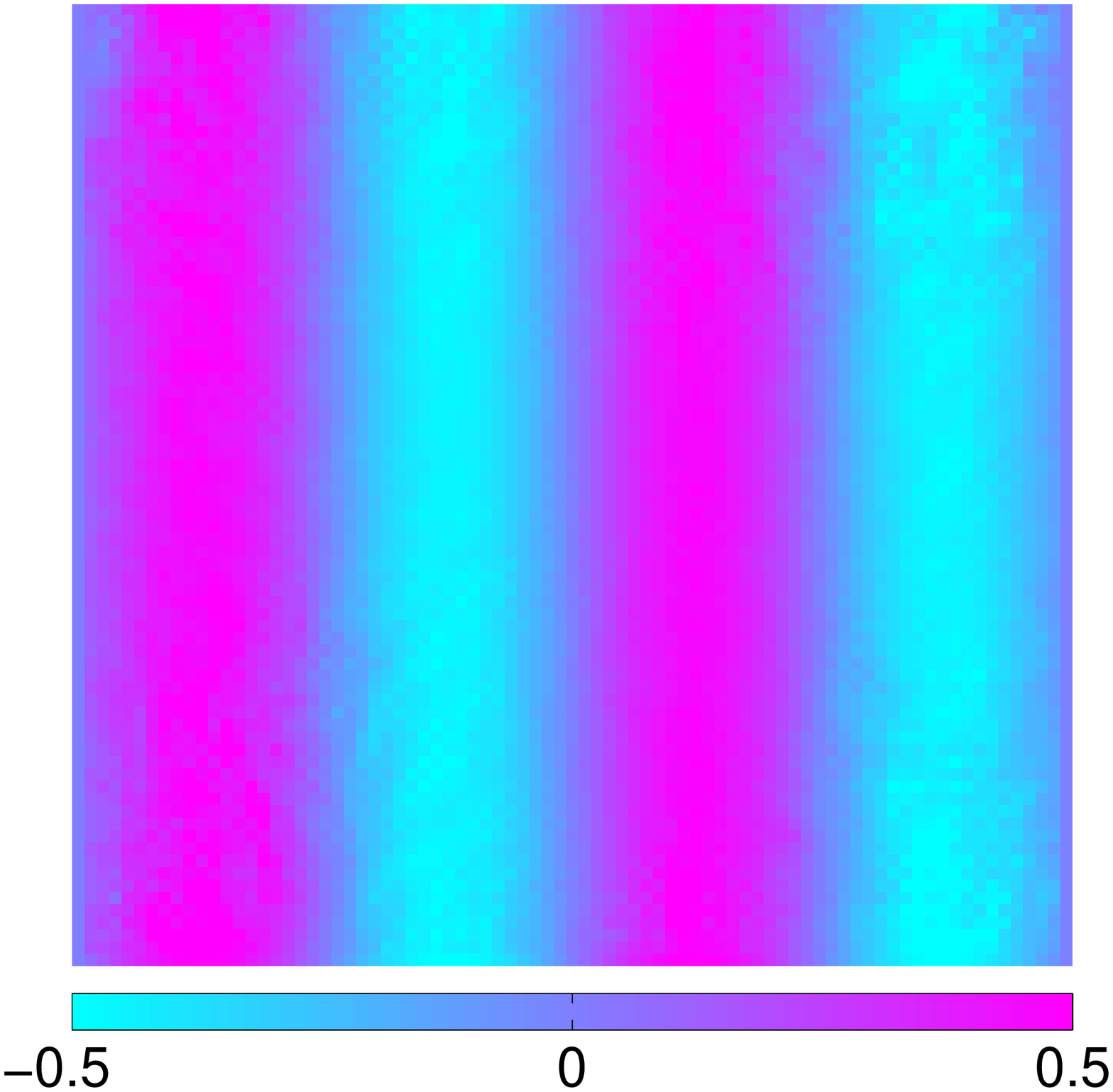}
    \label{ex1ntau}
    }
     \subfigure[$\zeta$ at \{$y=0$\}]{
     \includegraphics[trim=10mm 5mm 10mm 0mm,clip,width=35mm,height=38mm]{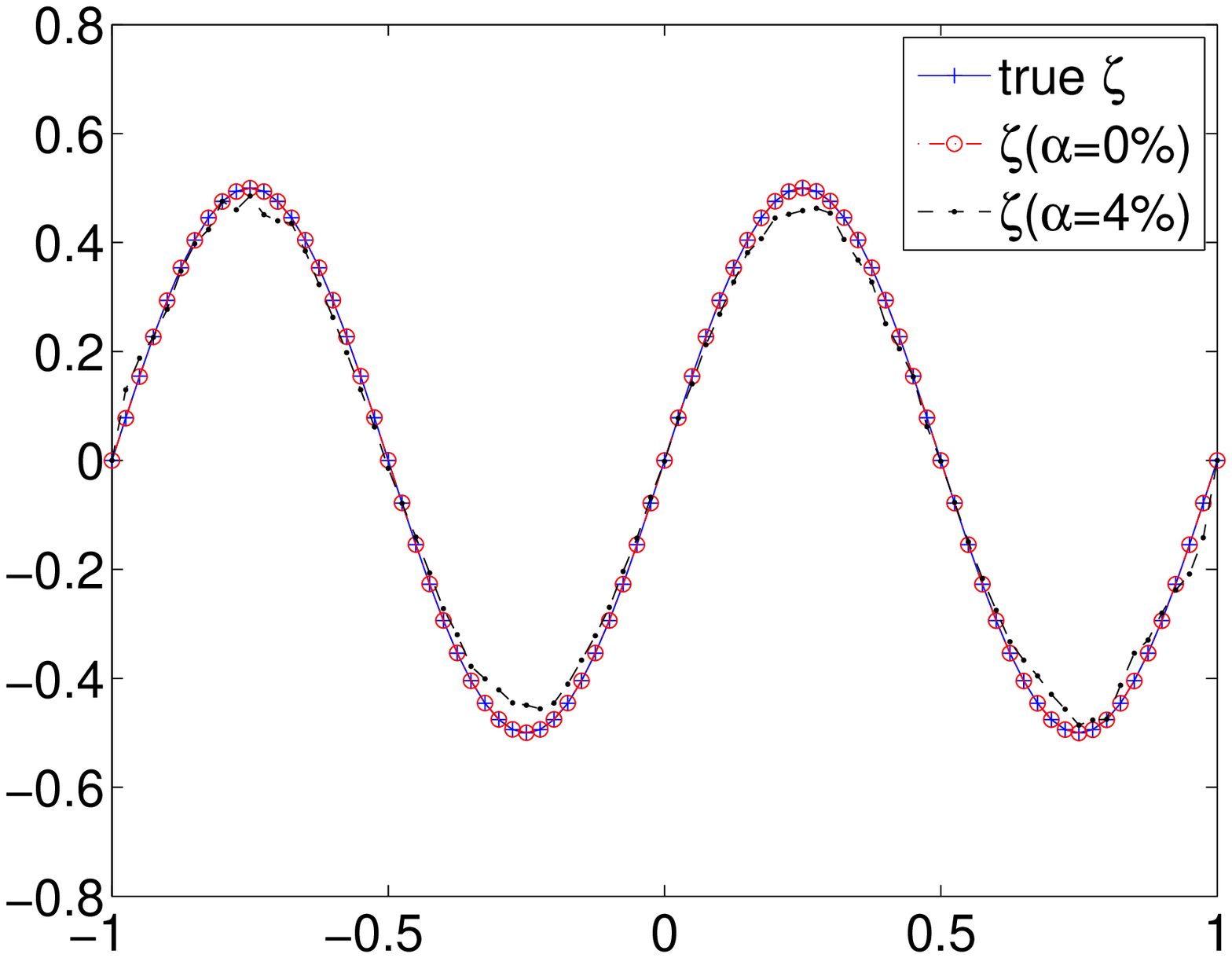} 
     \label{ex1rtau}
     }
     \subfigure[true $|\gamma|^{\frac{1}{2}}$]{
      \includegraphics[width=37mm,height=35mm]{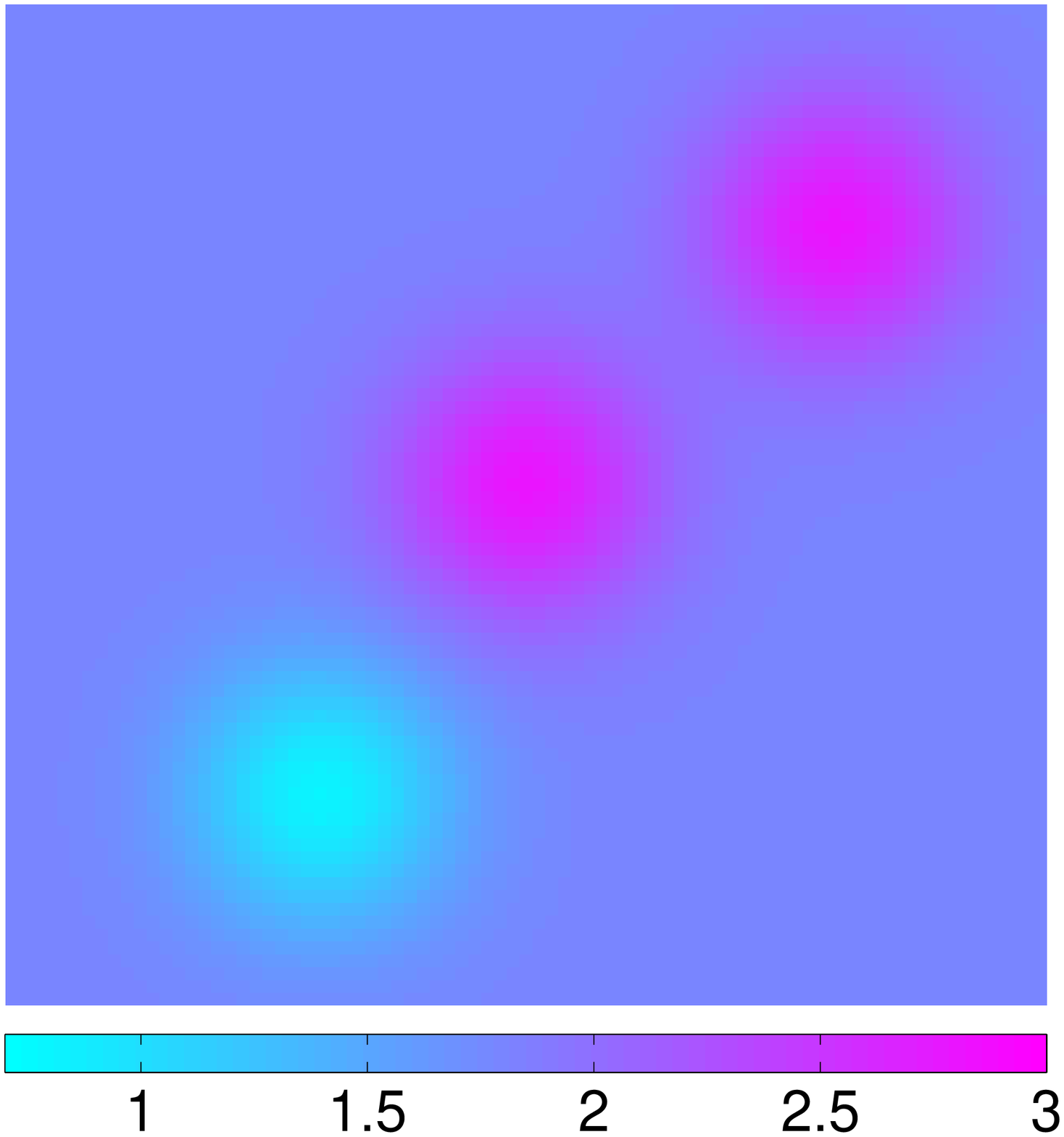}
      \label{ex1tbeta}
      } 
   \subfigure[$|\gamma|^{\frac{1}{2}}$ $(\alpha=0\%)$]{  
     \includegraphics[width=37mm,height=35mm]{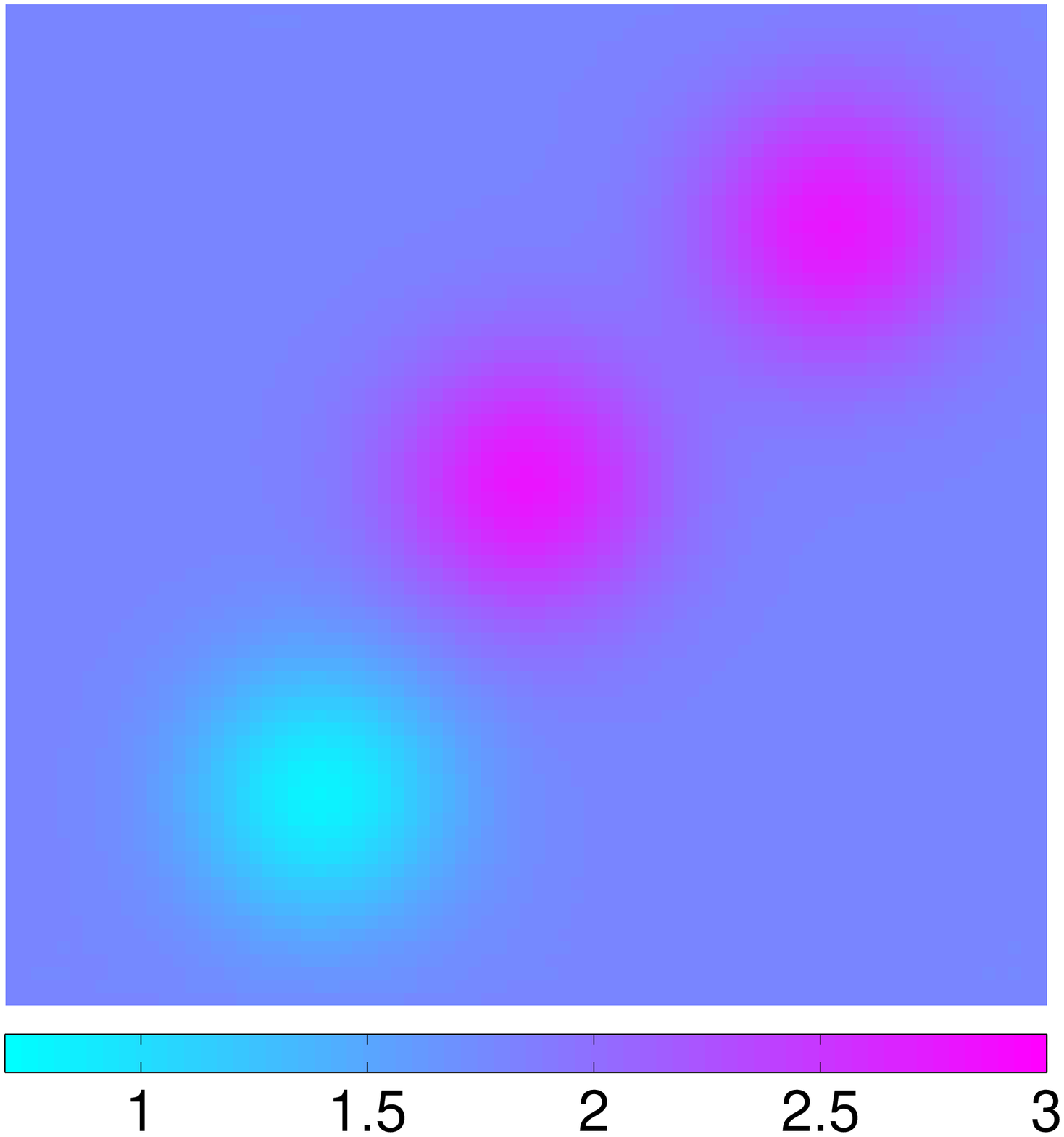}
     \label{ex1cbeta}
     }
   \subfigure[$|\gamma|^{\frac{1}{2}}$ ($\alpha=4\%$)]{ 
    \includegraphics[width=37mm,height=35mm]{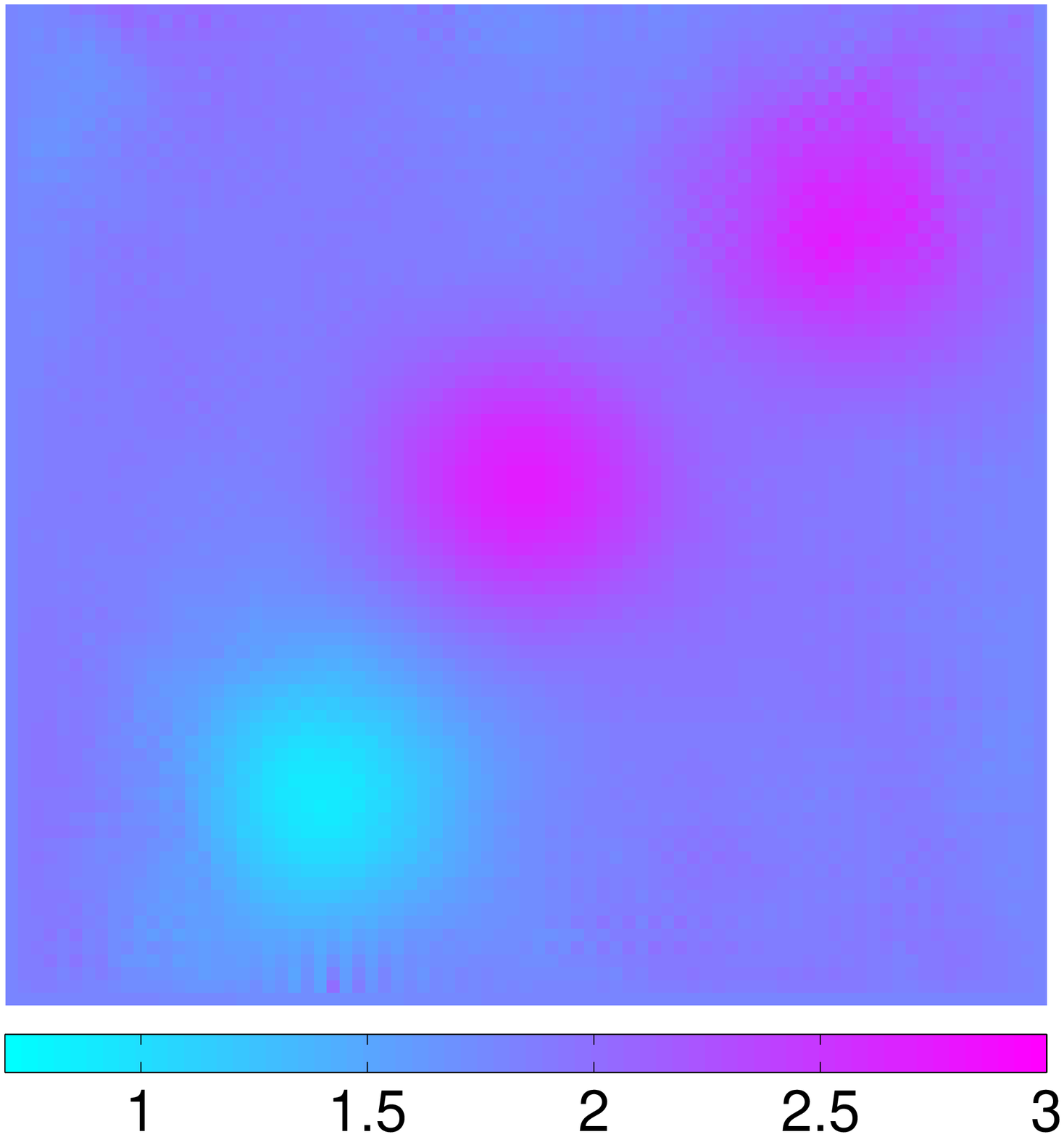}
    \label{ex1nbeta}
    }
   \subfigure[$|\gamma|^{\frac{1}{2}}$ at $\{y=0\}$]{ 
     \includegraphics[trim=10mm 5mm 10mm 0mm,clip,width=35mm,height=38mm]{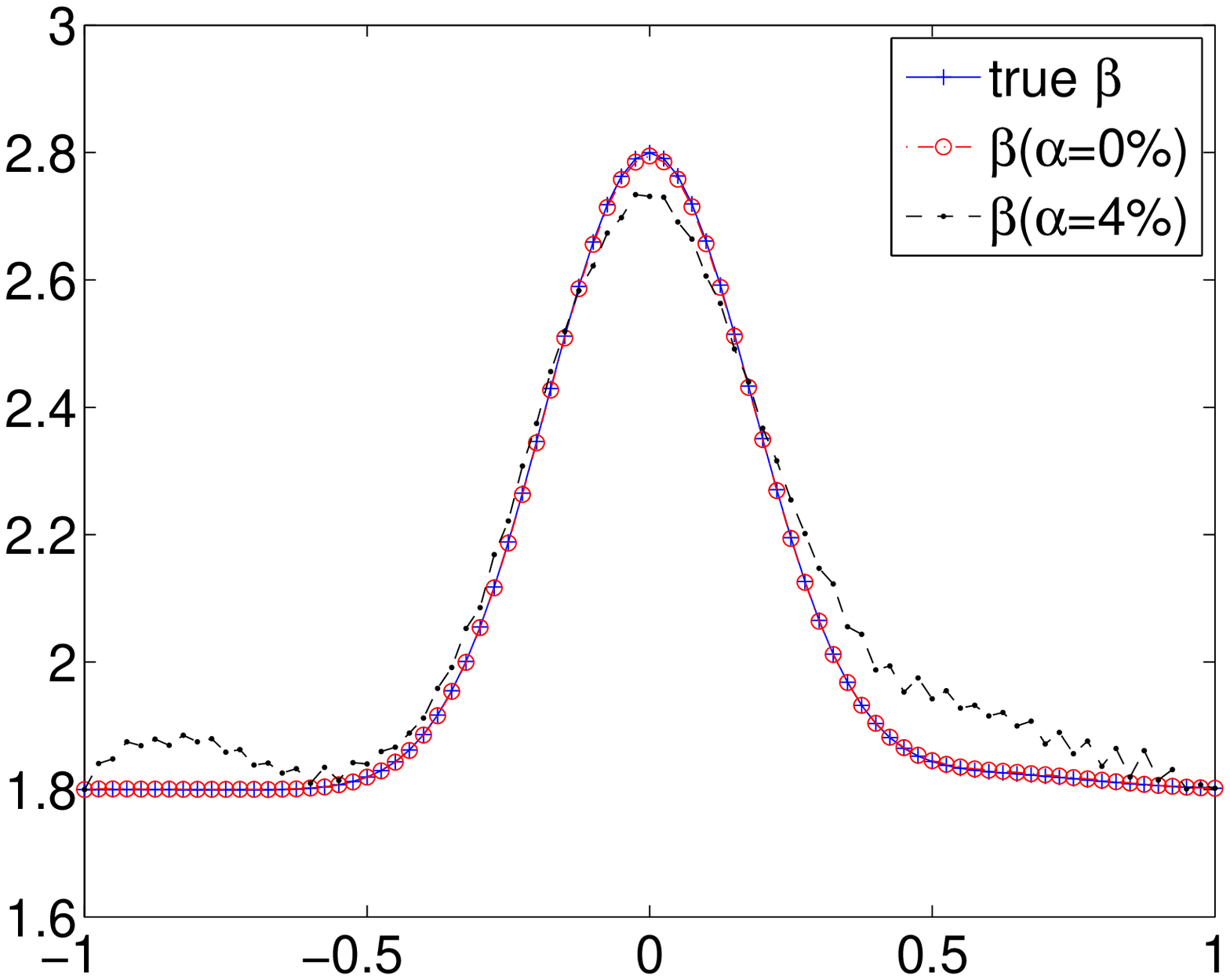} 
     \label{ex1rbeta}
     }
    \caption{Experiment 1. \subref{ex1txi}\&\subref{ex1ttau}\&\subref{ex1tbeta}: true values of $(\xi, \zeta, \beta)$. \subref{ex1cxi}\&\subref{ex1ctau}\&\subref{ex1cbeta}: reconstructions with noiseless data. \subref{ex1nxi}\&\subref{ex1ntau}\&\subref{ex1nbeta}: reconstructions with noisy data($\alpha=4\%$). \subref{ex1rxi}\&\subref{ex1rtau}\&\subref{ex1rbeta}: cross sections along $\{y=0\}$.}
\label{E1}
\end{figure}

\paragraph{Experiment 1.} In the first experiment, we intend to reconstruct the smooth coefficients $\xi$, $\zeta$ and $\beta$ defined in \eqref{3coef} and given by,
\begin{align}
\left\{\begin{array}{lll}
\xi = 2+ \sin(\pi x)\sin(\pi y)\\
 \zeta = 0.5\sin(2\pi x)\\
\beta = 1.8+e^{-15(x^2+y^2)} +e^{-15((x-0.6)^2+(y-0.5)^2)} - e^{-15((x+0.4)^2+(y+0.6)^2)} .
\end{array}\right.
\end{align}
We consider five different illuminations $(g_1,g_2,g_3,g_4,g_5)$ that are defined as follows,
\begin{align}\label{5 bc}
(g_1,g_2,g_3,g_4,g_5)(\mathbf{x}) = (x+y,y+0.1y^2,3x^2+2y^2,x^2-0.5y^2,xy) \quad \mathbf{x}\in \partial[-1,1]^2
\end{align}
where $g_1,g_2$ are used generating the solutions satisfying Lemma \ref{cond cgo}.$A$. We performed two sets of reconstructions using clean and noisy synthetic data respectively. The $l_2$-regularization procedure is used in this simulation. For the noisy data, the noise level $\alpha=4\%$. The results of the numerical experiment are shown in Figure \ref{E1}. The relative $L^2$ errors in the reconstructions are $\mathcal{E}^C_{\xi}=0.1\%$, $\mathcal{E}^N_{\xi}=4.0\%$, $\mathcal{E}^C_{\zeta}=0.6\%$, $\mathcal{E}^N_{\zeta}=11.8\%$, $\mathcal{E}^C_{\beta}=0.2\%$ and $\mathcal{E}^N_{\beta}=3.7\%$. 

\begin{figure}[htp]
  \centering
  \subfigure[$|\gamma|^{\frac{1}{2}}$ ($\alpha=20\%$)]{ 
	    \includegraphics[width=37mm,height=35mm]{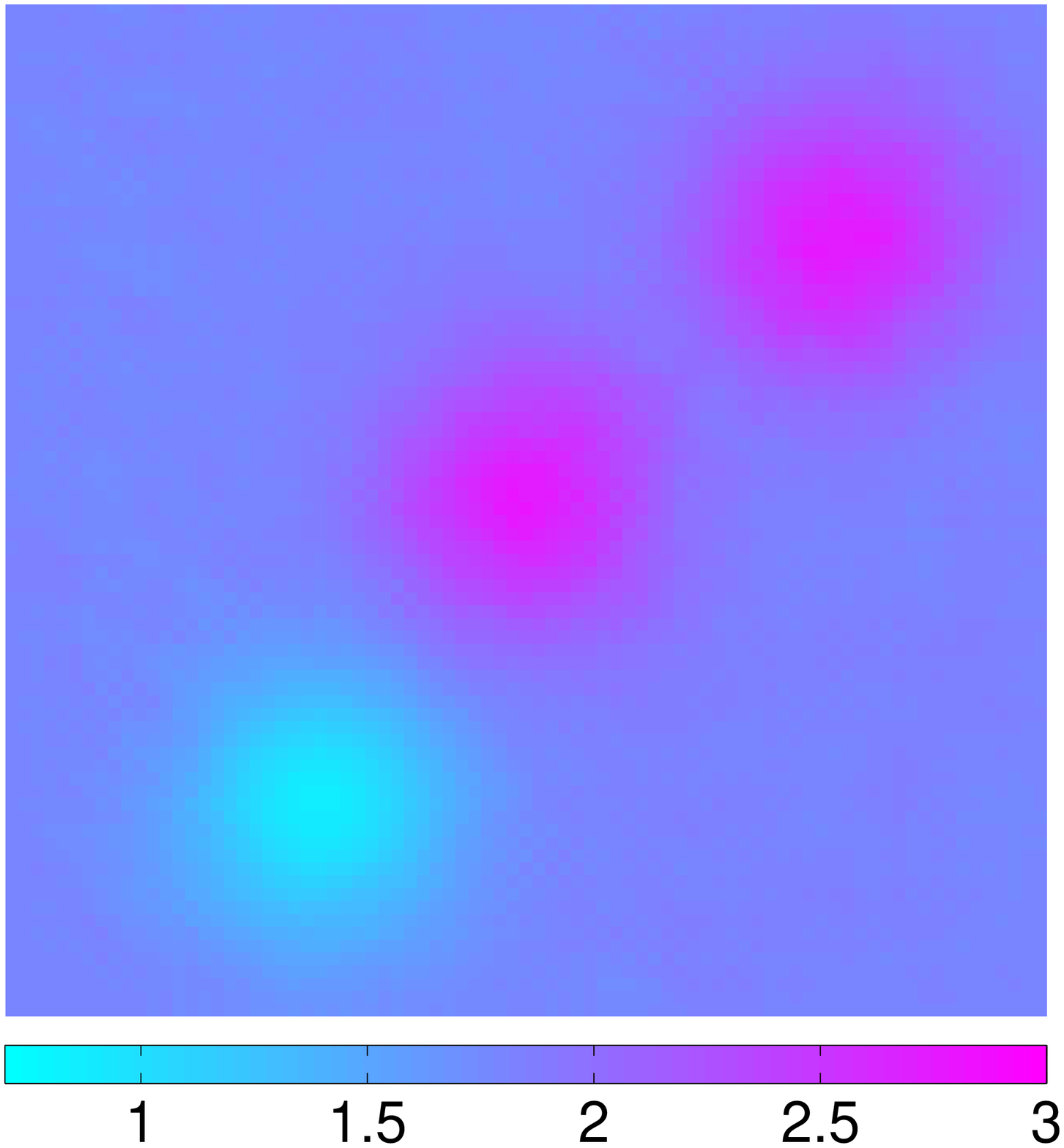} 
      \label{beta 20perc}
  }
  \qquad
	\subfigure[$|\gamma|^{\frac{1}{2}}$ at $\{y=0\}$]{ 
		  \includegraphics[trim=10mm 5mm 10mm 0mm,clip,width=35mm,height=38mm]{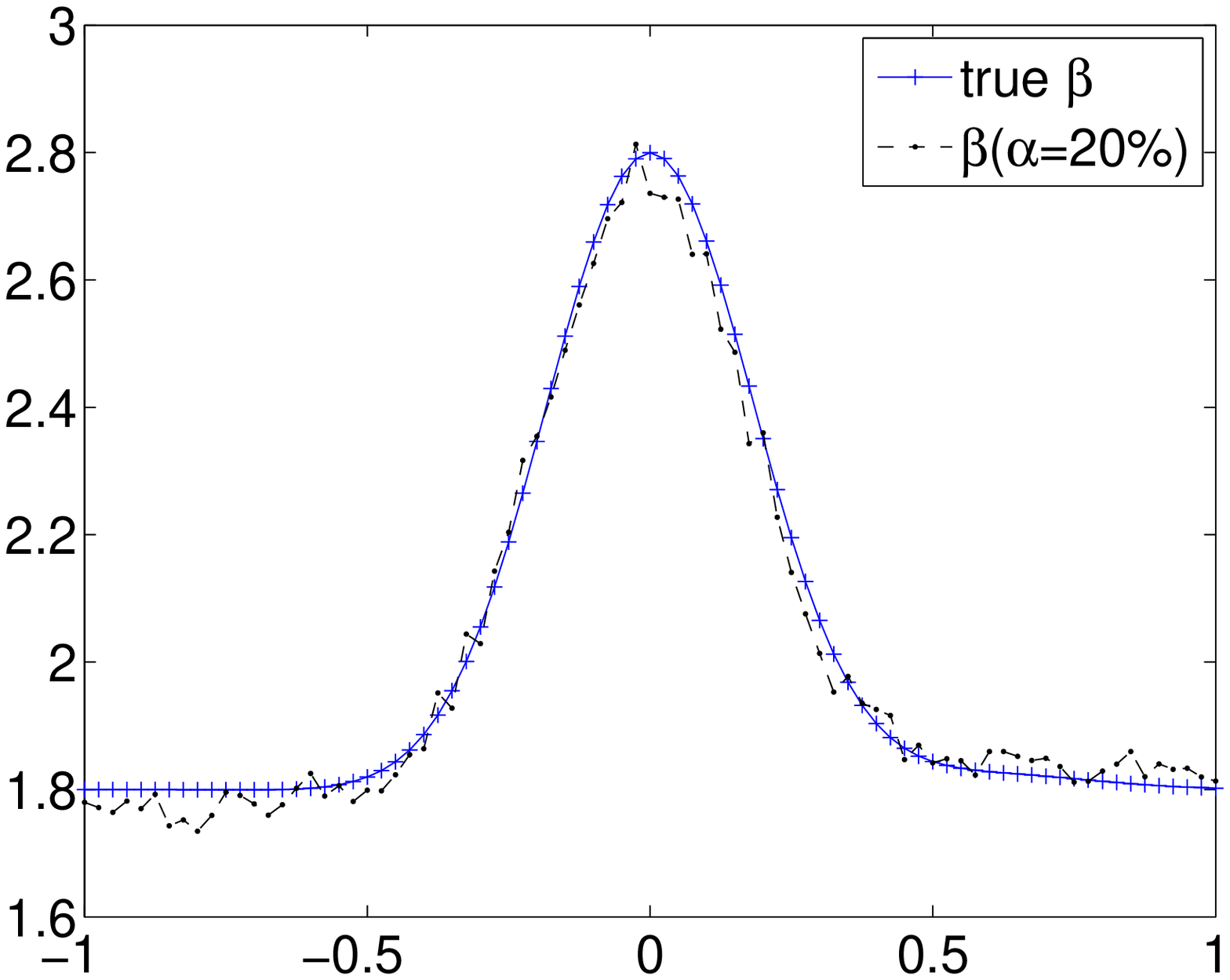}   
      \label{cross beta 20perc}
  } 
  \caption{Reconstruction of $\beta$ with true anisotropy. \subref{beta 20perc}: reconstructed $\beta$ using true anisotropy and noisy data($\alpha=20\%$). \subref{cross beta 20perc}: cross-section along $y=0$. }
  \label{E1anis} 
\end{figure}

\bigbreak
\noindent{\emph{Reconstruction of $\beta$ with (known) true anisotropic part $\tilde\gamma$}.} We now use the true $\xi$ and $\zeta$ to reconstruct $\beta$ with noisy data ($\alpha=20\%$). Figure \ref{E1anis} displays the numerical results. The reconstruction is quite robust to noise when the anisotropy is known: the $L^2$ relative error is $1.6\%$. Comparing Fig.\ref{ex1nbeta}\&\subref{ex1rbeta} with Fig.\ref{beta 20perc}\&\subref{cross beta 20perc}, it is clear that the reconstruction of the isotropy $\beta$ is more stable than that of the anisotropy $\tilde\gamma$. This is consistent with the better stability estimates obtained in Theorem \ref{stability}.

\begin{figure}
  \centering
 \subfigure[true $\xi$]{
      \includegraphics[width=37mm,height=35mm]{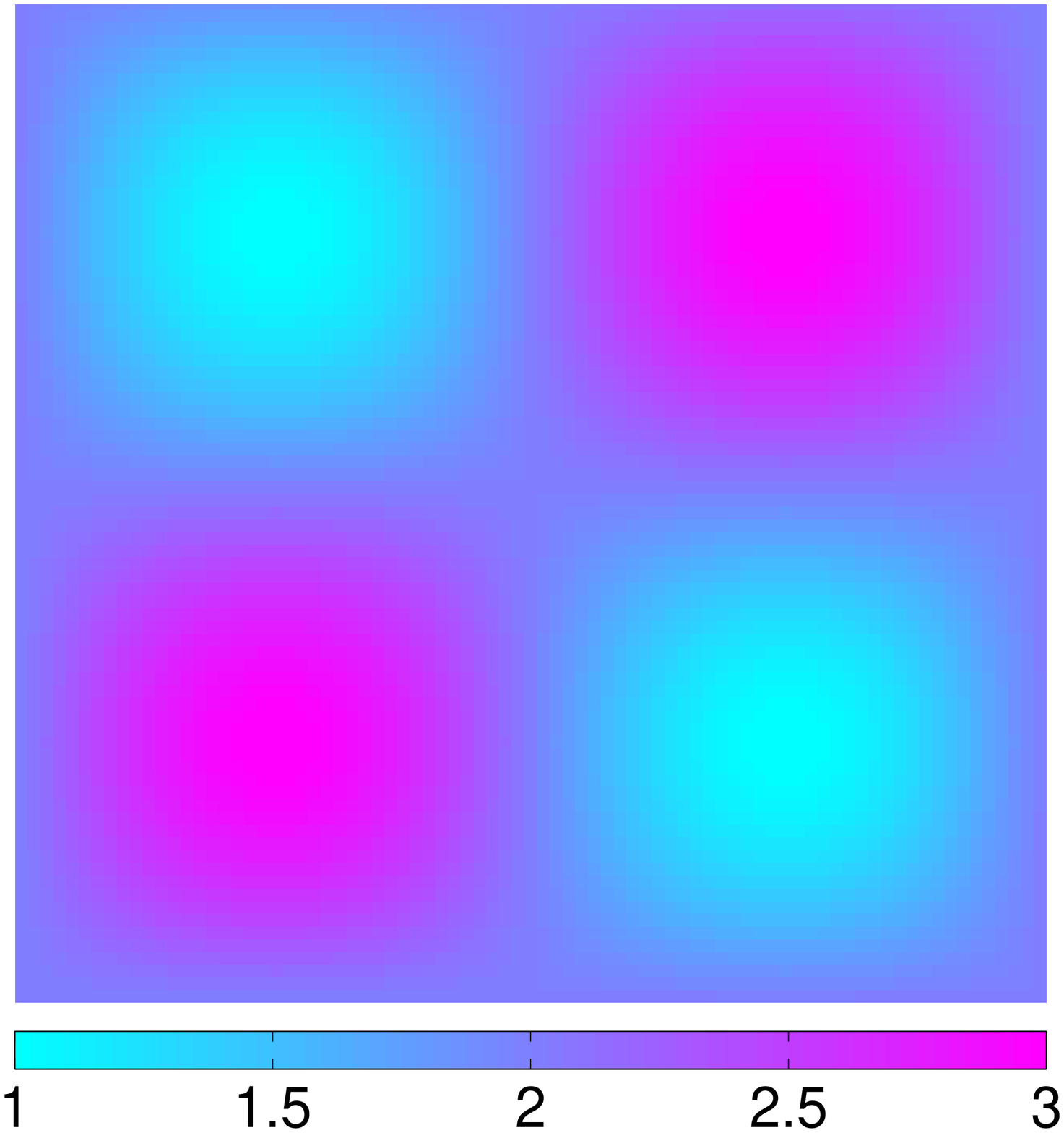}
      \label{ex2txi}
      } 
  \subfigure[$\xi$ ($\alpha=0\%$)]{    
     \includegraphics[width=37mm,height=35mm]{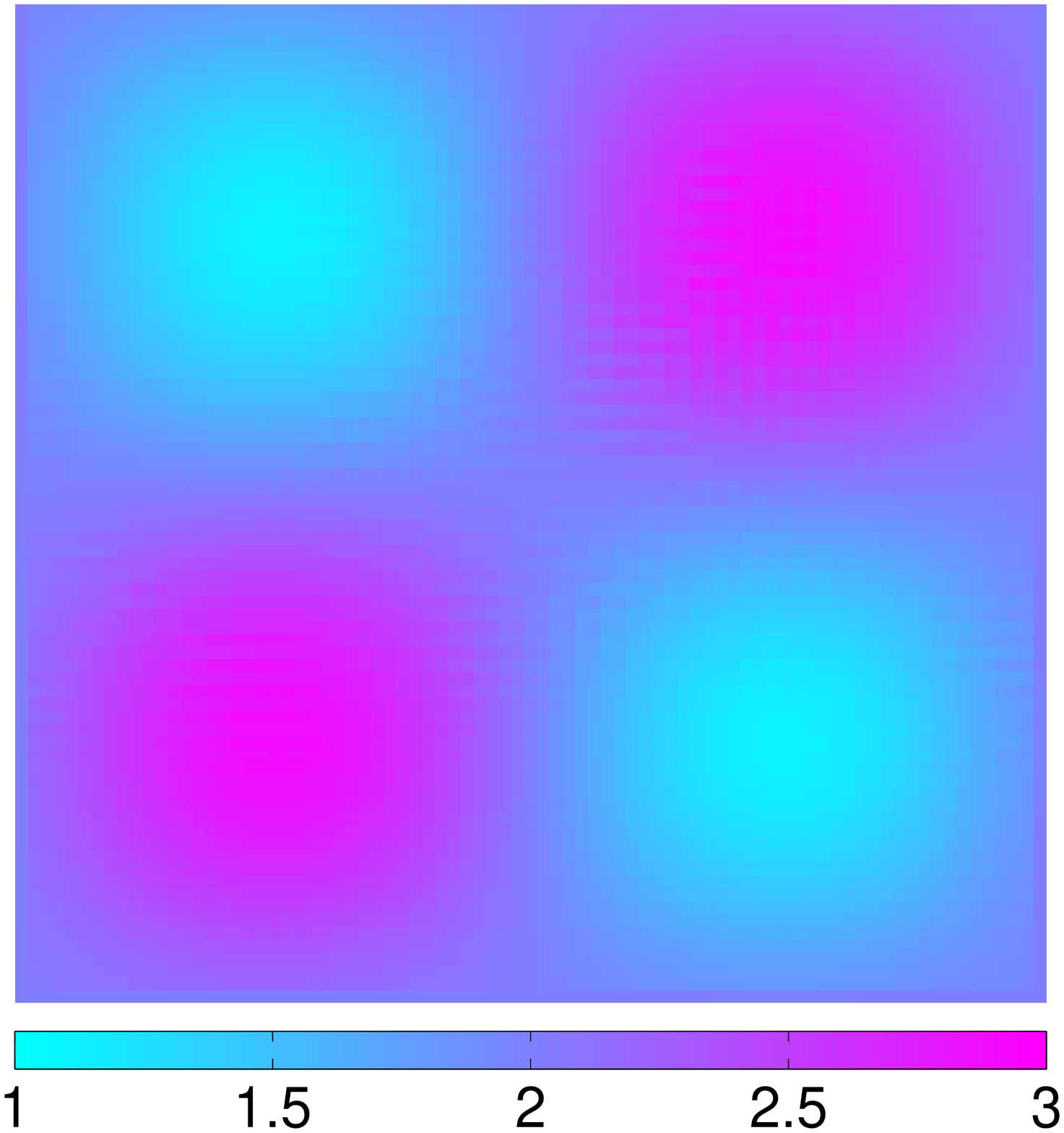}
     \label{ex2cxi}
     }
   \subfigure[$\xi$ ($\alpha=4\%$)]{
    \includegraphics[width=37mm,height=35mm]{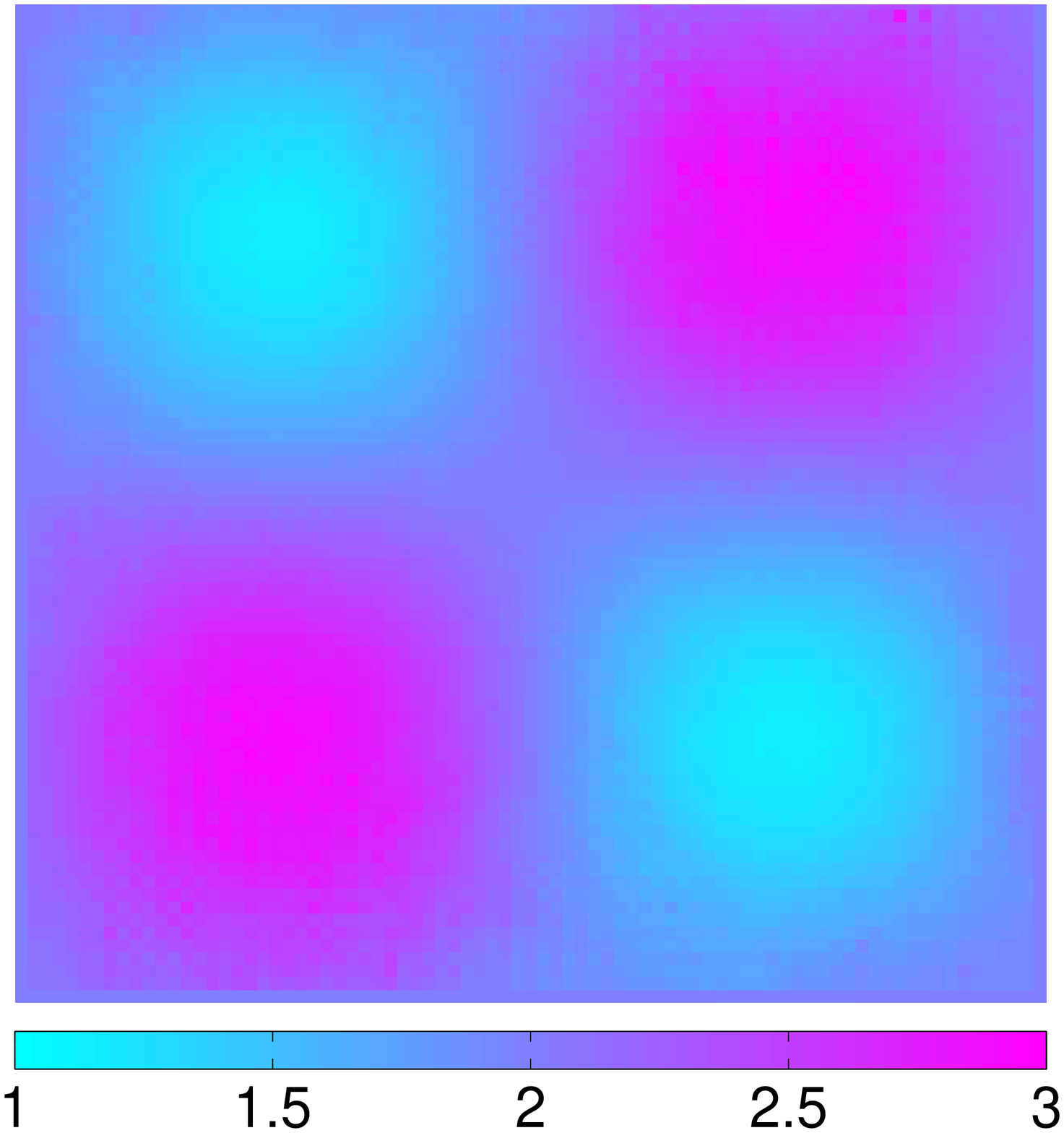}
    \label{ex2nxi}
    }
    \subfigure[$\xi$ at \{$y=0$\}]{
     \includegraphics[trim=10mm 5mm 10mm 0mm,clip,width=35mm,height=38mm]{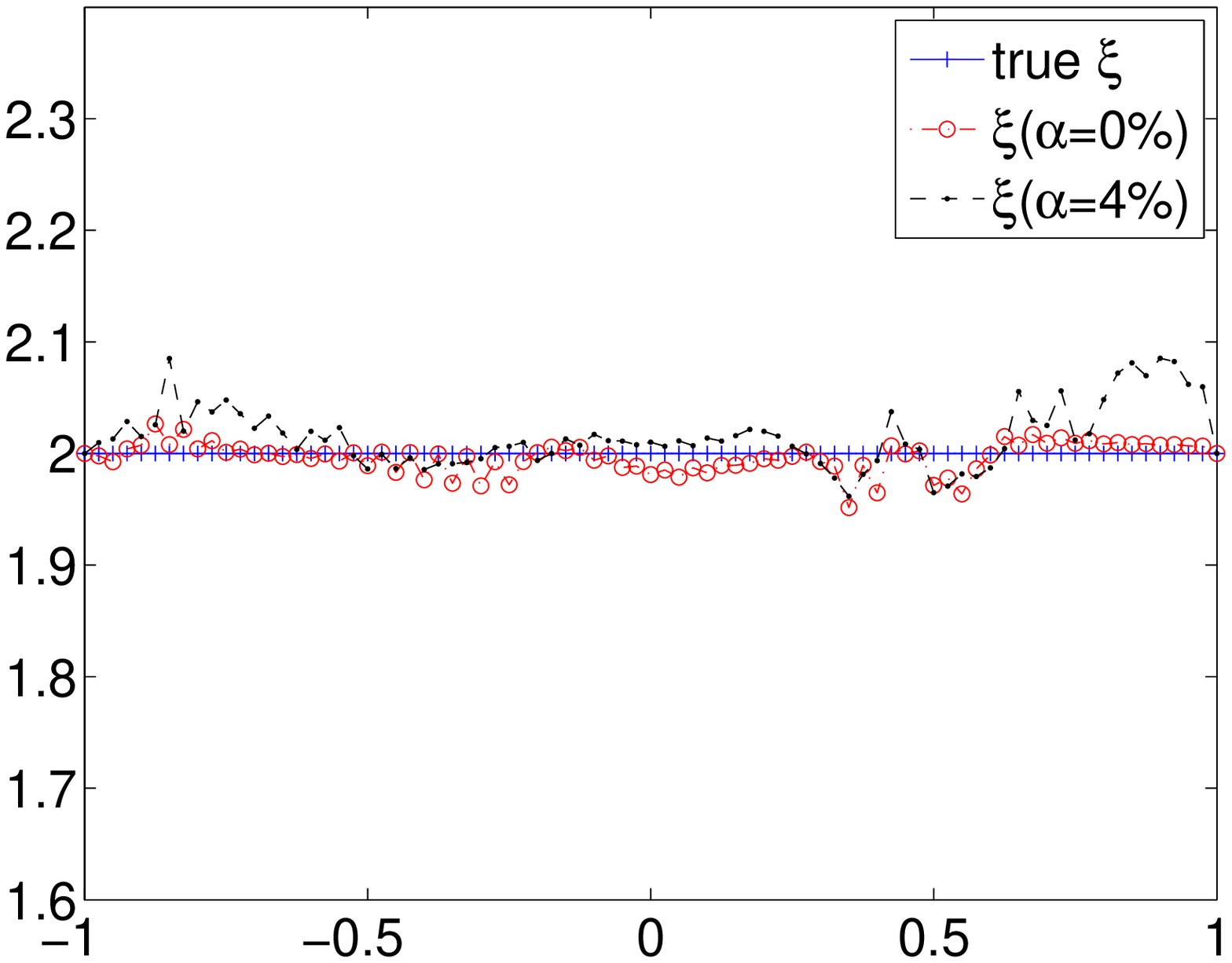} 
     \label{ex2rxi}
     }
     
     \subfigure[true $\zeta$]{
      \includegraphics[width=37mm,height=35mm]{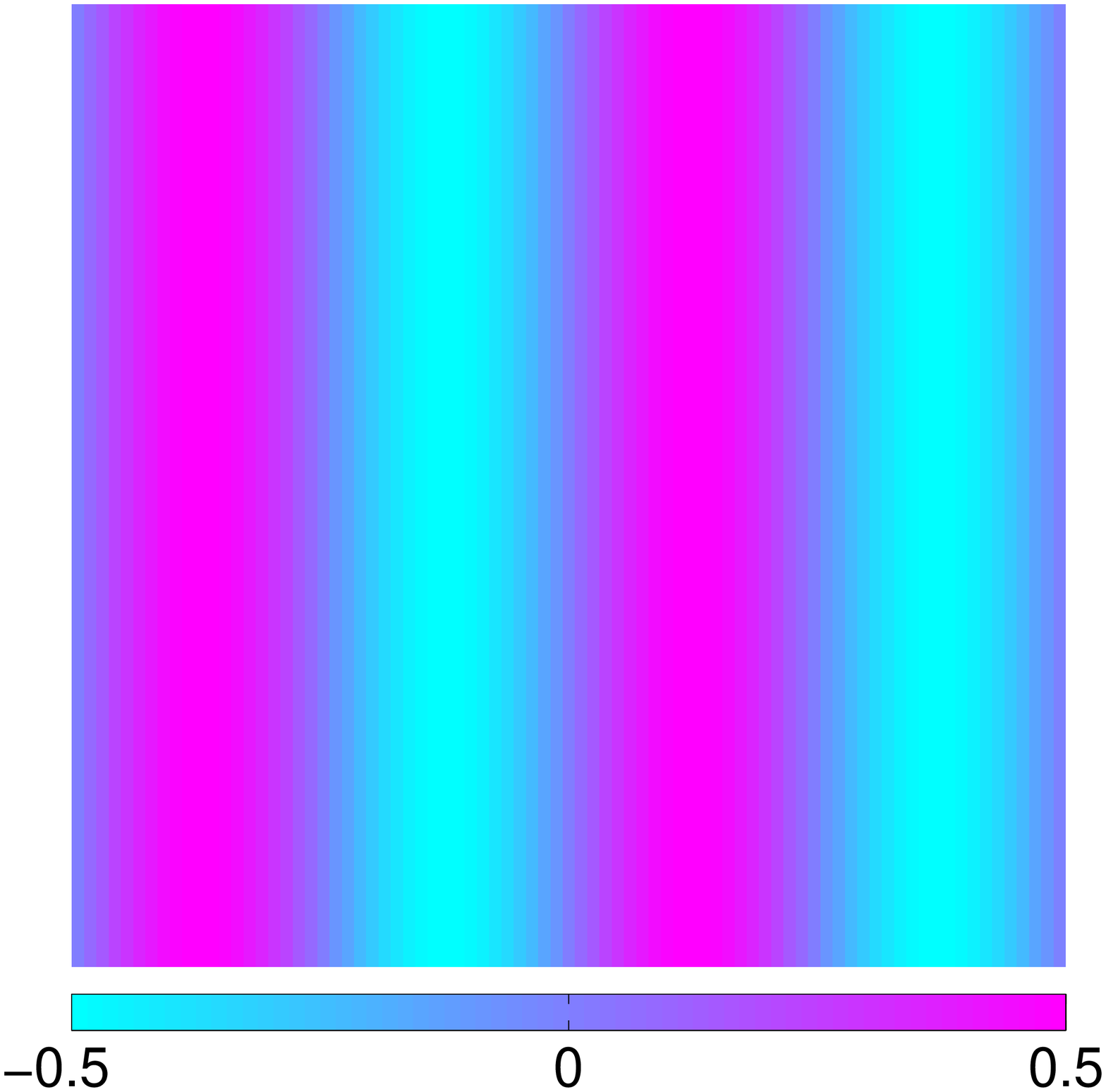}
      \label{ex2ttau}
      } 
 \subfigure[$\zeta$ ($\alpha=0\%$)]{ 
     \includegraphics[width=37mm,height=35mm]{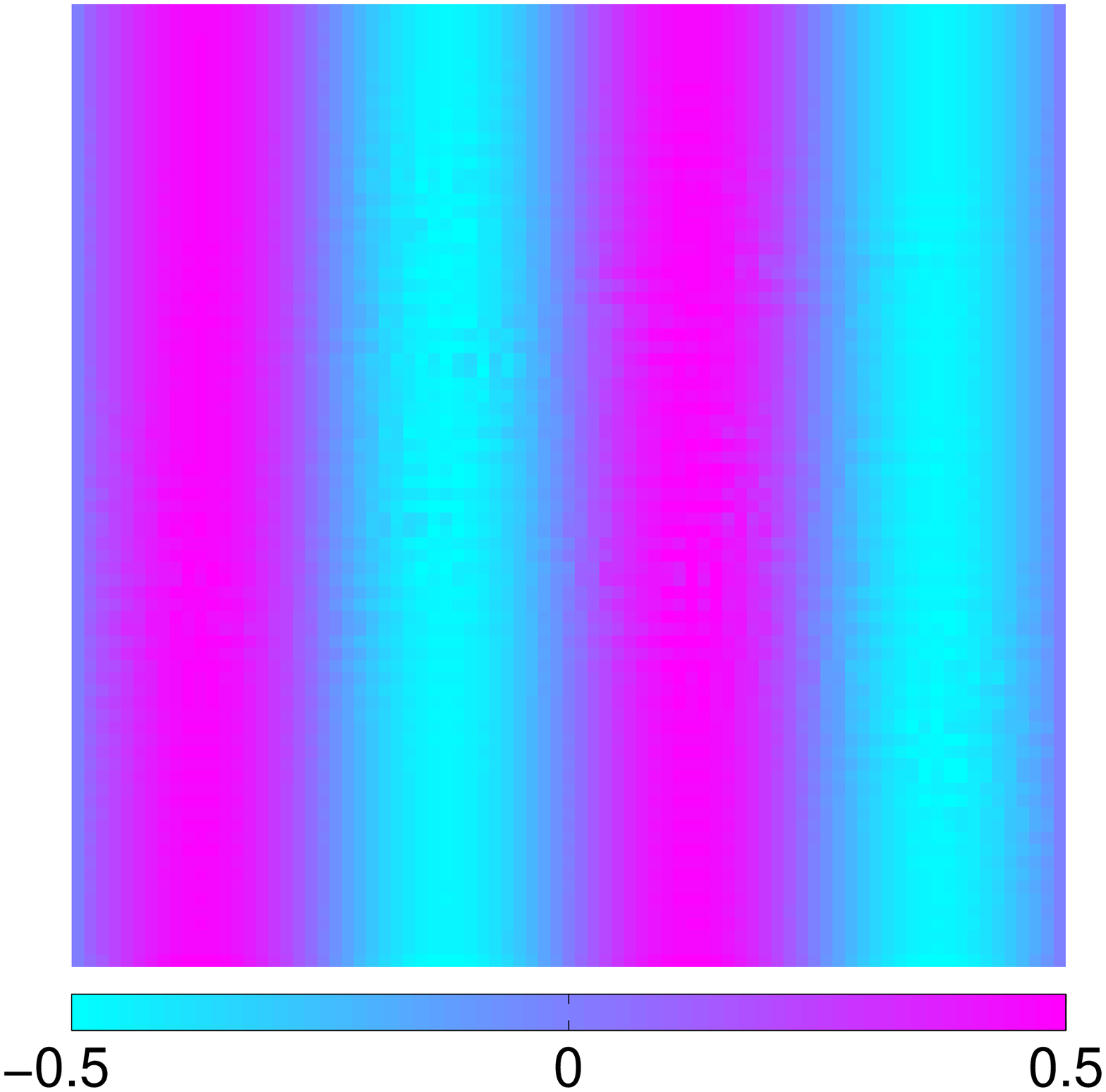}
     \label{ex2ctau}
     }
   \subfigure[$\zeta$ ($\alpha=4\%$)]{ 
    \includegraphics[width=37mm,height=35mm]{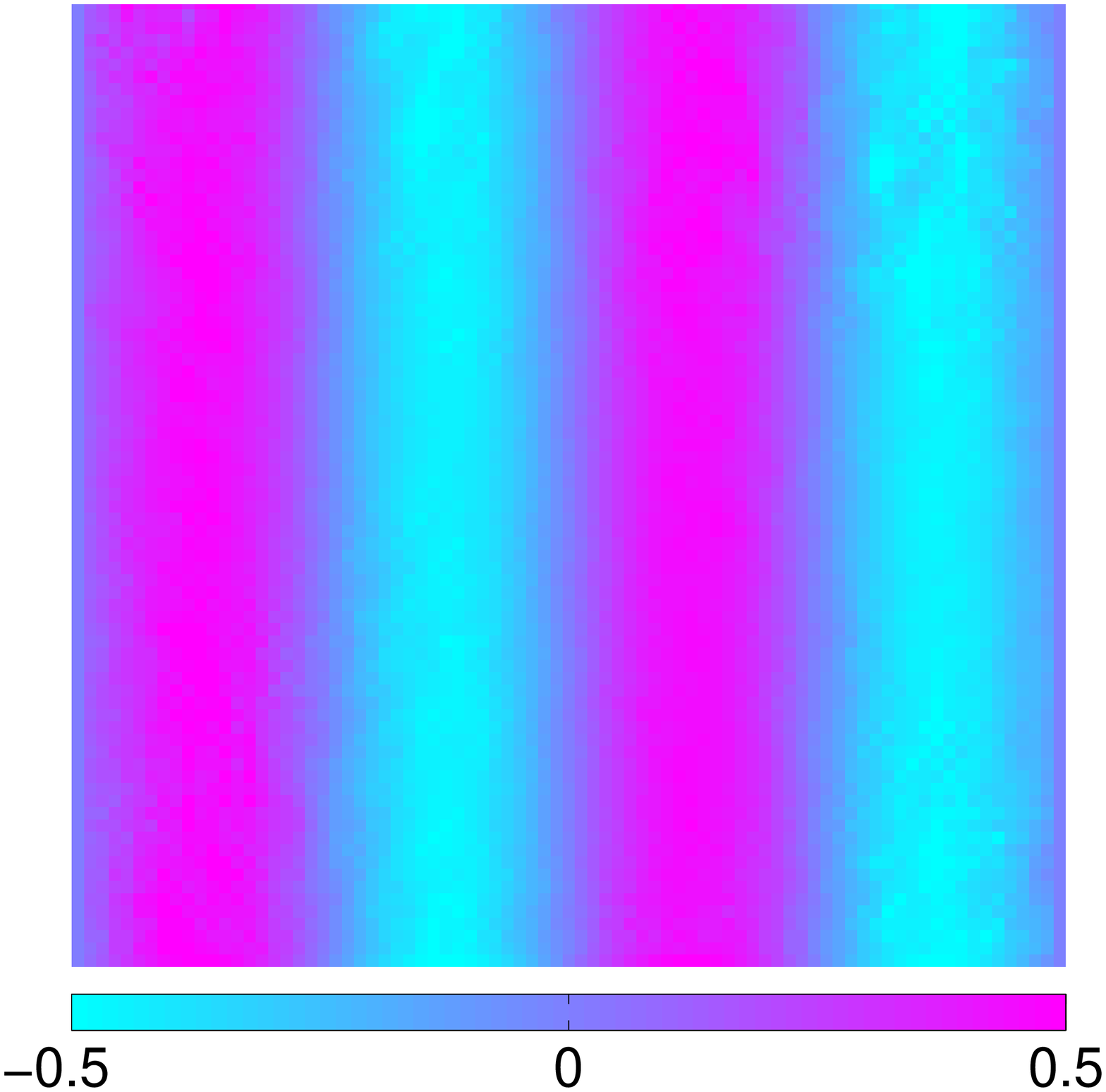}
    \label{ex2ntau}
    }
     \subfigure[$\zeta$ at \{$y=0$\}]{
     \includegraphics[trim=10mm 5mm 10mm 0mm,clip,width=35mm,height=38mm]{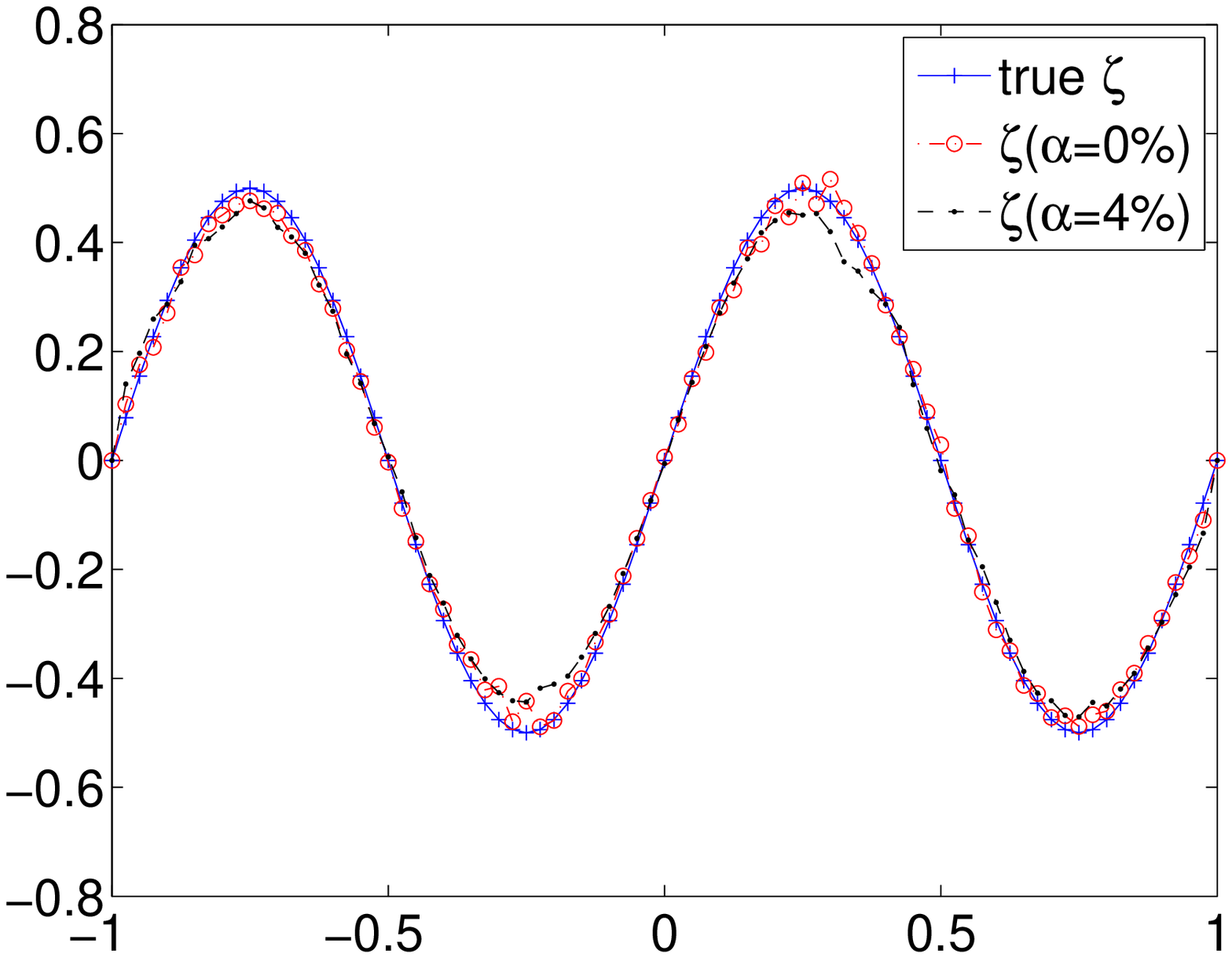} 
     \label{ex2rtau}
     }
       
     \subfigure[true $|\gamma|^{\frac{1}{2}}$]{
      \includegraphics[width=37mm,height=35mm]{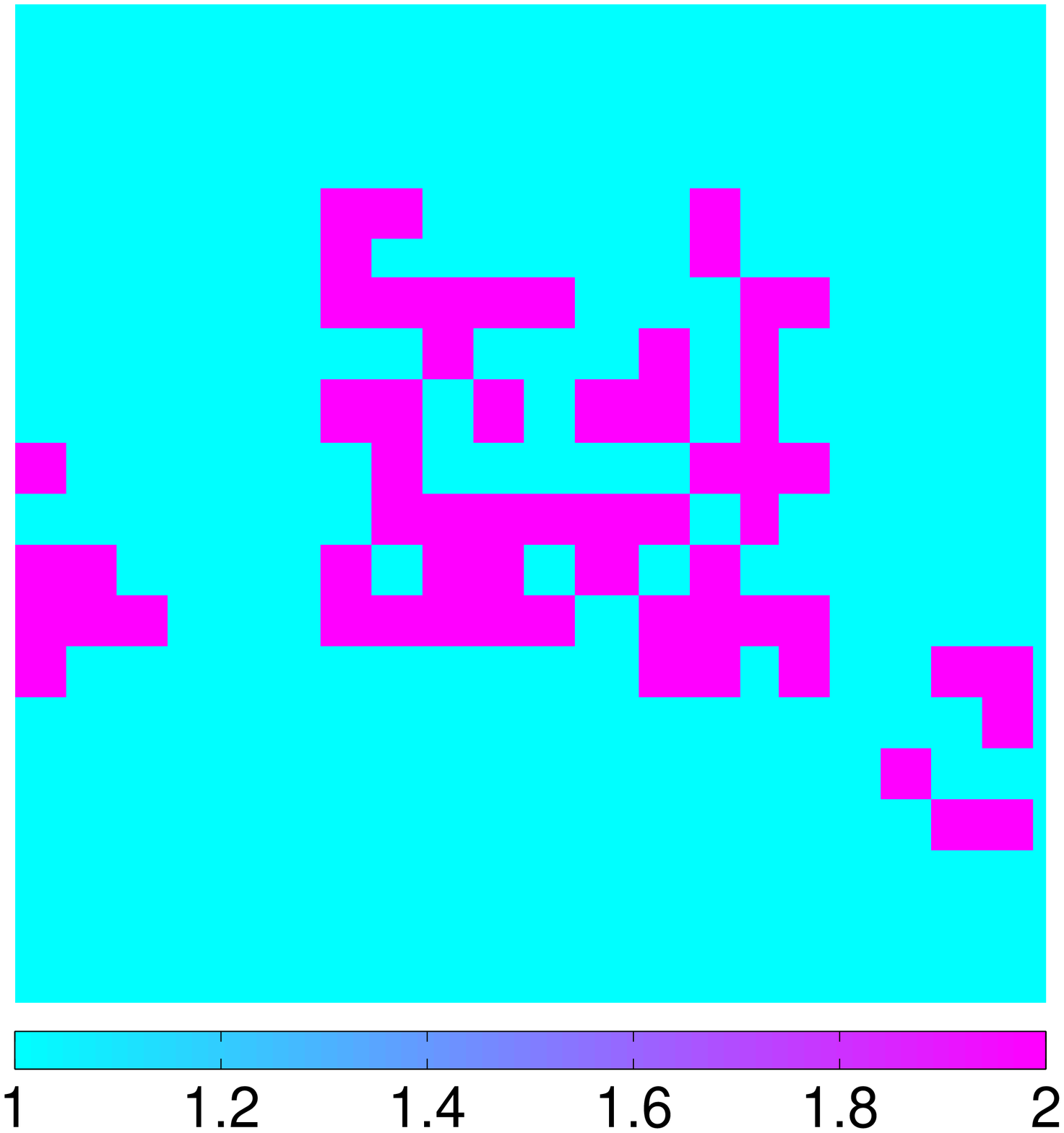}
      \label{ex2tbeta}
      } 
   \subfigure[$|\gamma|^{\frac{1}{2}}$ $(\alpha=0\%)$]{  
     \includegraphics[width=37mm,height=35mm]{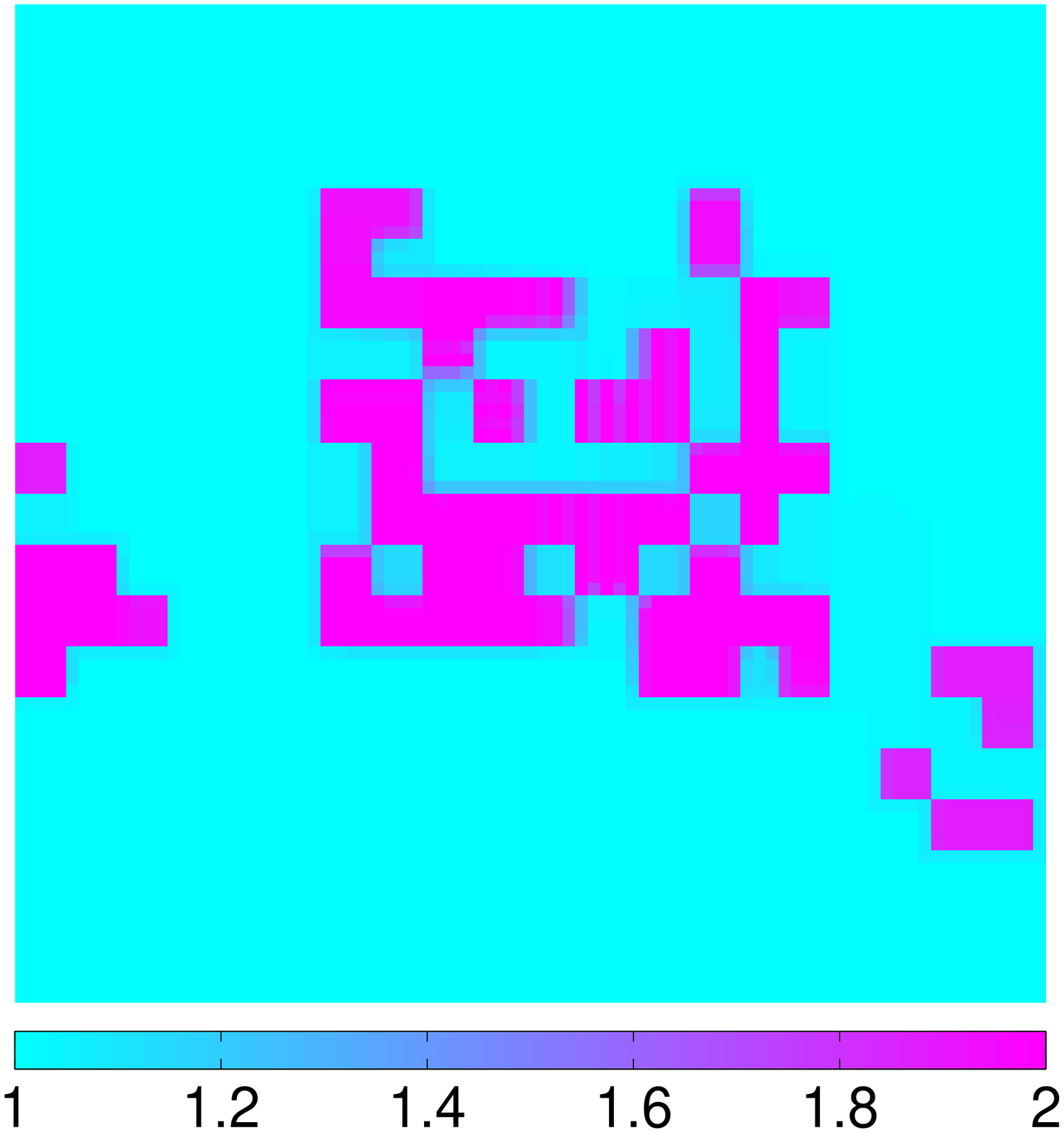}
     \label{ex2cbeta}
     }
   \subfigure[$|\gamma|^{\frac{1}{2}}$ ($\alpha=4\%$)]{ 
    \includegraphics[width=37mm,height=35mm]{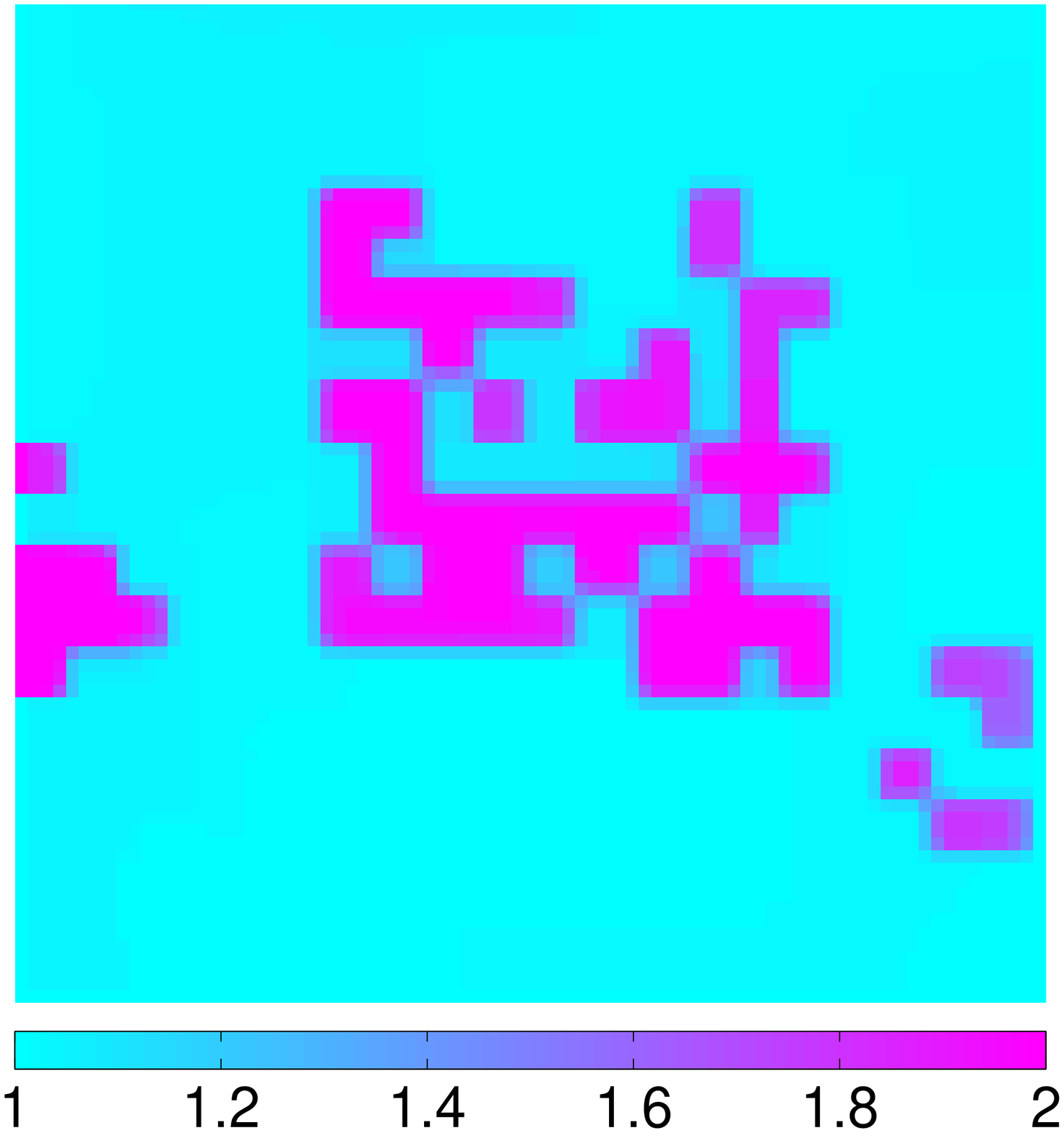}
    \label{ex2nbeta}
    }
   \subfigure[$|\gamma|^{\frac{1}{2}}$ at $\{y=0\}$]{ 
     \includegraphics[trim=10mm 5mm 10mm 0mm,clip,width=35mm,height=38mm]{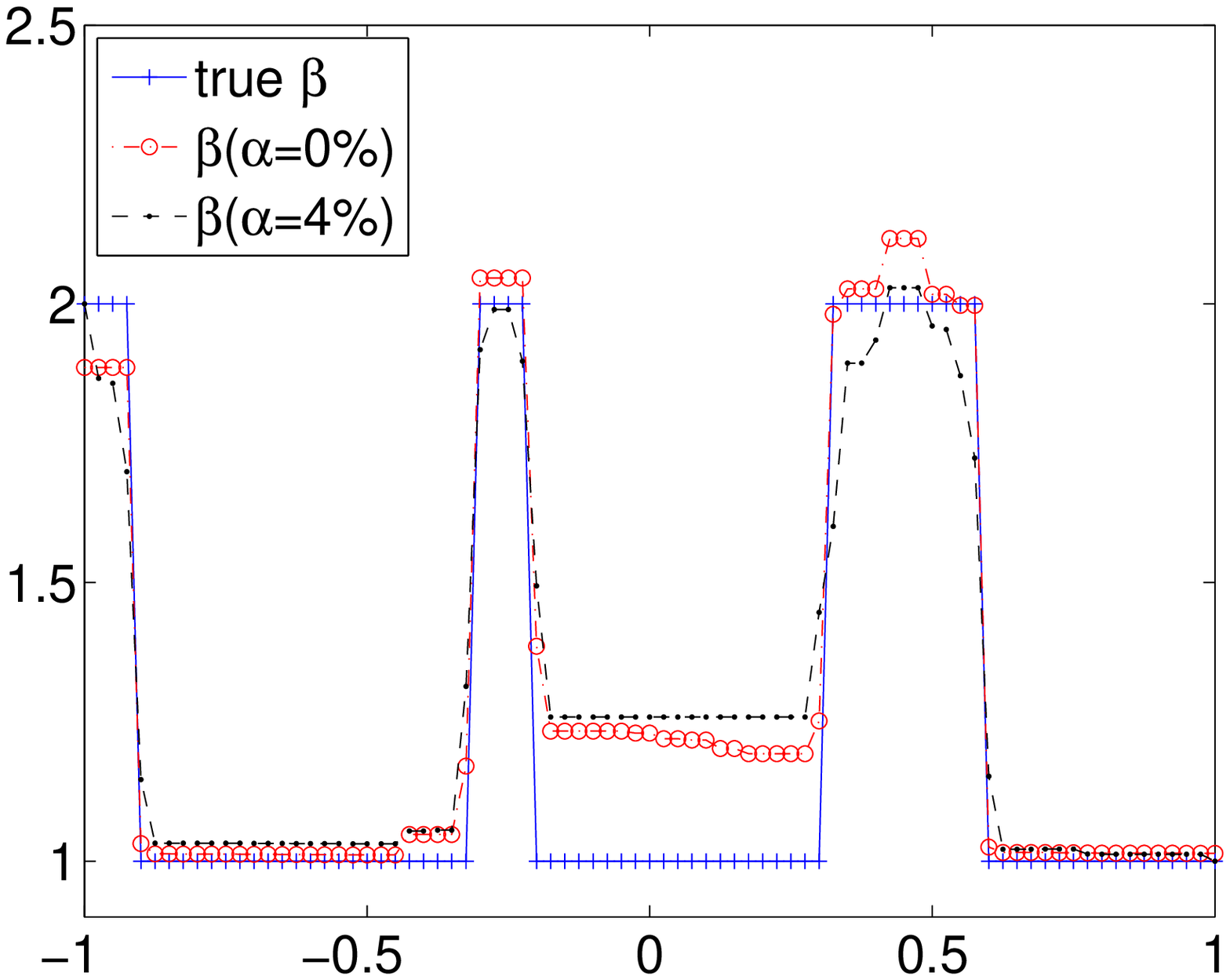} 
     \label{ex2rbeta}
     }

\caption{Experiment 2. \subref{ex2txi}\&\subref{ex2ttau}\&\subref{ex2tbeta}: true values of $(\xi, \zeta, \beta)$. \subref{ex2cxi}\&\subref{ex2ctau}\&\subref{ex2cbeta}: reconstructions with noiseless data. \subref{ex2nxi}\&\subref{ex2ntau}\&\subref{ex2nbeta}: reconstructions with noisy data($\alpha=4\%$). \subref{ex2rxi}\&\subref{ex2rtau}\&\subref{ex2rbeta}: cross sections along $\{y=0\}$.}
\label{E2}
\end{figure}

\paragraph{Experiment 2.} In this experiment, we intend to reconstruct the isotropy given by
\begin{align}
\beta(\mathbf{x})=
\left\{\begin{array}{lll}
1+(\text{sign}(\text{random})+1), \quad &\mathbf{x}\in X_{ij}, \enskip 1\leq i,j\leq 10\\
1+(\text{sign}(\text{random})+1), \quad &\mathbf{x}\in X'_{ij}\cup X''_{ij}, \enskip 1\leq i\leq 3,1\leq j\leq 5\\
1, \enskip &\text{otherwise}
\end{array}\right.
\end{align}
where random is a random number in $[-1, 1]$, $X_{ij}=[0.1(i-1)-0.4, 0.1i-0.4]\times[0.1(j-1)-0.4, 0.1j-0.4]$, $X'_{ij}=[0.1(i-1)-1, 0.1i-1]\times[0.1(j-1)-0.4, 0.1j-0.4]$ and $X''_{ij}=[0.1(i-1)+0.7, 0.1i+0.7]\times[0.1(j-1)-0.8, 0.1j-0.8]$. The anisotropy characterized by $(\xi,\zeta)$ is the same as Experiment 1. The measurements are constructed with the $5$ illuminations given by \eqref{5 bc}. Reconstructions with noise-free and noisy data are performed with a $l_2$ regularization for the anisotropy and $l_1$ regularization using the split Bregman iteration method for the isotropic component. The noise level $\alpha=4\%$. The numerical results of the numerical experiment are shown in Figure \ref{E2}. The relative $L^2$ errors in the reconstructions are $\mathcal{E}^C_{\xi}=2.8\%$, $\mathcal{E}^N_{\xi}=3.7\%$, $\mathcal{E}^C_{\zeta}=6.9\%$, $\mathcal{E}^N_{\zeta}=11.8\%$, $\mathcal{E}^C_{\beta}=5.1\%$ and $\mathcal{E}^N_{\beta}=8.2\%$, respectively.

\begin{figure}
  \centering
   \subfigure[true $\xi$]{
      \includegraphics[width=37mm,height=35mm]{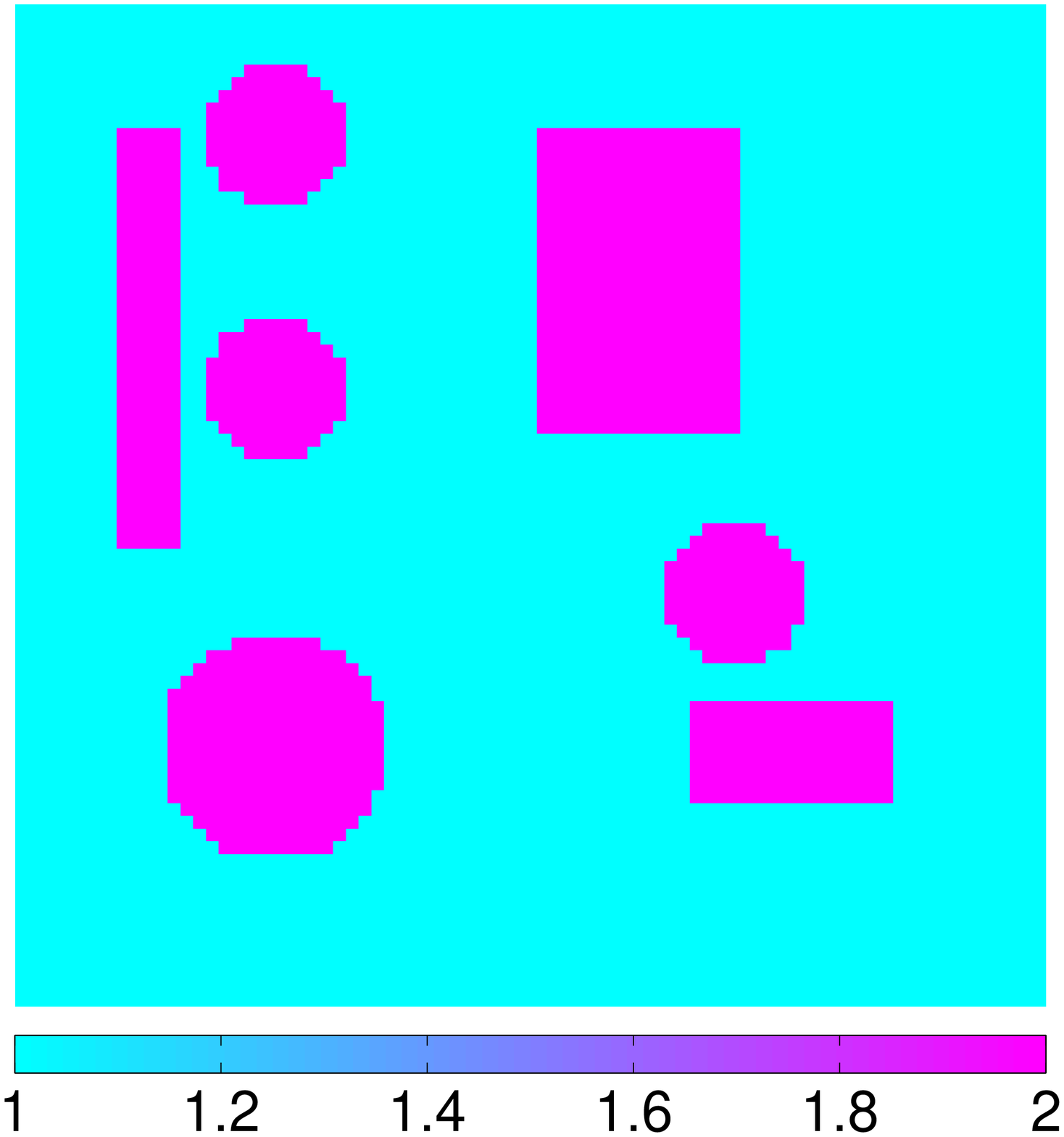}
      \label{ex3txi}
      } 
  \subfigure[$\xi$ ($\alpha=0\%$)]{    
     \includegraphics[width=37mm,height=35mm]{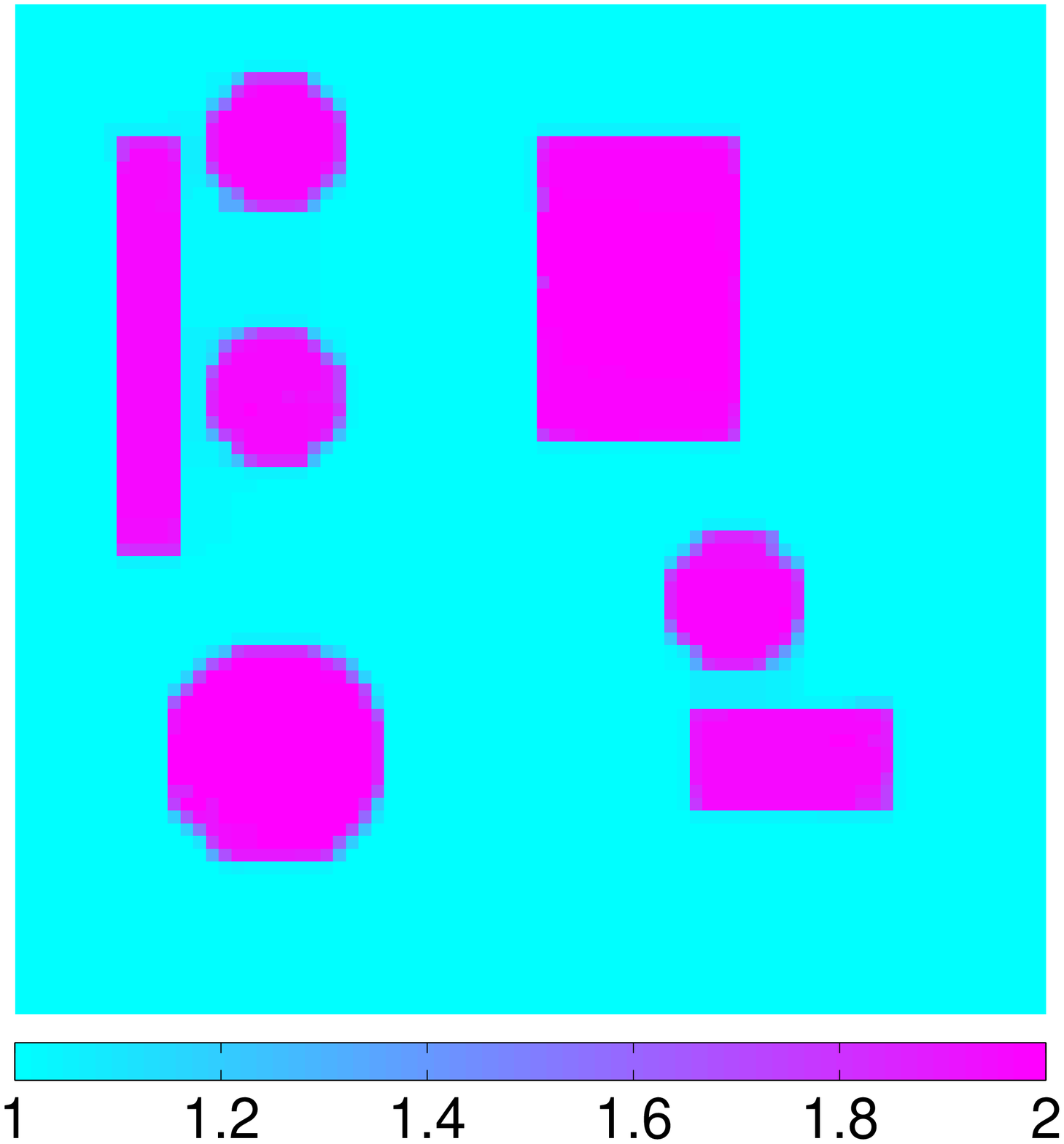}
      \label{ex3cxi}
     }
   \subfigure[$\xi$ ($\alpha=4\%$)]{
    \includegraphics[width=37mm,height=35mm]{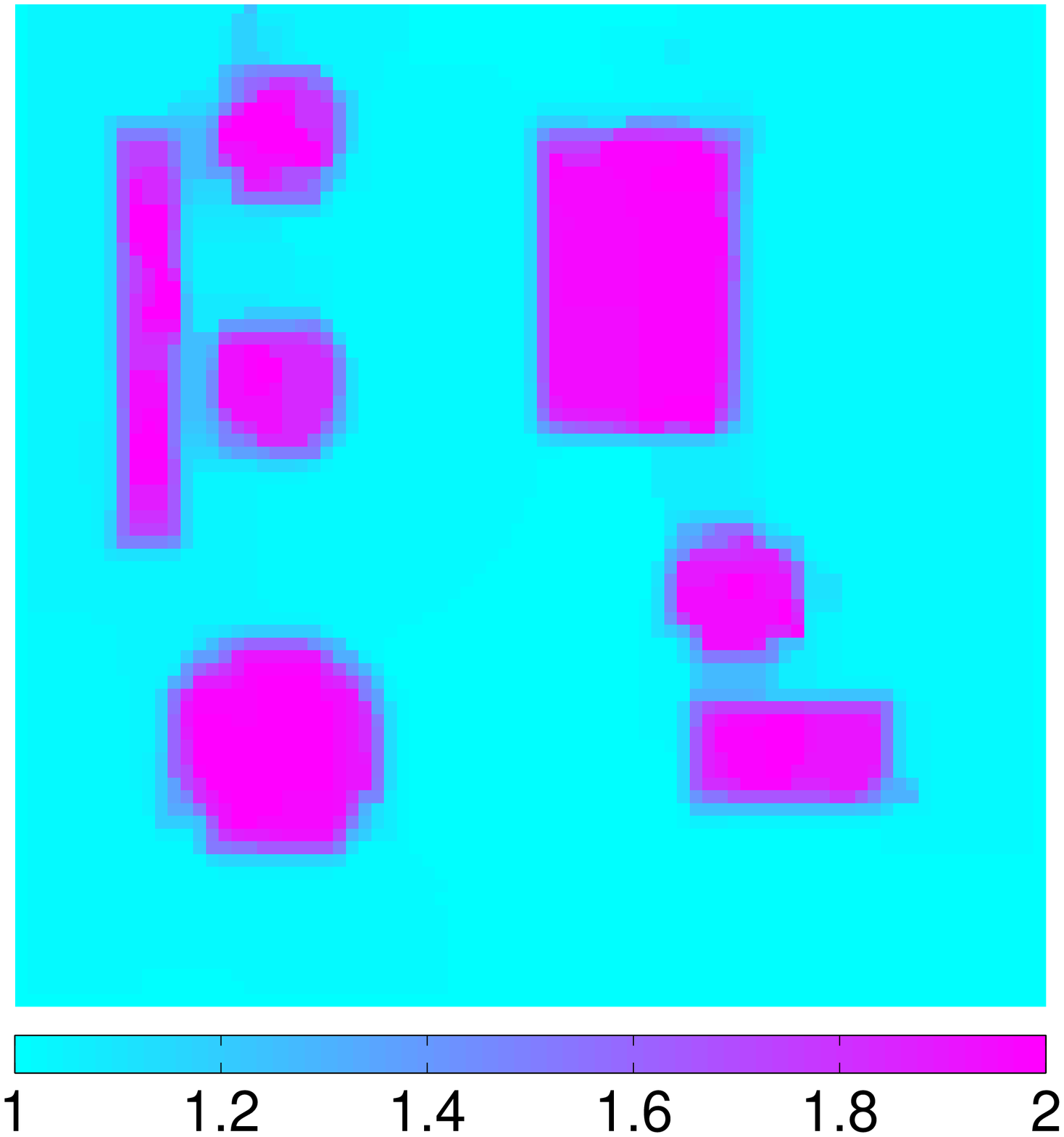}
     \label{ex3nxi}
    }
    \subfigure[$\xi$ at \{$y=-0.5$\}]{
     \includegraphics[trim=10mm 5mm 10mm 0mm,clip,width=35mm,height=38mm]{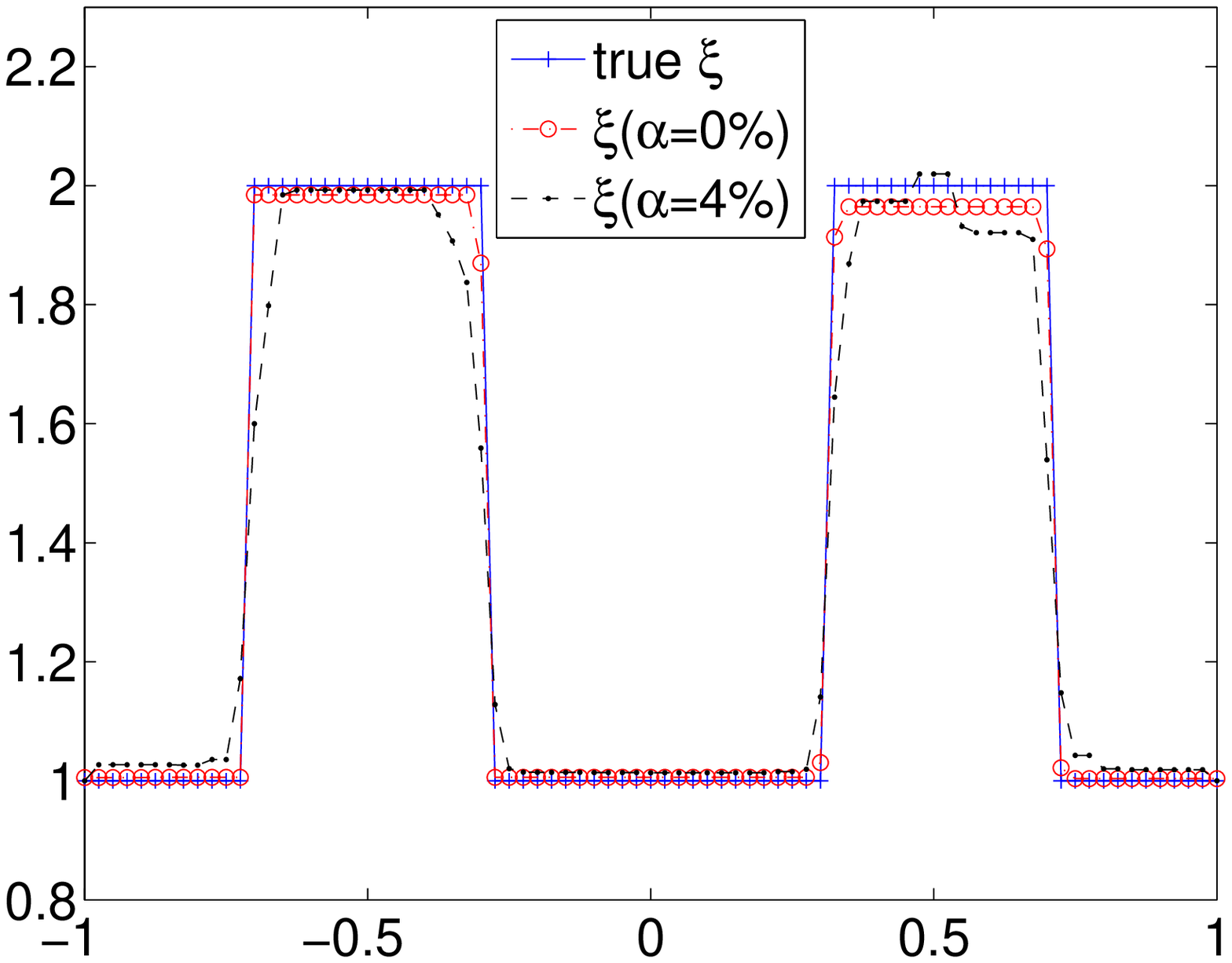} 
      \label{ex3rxi}
     }
     
     \subfigure[true $\zeta$]{
      \includegraphics[width=37mm,height=35mm]{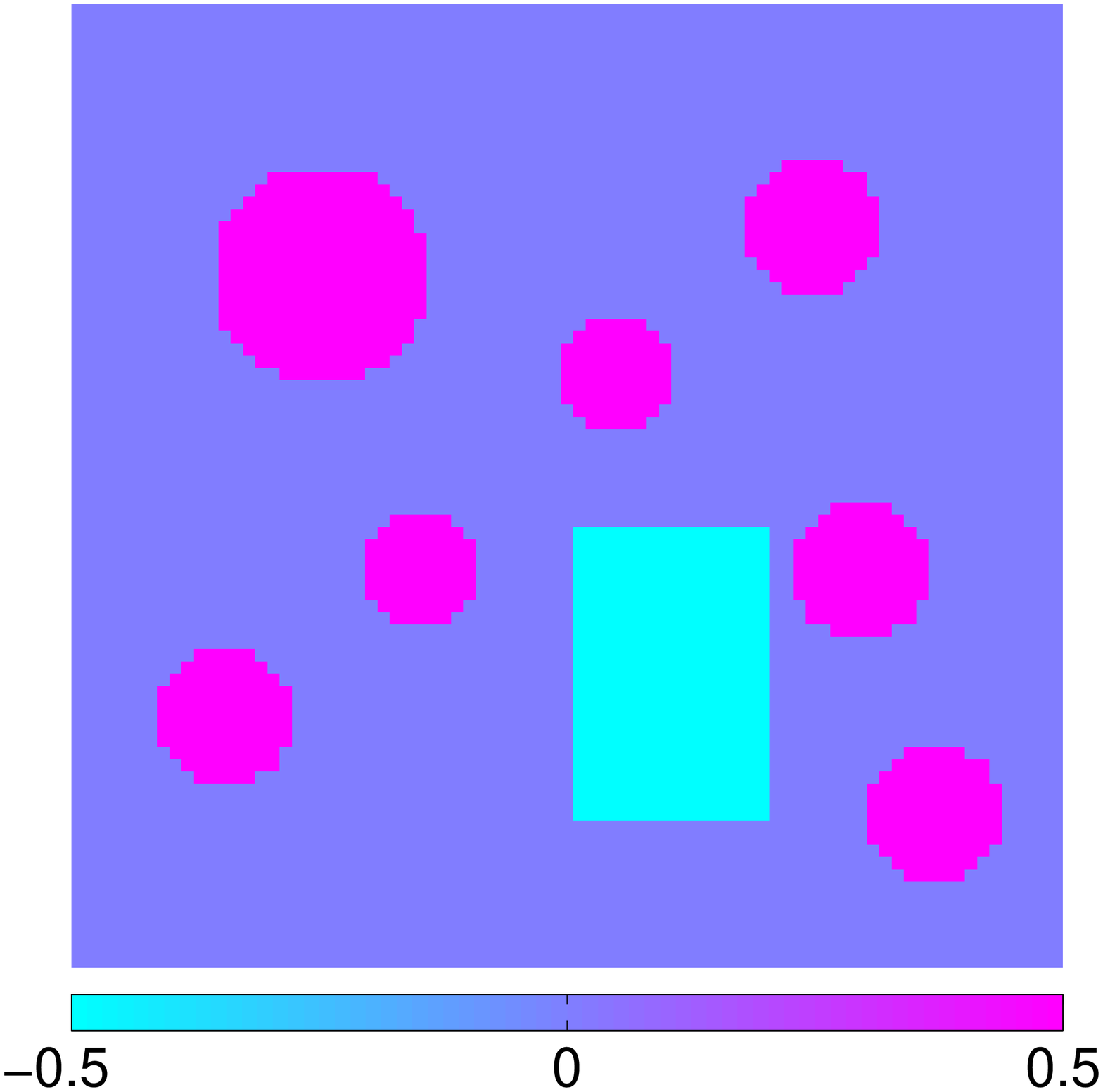}
       \label{ex3ttau}
      } 
 \subfigure[$\zeta$ ($\alpha=0\%$)]{ 
     \includegraphics[width=37mm,height=35mm]{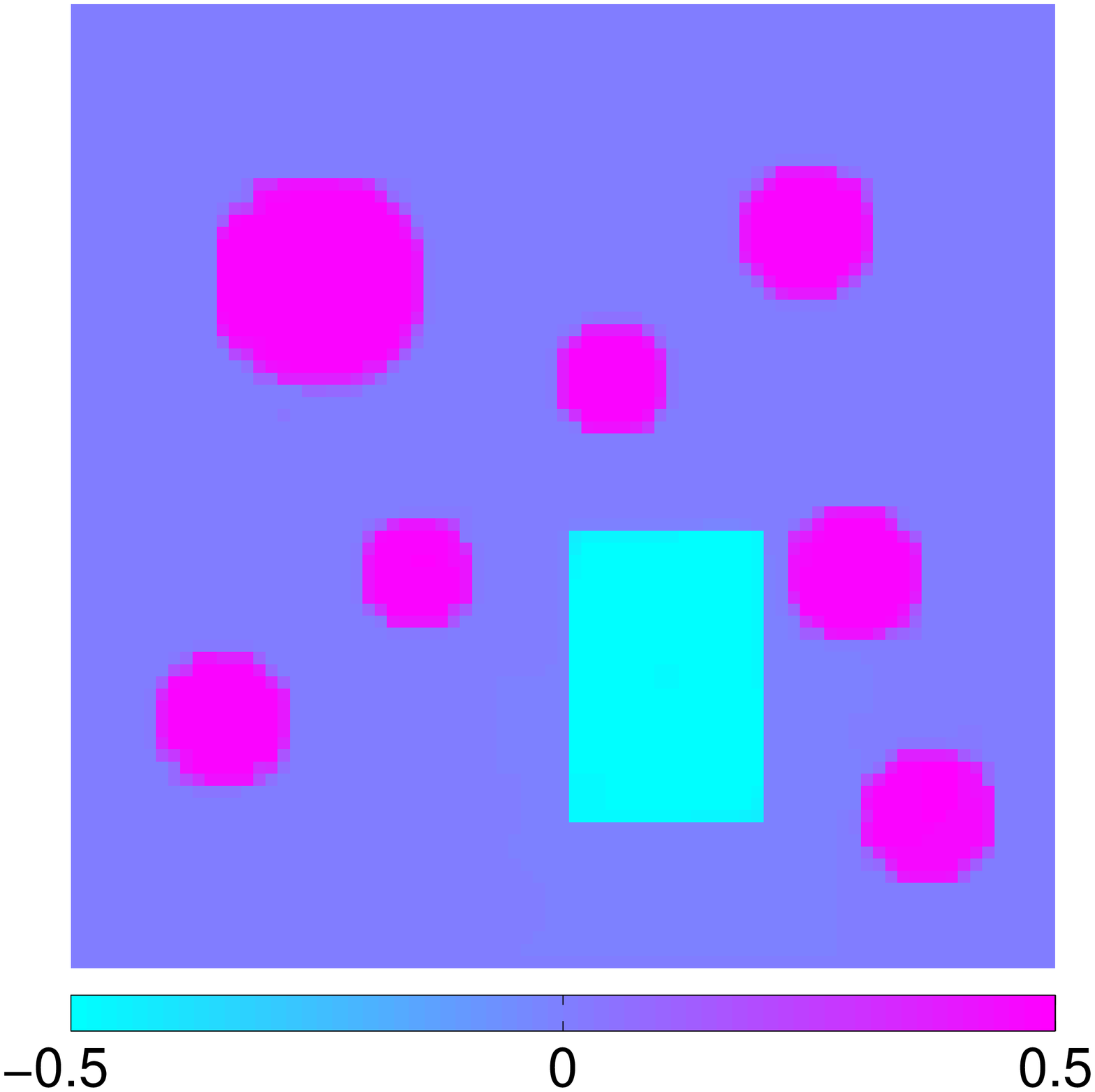}
     \label{ex3ctau}
     }
   \subfigure[$\zeta$ ($\alpha=4\%$)]{ 
    \includegraphics[width=37mm,height=35mm]{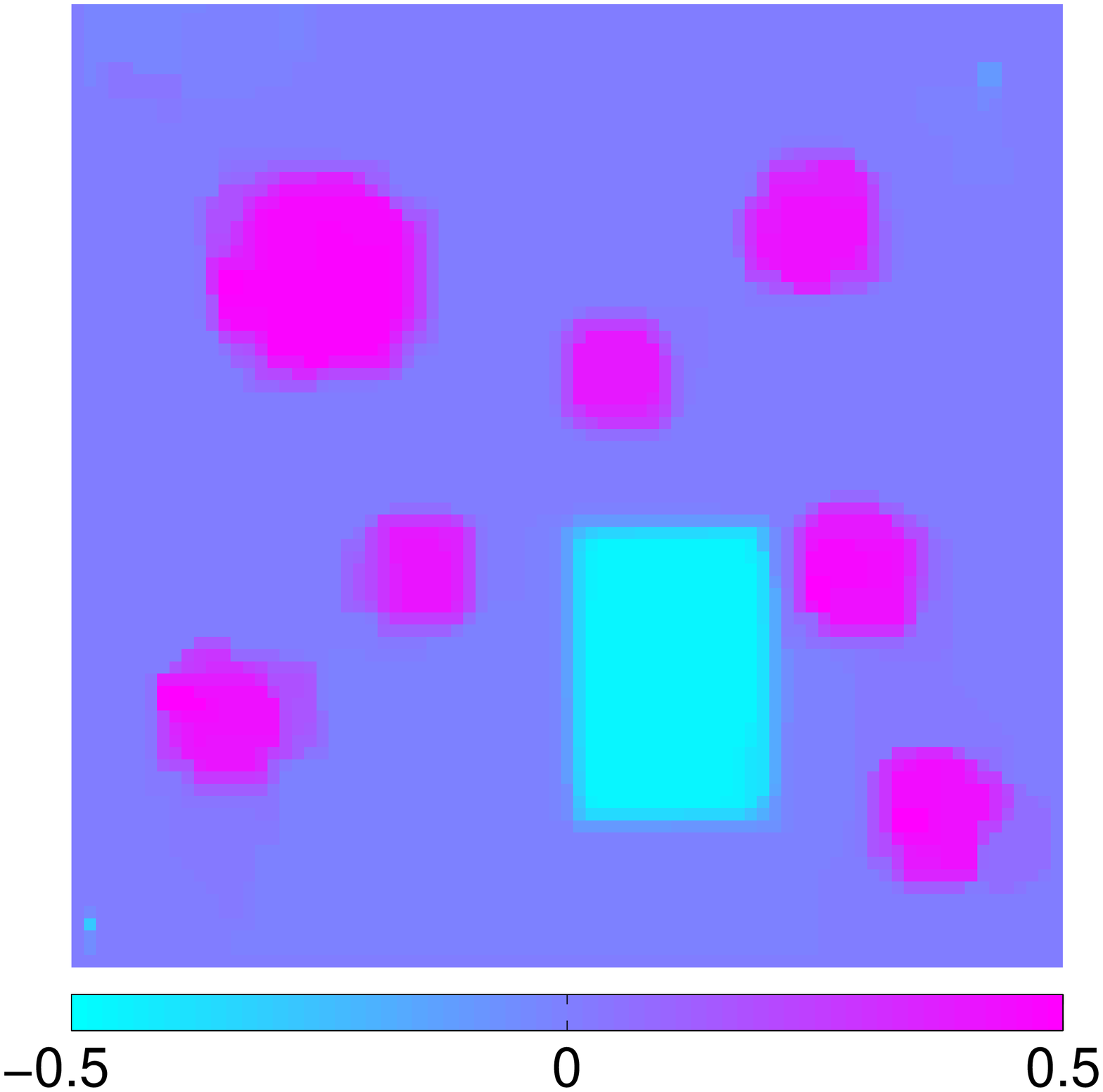}
    \label{ex3ntau}
    }
     \subfigure[$\zeta$ at \{$y=-0.5$\}]{
     \includegraphics[trim=10mm 5mm 10mm 0mm,clip,width=35mm,height=38mm]{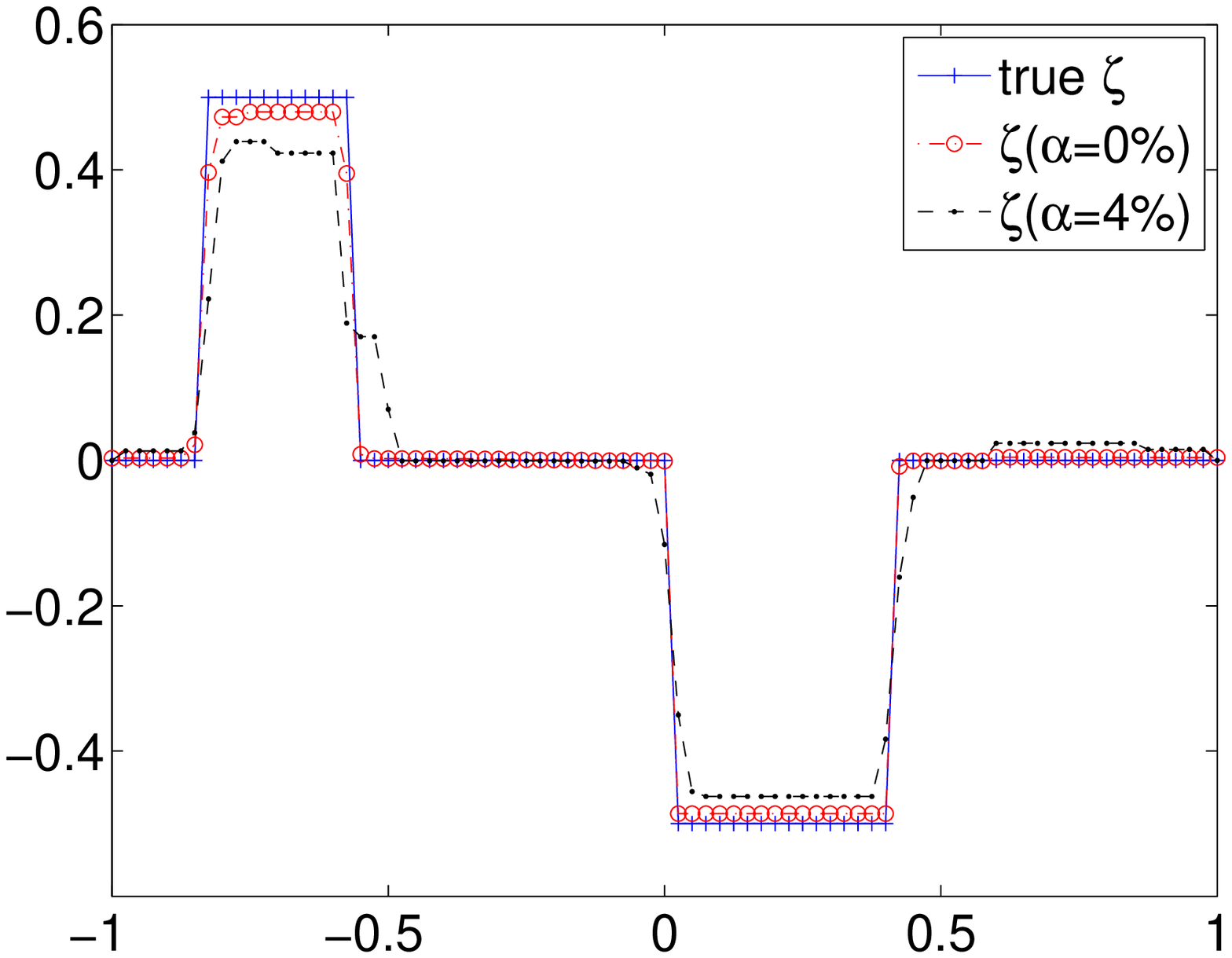} 
     \label{ex3rtau}
     }
       
     \subfigure[true $|\gamma|^{\frac{1}{2}}$]{
      \includegraphics[width=37mm,height=35mm]{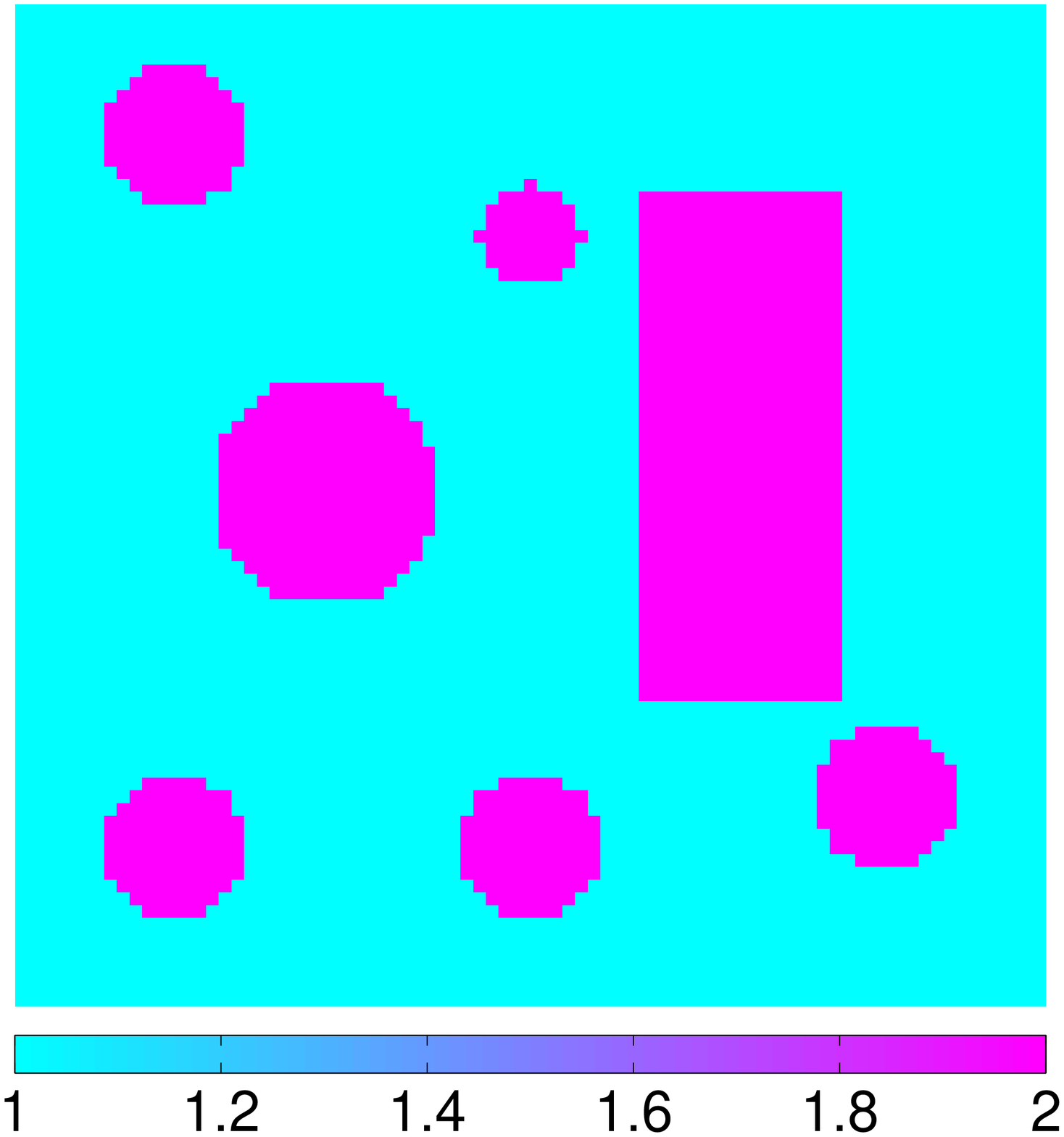}
      \label{ex3tbeta}
      } 
   \subfigure[$|\gamma|^{\frac{1}{2}}$ $(\alpha=0\%)$]{  
     \includegraphics[width=37mm,height=35mm]{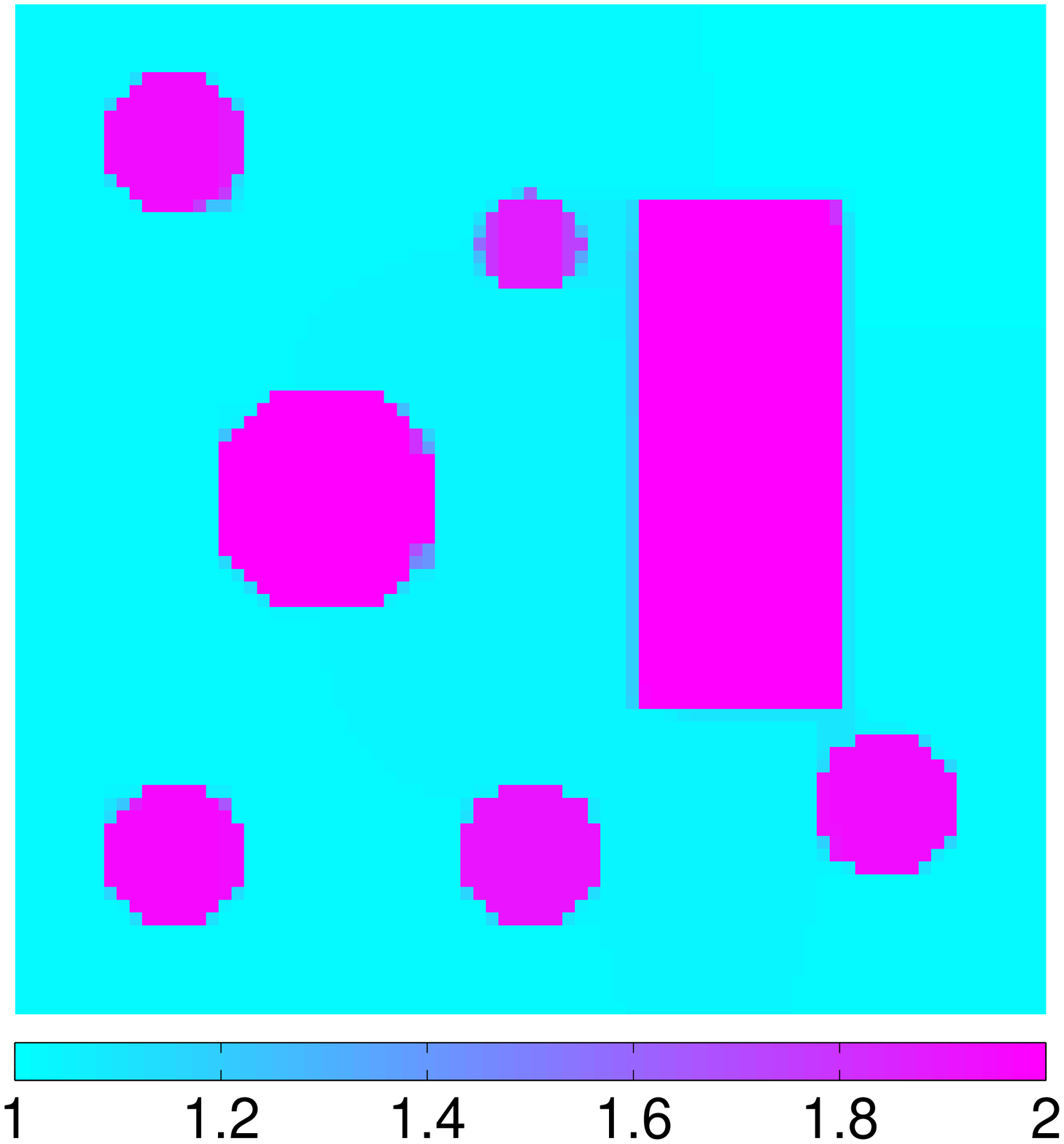}
     \label{ex3cbeta}
     }
   \subfigure[$|\gamma|^{\frac{1}{2}}$ ($\alpha=4\%$)]{ 
    \includegraphics[width=37mm,height=35mm]{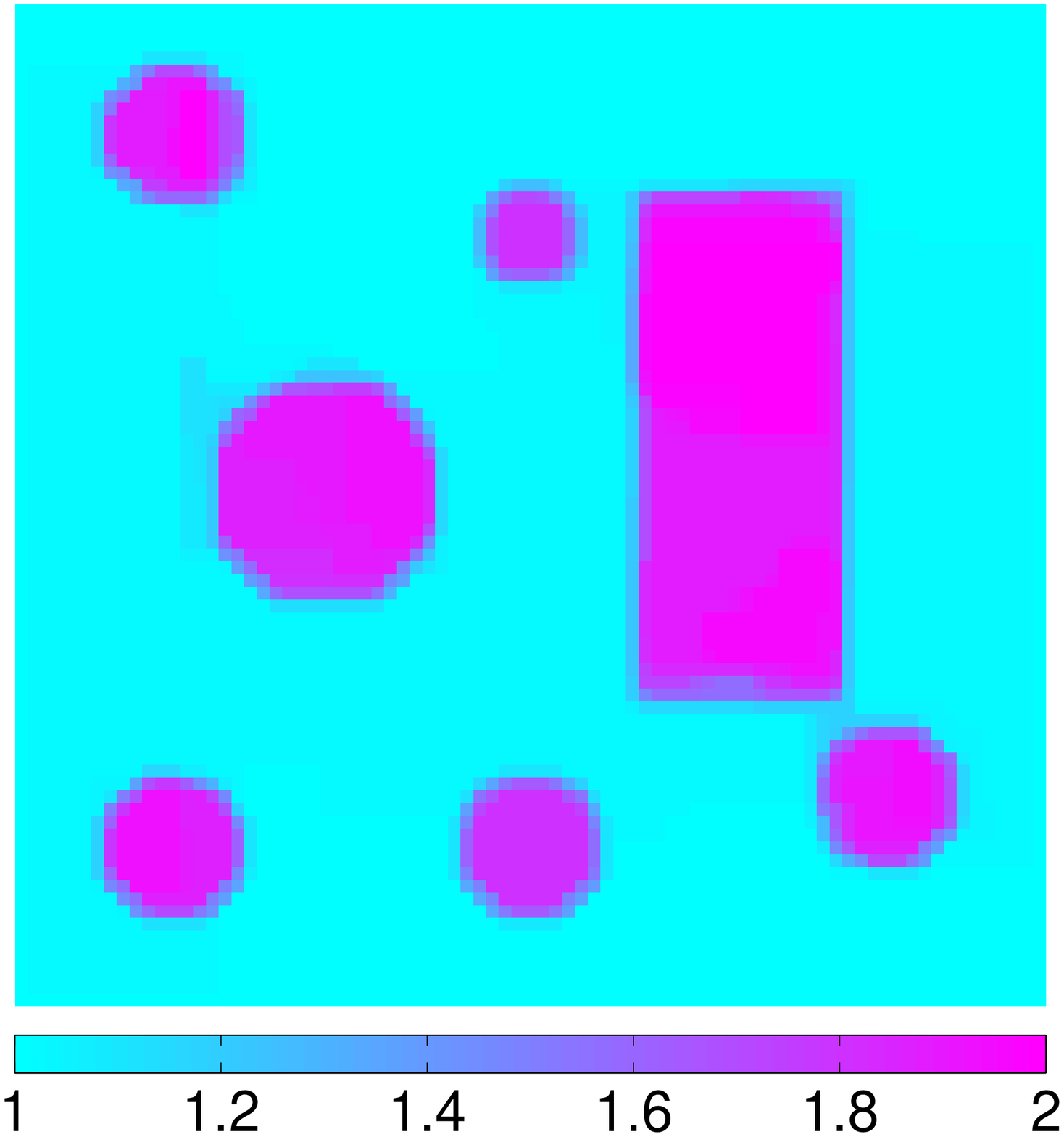}
    \label{ex3nbeta}
    }
   \subfigure[$|\gamma|^{\frac{1}{2}}$ at $\{y=-0.5\}$]{ 
     \includegraphics[trim=10mm 5mm 10mm 0mm,clip,width=35mm,height=38mm]{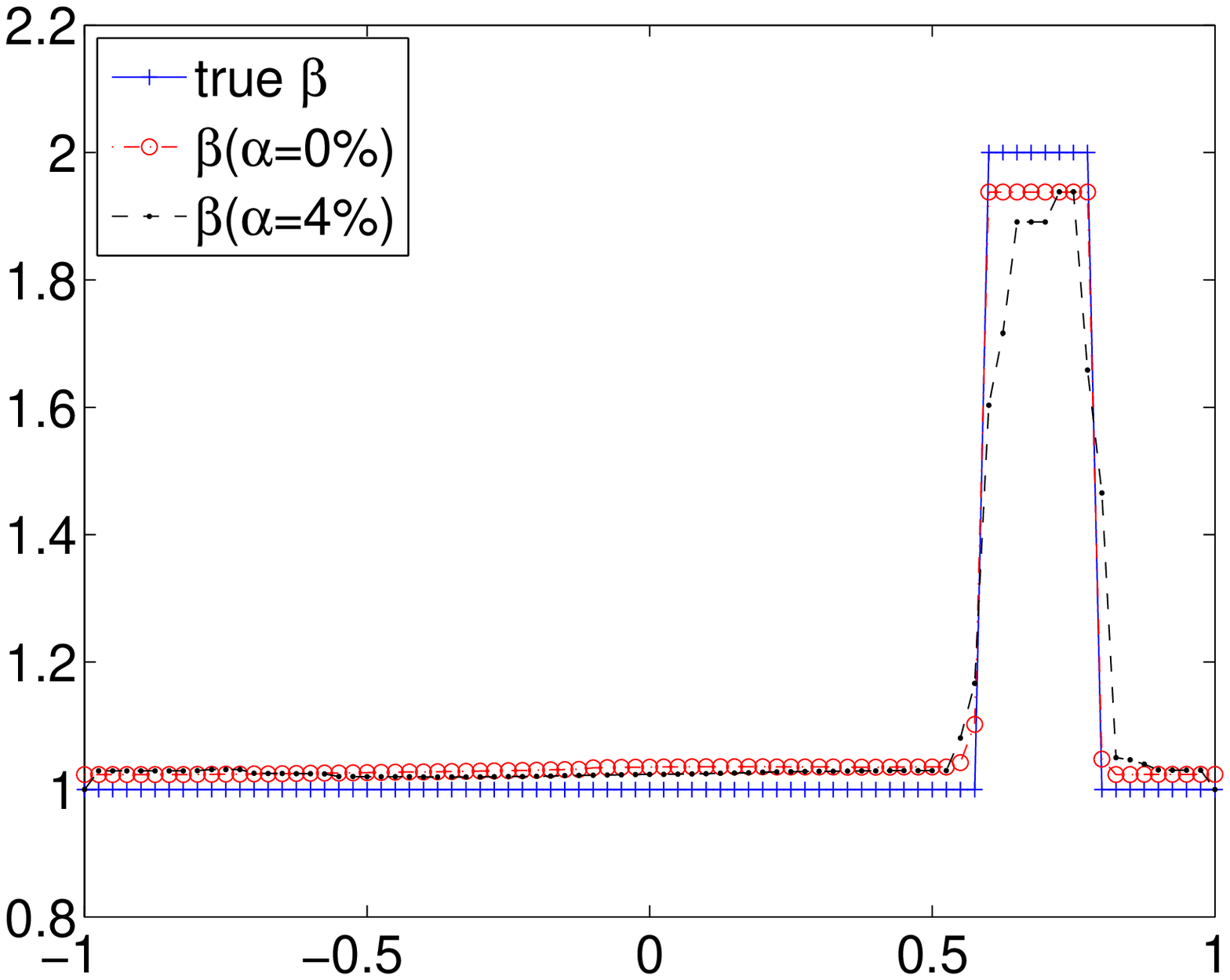} 
     \label{ex3rbeta}
     }
  
\caption{Experiment 3. \subref{ex3txi}\&\subref{ex3ttau}\&\subref{ex3tbeta}: true values of $(\xi, \zeta, \beta)$. \subref{ex3cxi}\&\subref{ex3ctau}\&\subref{ex3cbeta}: reconstructions with noiseless data. \subref{ex3nxi}\&\subref{ex3ntau}\&\subref{ex3nbeta}: reconstructions with noisy data($\alpha=4\%$). \subref{ex3rxi}\&\subref{ex3rtau}\&\subref{ex3rbeta}: cross sections along $\{y=-0.5\}$.}
\label{E3}
\end{figure}

\paragraph{Experiment 3.}In this experiment, we attempt to reconstruct coefficients with discontinuities. To simplify the implementation, we only consider piecewise constant coefficients. Here we use the same illuminations as in Experiment 1. Reconstructions with both noiseless and noisy data are performed with $l_1$ regularization using the split Bregman iteration method for both the anisotropic and isotropic components. The noise level $\alpha=4\%$. The results of the numerical experiment are shown in Figure \ref{E3}. From the figures, we observe that the singularities of the coefficients create minor artifacts on the reconstructions and the error in the reconstruction is larger at the discontinuities than in the rest of the domain. The relative $L^2$ errors in the reconstructions are $\mathcal{E}^C_{\xi}=3.9\%$, $\mathcal{E}^N_{\xi}=9.6\%$, $\mathcal{E}^C_{\zeta}=13.4\%$, $\mathcal{E}^N_{\zeta}=31.9\%$, $\mathcal{E}^C_{\beta}=3.7\%$ and $\mathcal{E}^N_{\beta}=8.2\%$.

\subsection{Experiments with control over part of the boundary}

The previous experiments show that the reconstruction of both smooth and discontinuous coefficients is very accurate and robust to noise when one can fully prescribe boundary conditions ensuring conditions $A$\&$B$ of Lemma \ref{cond cgo}. In practice, one does not always have access to the whole boundary, and instead may have to prescribe boundary conditions on only a small part of the domain. In the next series of  experiments, we assume to only have control over the bottom boundary of the square domain $X$, call it $\partial X_B = [-1,1]\times\{-1\}$. Over the rest of the boundary, we successively impose homogeneous Dirichlet boundary conditions (Experiment 4), then homogeneous Neumann boundary conditions (Experiment 5). In two spatial dimensions, either case forces all conductivity solutions to have their gradients to be pairwise collinear (normal to the boundary for Dirichlet conditions, tangential to the boundary for Neumann conditions). This violates both conditions of Lemma \ref{cond cgo}, and we expect reconstructions to do poorly near the uncontrolled part of the boundary. Note that we can predict the accuracy of the reconstruction from the measured data since the constants of independence appearing in  $A$\&$B$ in Lemma \ref{cond cgo} can be estimated from the measurements $\{H_j\}_j$. For instance, of two measurements are not sufficiently linearly independent, then additional measurements may be considered before the reconstruction formulas are applied. 

Note that in higher spatial dimensions, the practically more relevant homogeneous Neumann conditions should lead to better reconstructions as these conditions impose less constraints on gradients than homogeneous Dirichlet conditions. 

\paragraph{Experiment 4.} We now repeat Experiment 3 using illuminations that are only non-zero on the bottom boundary of the domain. \\

\begin{figure}
  \centering
  \subfigure[$\xi$ ($\alpha=0\%$)]{ 
    \includegraphics[width=37mm,height=35mm]{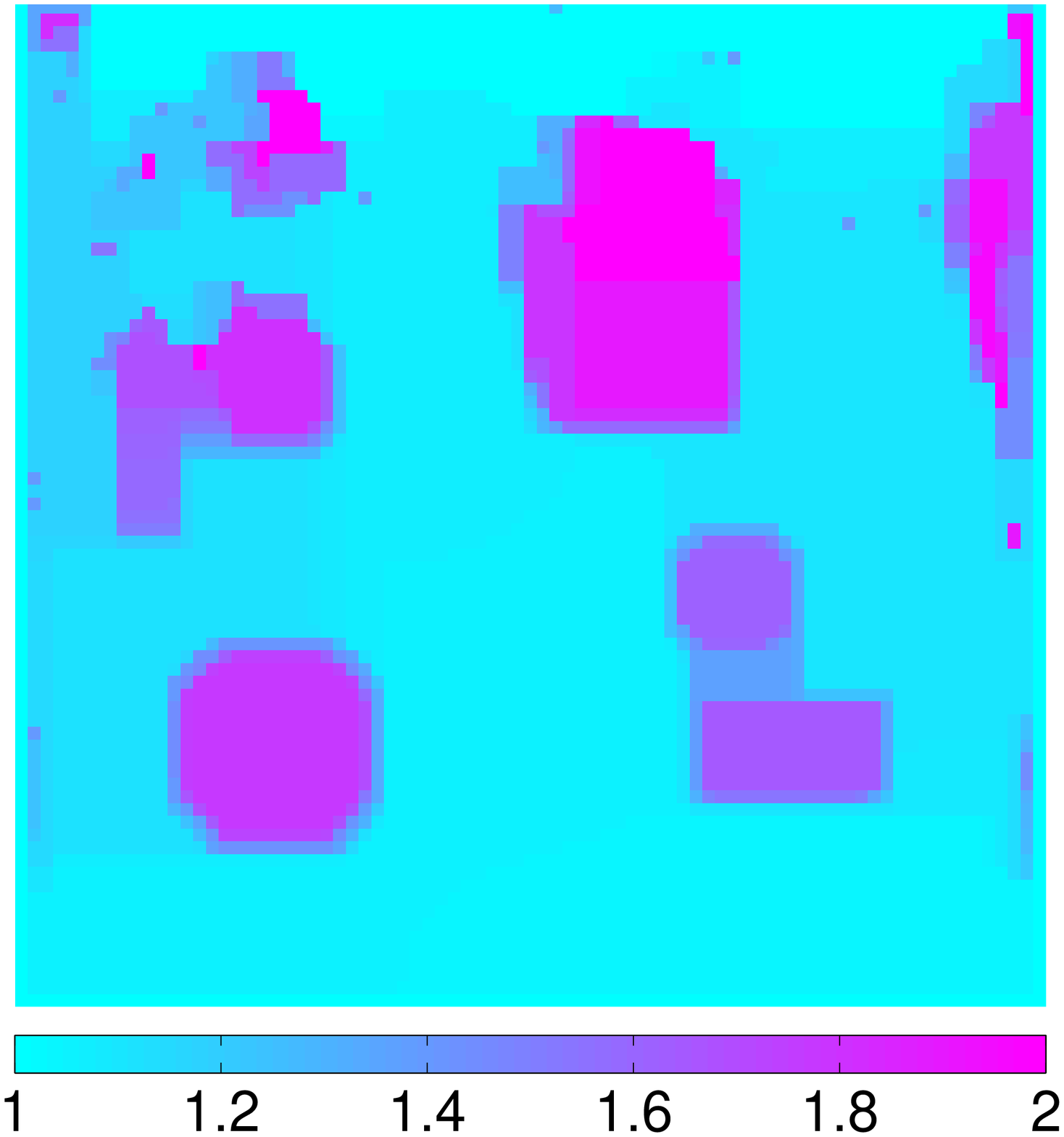}
    \label{smallxi}
    }   
    \subfigure[$\zeta$ ($\alpha=0\%$)]{ 
    \includegraphics[width=37mm,height=35mm]{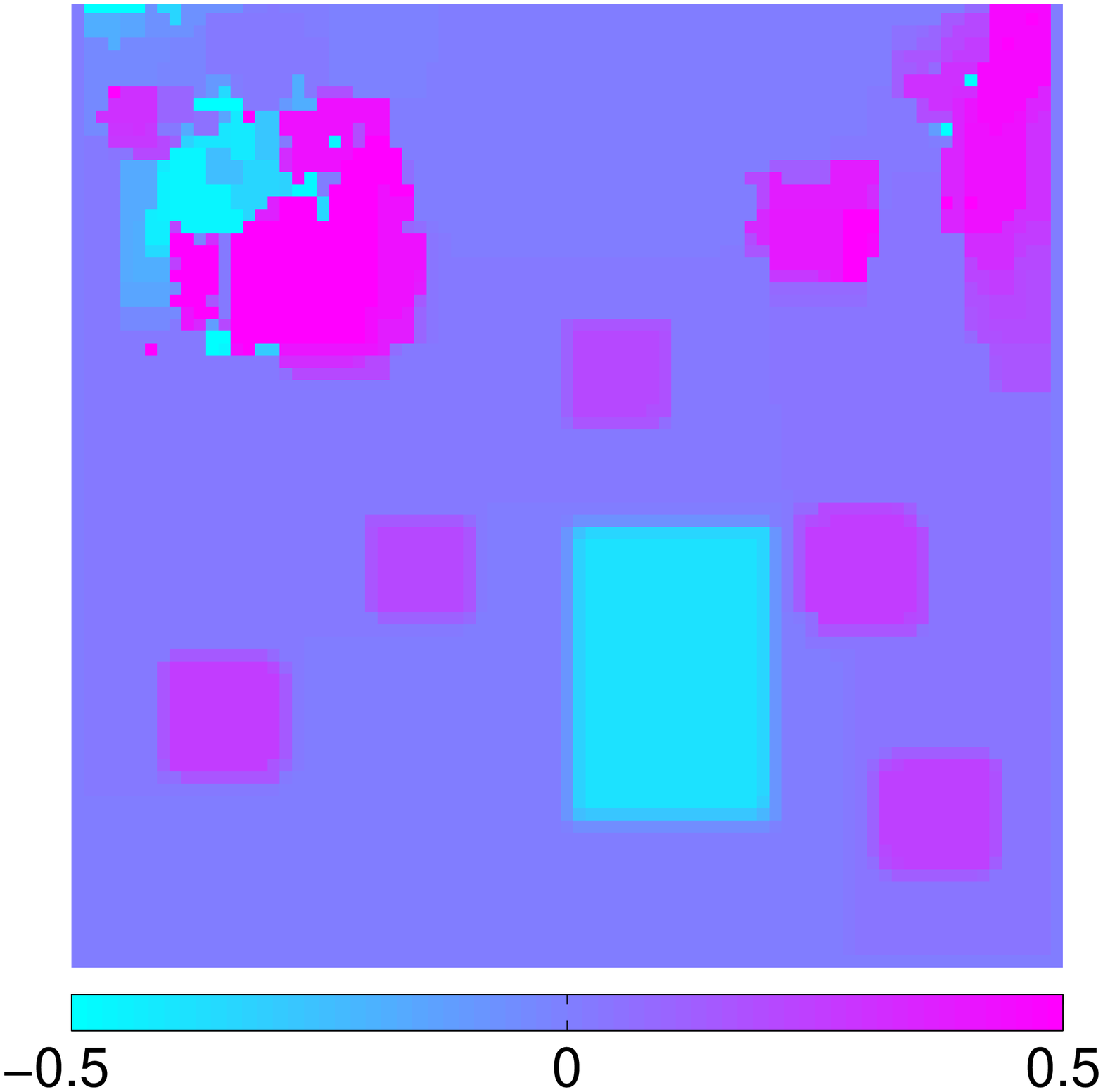}
    \label{smalltau}
    }
    \subfigure[$\zeta$ ($\alpha=0\%$)]{ 
    \includegraphics[width=37mm,height=35mm]{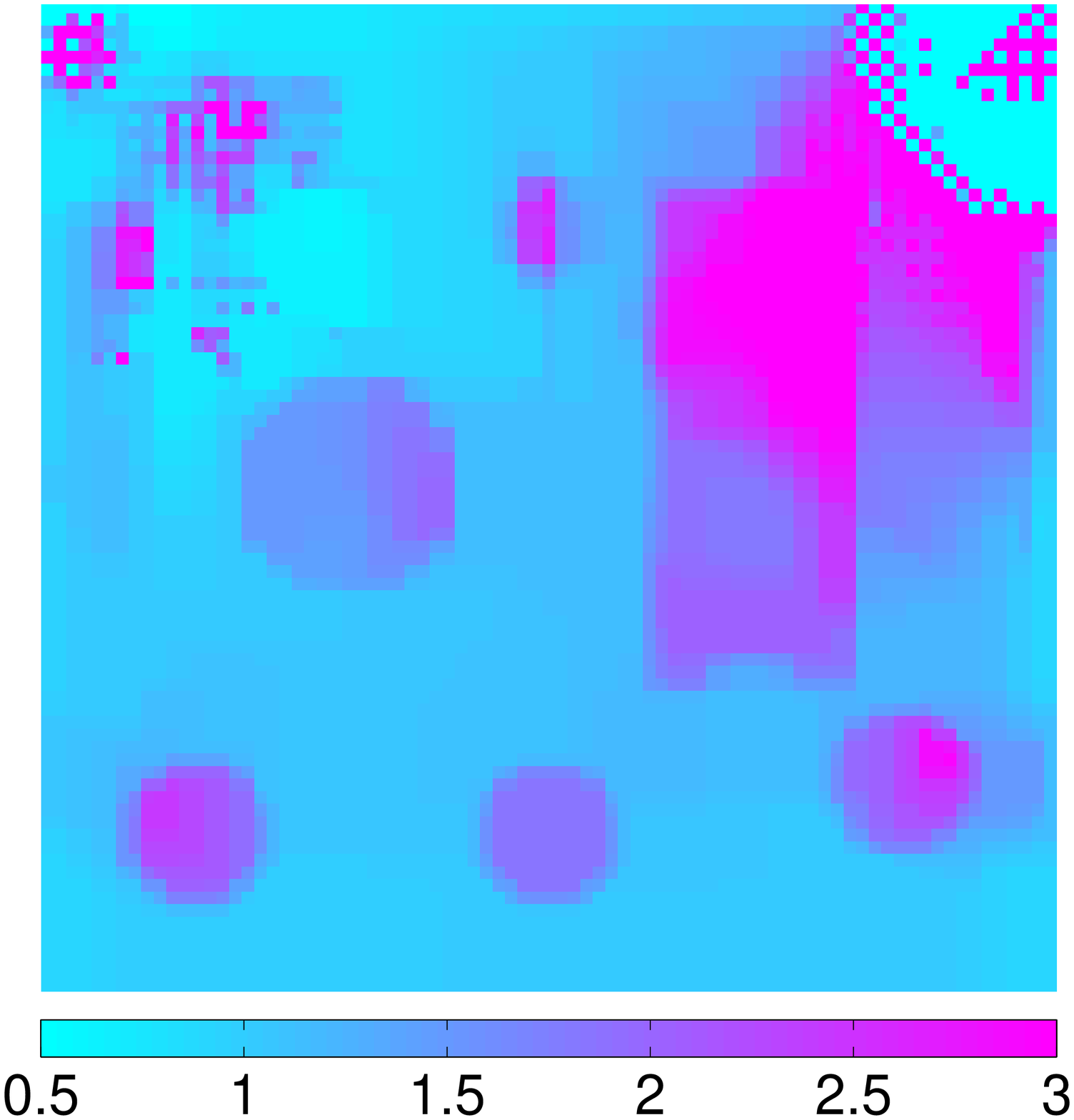}
    \label{smallbeta}
    }
    \subfigure[$\log|\det(\nabla u_i,\nabla u_j)|$]{ 
    \includegraphics[trim=15mm 5mm 14mm 0mm,clip,width=35mm,height=38mm]{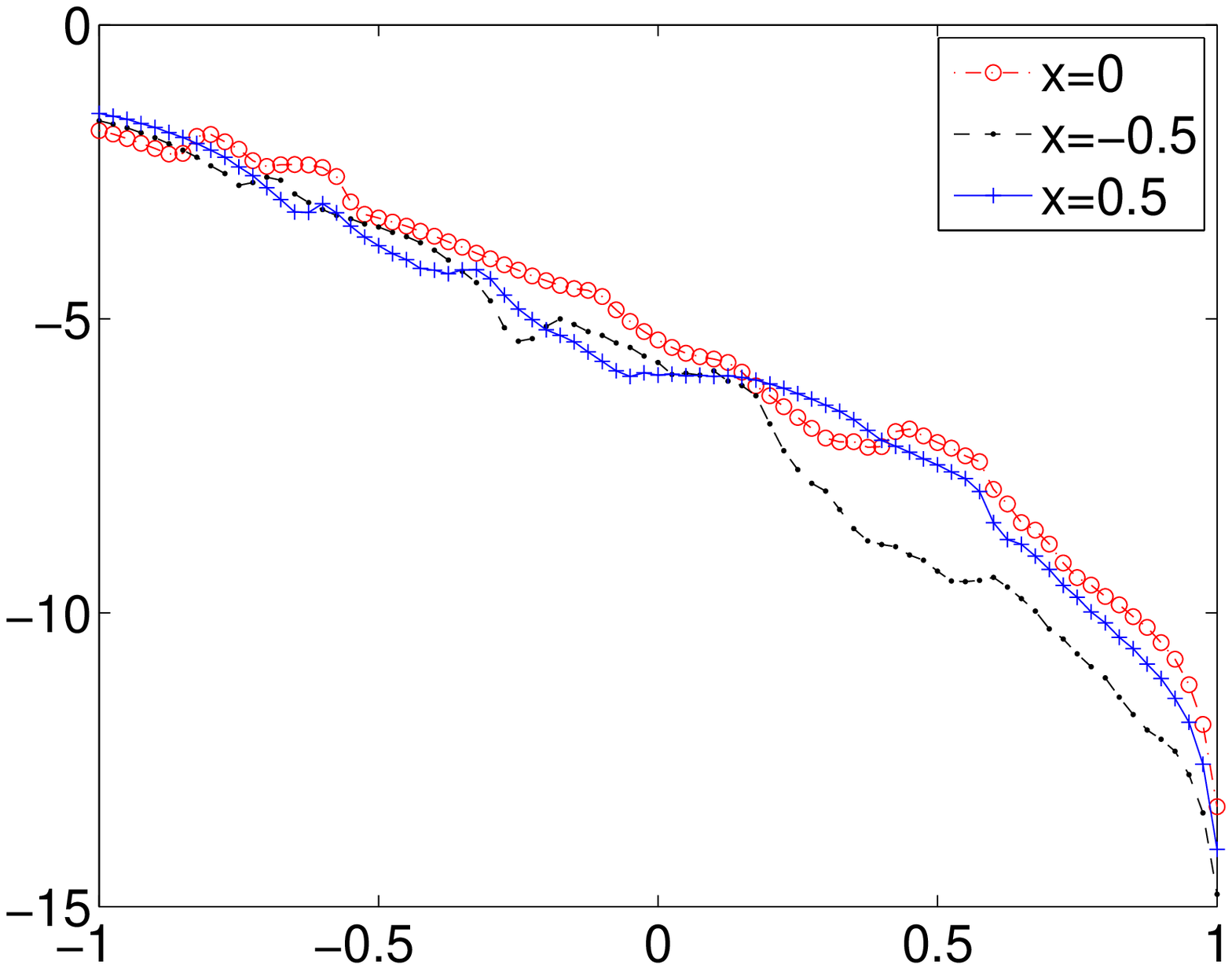}
         \label{logdert}    
    }     
         \caption{Simulations on $X$. \subref{smallxi}\&\subref{smalltau}\&\subref{smallbeta}: reconstructions with noiseless data. \subref{logdert}: cross section of $\max\limits_{1\leq i<j\leq 5}\log|\det(\nabla u_i,\nabla u_j)|$ along $\{x=0\}$, $\{x=-0.5\}$ and $\{x=0.5\}$.}
\label{E4small}
\end{figure}

\noindent{\emph{Reconstructions of the anisotropy $\tilde\gamma$} in $[-1,1]^2$.} We first perform the reconstructions of $\xi$ and $\zeta$. We use five illuminations given by Gaussian functions as follows,
\begin{align}
g_i(\mathbf{x}) = \left\{\begin{array}{ll}
(2\pi\cdot0.2^2)^{-\frac{1}{2}}\exp\{-\frac{1}{2\cdot0.2^2}(x+x_i)^2\}, \enskip &\mathbf{x}\in \partial X_B\\
0, \enskip &\mathbf{x}\in \partial X\setminus\partial X_B
\end{array}\right. \quad 1\leq i\leq 5
\end{align}
where $\{\mathbf{x}_i\}_{1\leq i\leq 5}=\{-0.8,-0.4,0,0.4,0.8\}$. Reconstructions with noise-free data are shown in Figure \ref{E4small}. From this simulation, we can see that even with noise-free data, the reconstruction degrades as one gets farther away from the controlled boundary $\partial X_B$, while it remains accurate near $\partial X_B$. 

\bigbreak
\noindent{\emph{Reconstructions of $\gamma$ in an extended domain.} From the numerical simulation in Figure \ref{E4small}, it is clear that the reconstruction procedure does not perform well for $\mathbf{x}$ far from $\partial X_B$. From Fig.\ref{logdert}, we can see that $\det(\nabla u_i,\nabla u_j)$ decays very rapidly, which means that Lemma \ref{cond cgo}.$A$ is not fulfilled. 

A way to scan a deeper part of the domain with conductivity solutions of linearly independent gradients is obtained by spreading out the various boundary conditions along the $x$-axis. To this end, we now extend the domain $X$ to $X'= [-3,3]\times [-1.2,4.8]$ and use a $\mathsf{N'+1\times N'+1}$ square grid with $\mathsf{N}'=240$. We use the following five Gaussian functions as illuminations,
\begin{align}
g_i(\mathbf{x}) = \left\{\begin{array}{ll}
(2\pi\cdot0.2^2)^{-\frac{1}{2}}\exp\{-\frac{1}{2\cdot0.2^2}(x+x_i)^2\}, \enskip &\mathbf{x}\in \partial X'_B\\
0, \enskip &\mathbf{x}\in \partial X'\setminus\partial X'_B
\end{array}\right. \quad 1\leq i\leq 5
\end{align}
where $\{\mathbf{x}_i\}_{1\leq i\leq 5}=\{-2.8,-1.5,0,1.5,2.8\}$. The reconstruction of the anisotropy $\tilde\gamma$ with noise free data is shown in Figure \ref{E4large}. In this setting, we see that the domain $X$ is now fully covered by conductivity solutions whose gradients fulfill condition $A$ from Lemma \ref{cond cgo}, and the reconstruction performs well everywhere on $X$. On the other hand, as expected, the reconstruction does not perform well outside $X$.

\begin{figure}[htb]
  \centering
    \subfigure[true $\xi$]{ 
     \includegraphics[width=36.5mm,height=35mm]{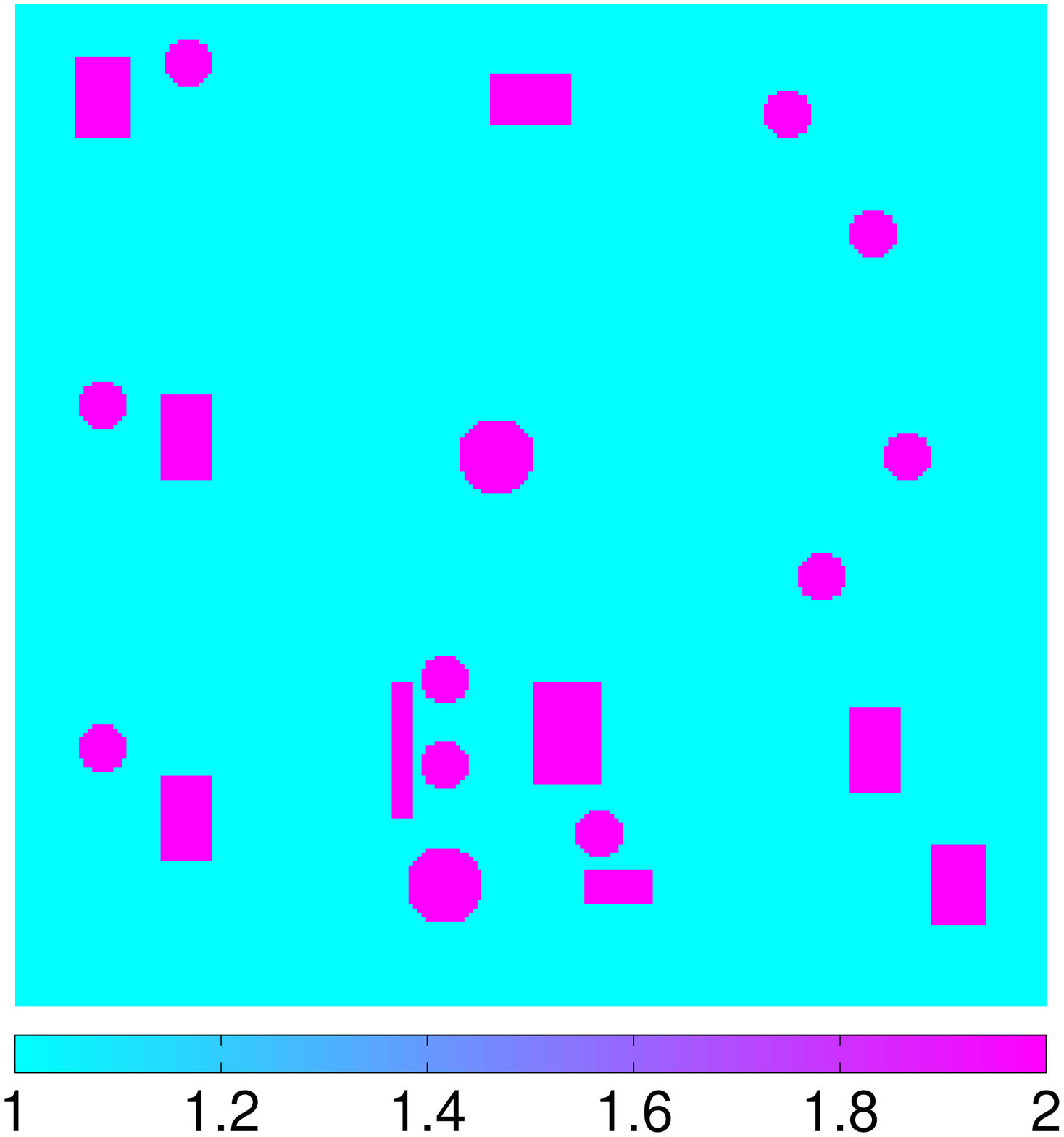}
     \label{ex4ltxi}
     }
      \subfigure[$\xi$ ($\alpha=0\%$)]{ 
    \includegraphics[width=36.5mm,height=35mm]{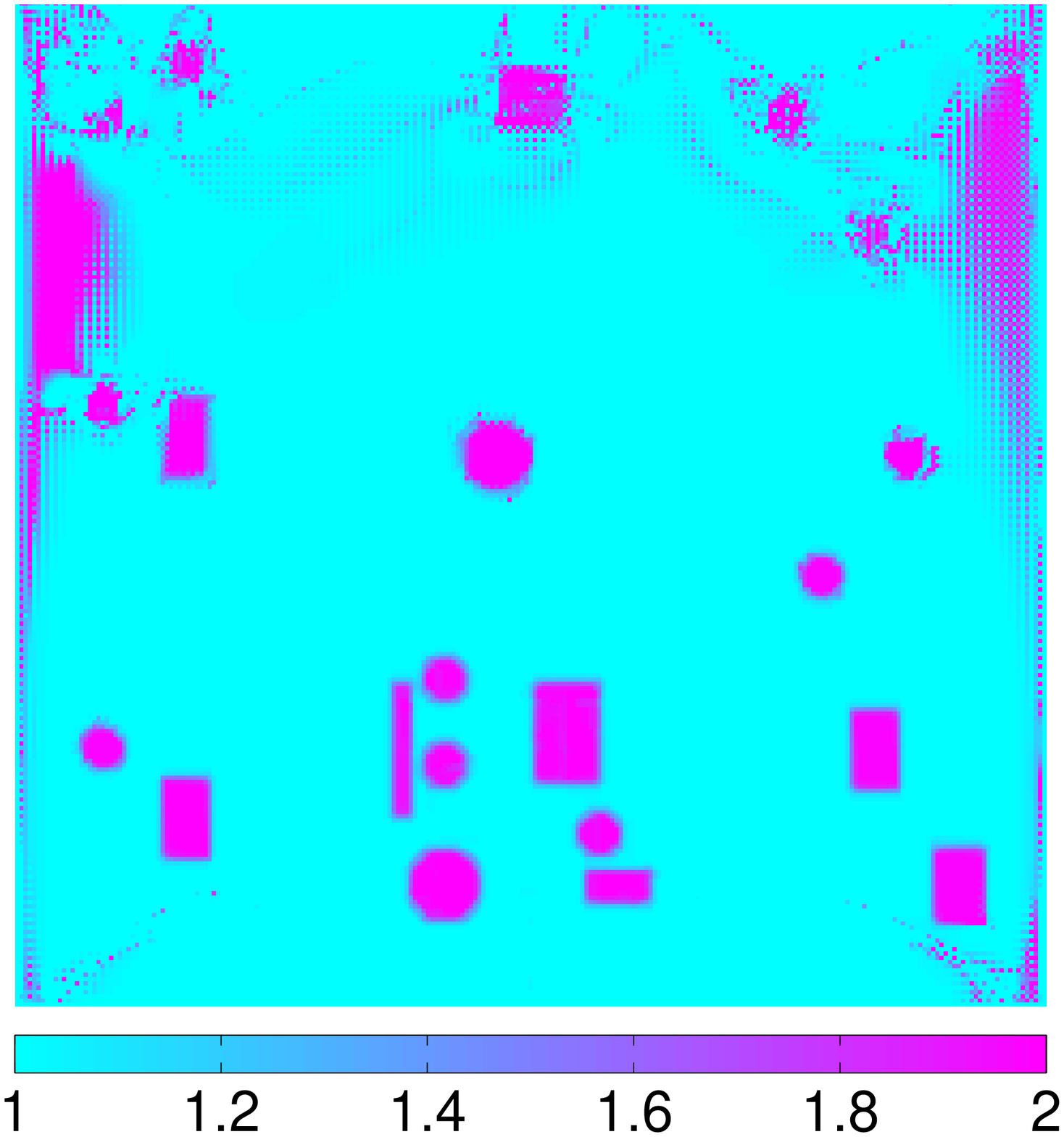}
    \label{ex4lcxi}
    }
    \subfigure[$\xi$ ($\alpha=4\%$)]{ 
      \includegraphics[width=36.5mm,height=35mm]{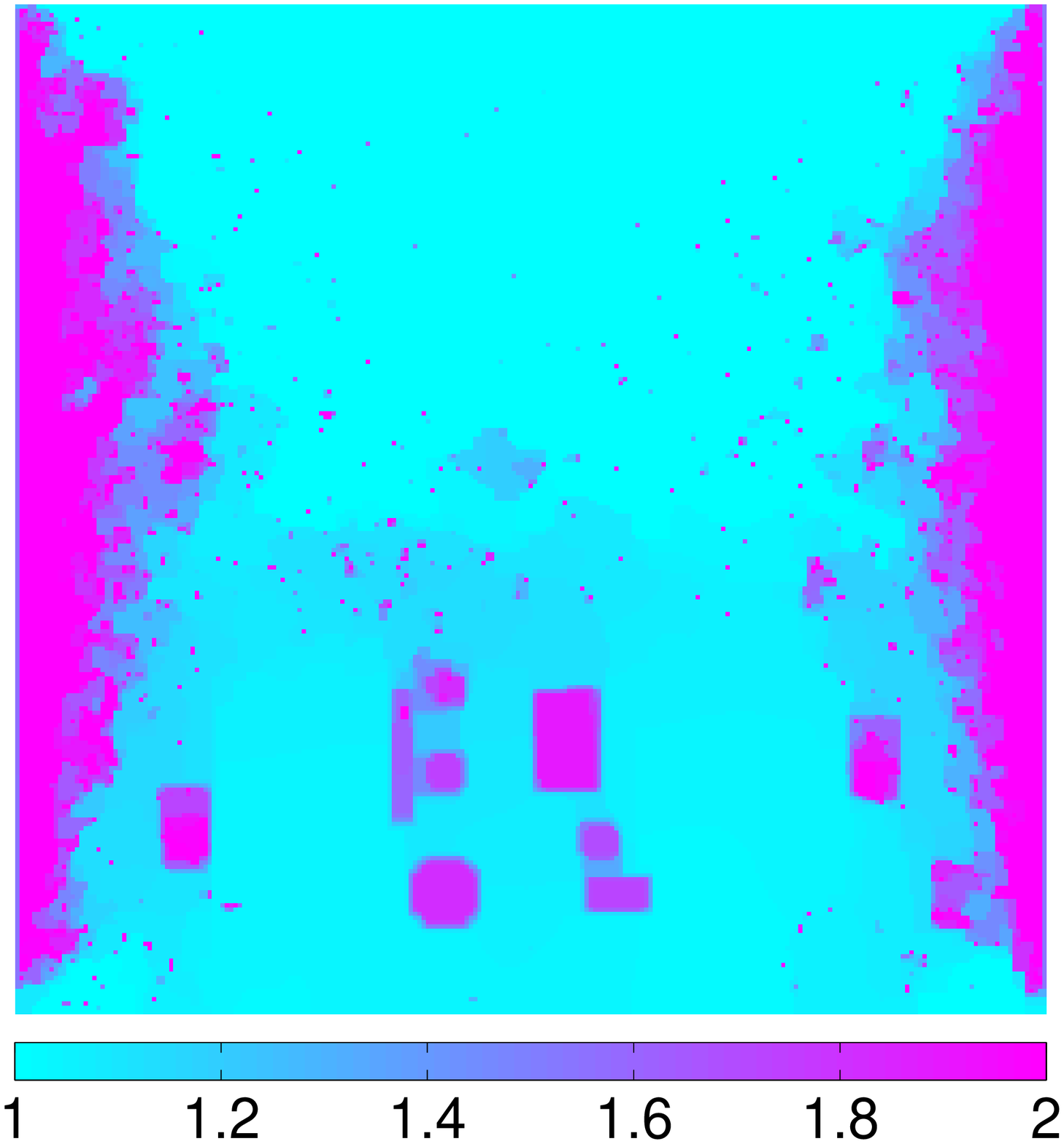} 
     \label{ex4lnxi}
     }
    
    \subfigure[true $\zeta$]{ 
      \includegraphics[width=36.5mm,height=35mm]{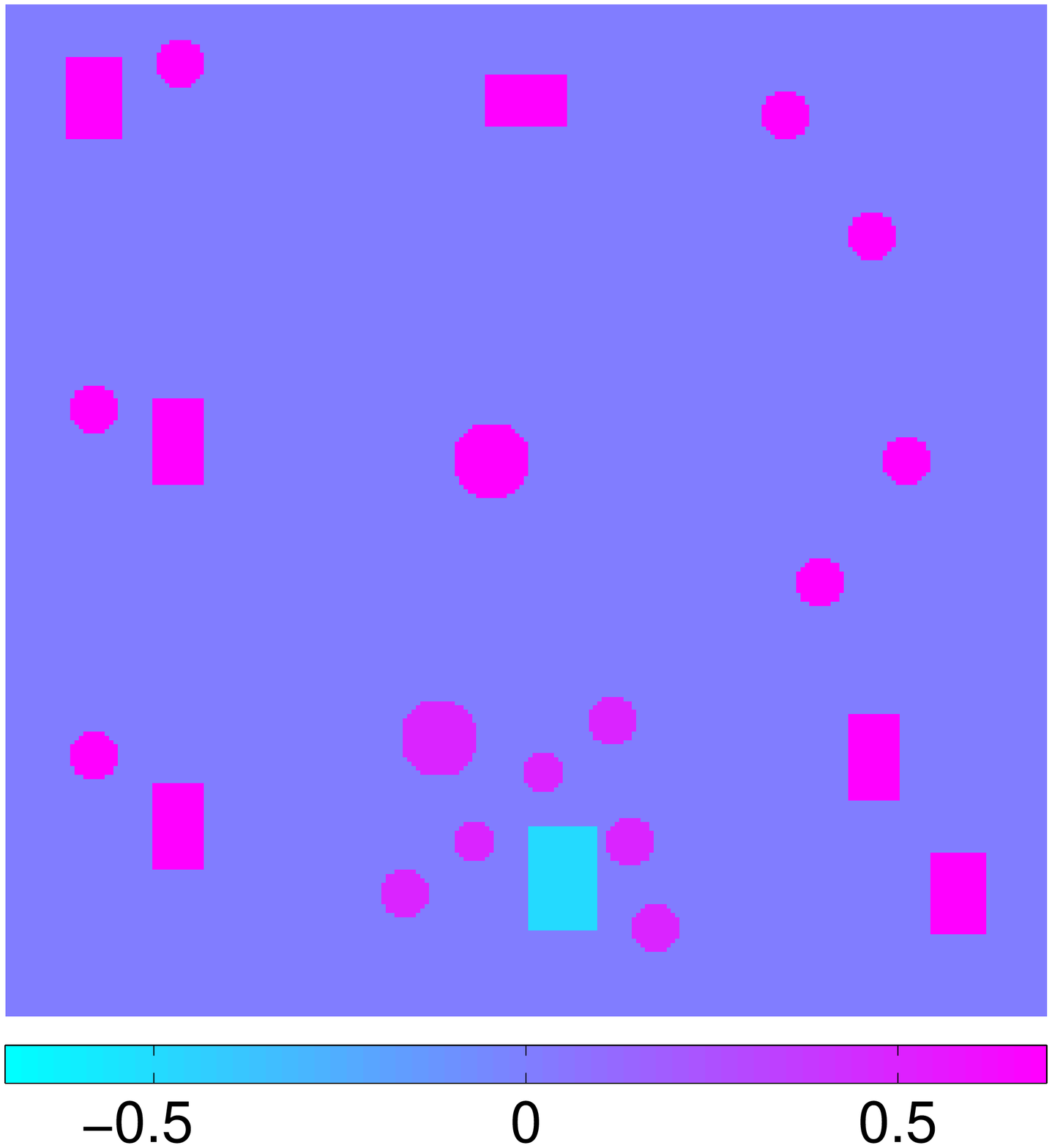} 
     \label{ex4lttau}
     }
     \subfigure[$\zeta$ ($\alpha=0\%$)]{
      \includegraphics[width=36.5mm,height=35mm]{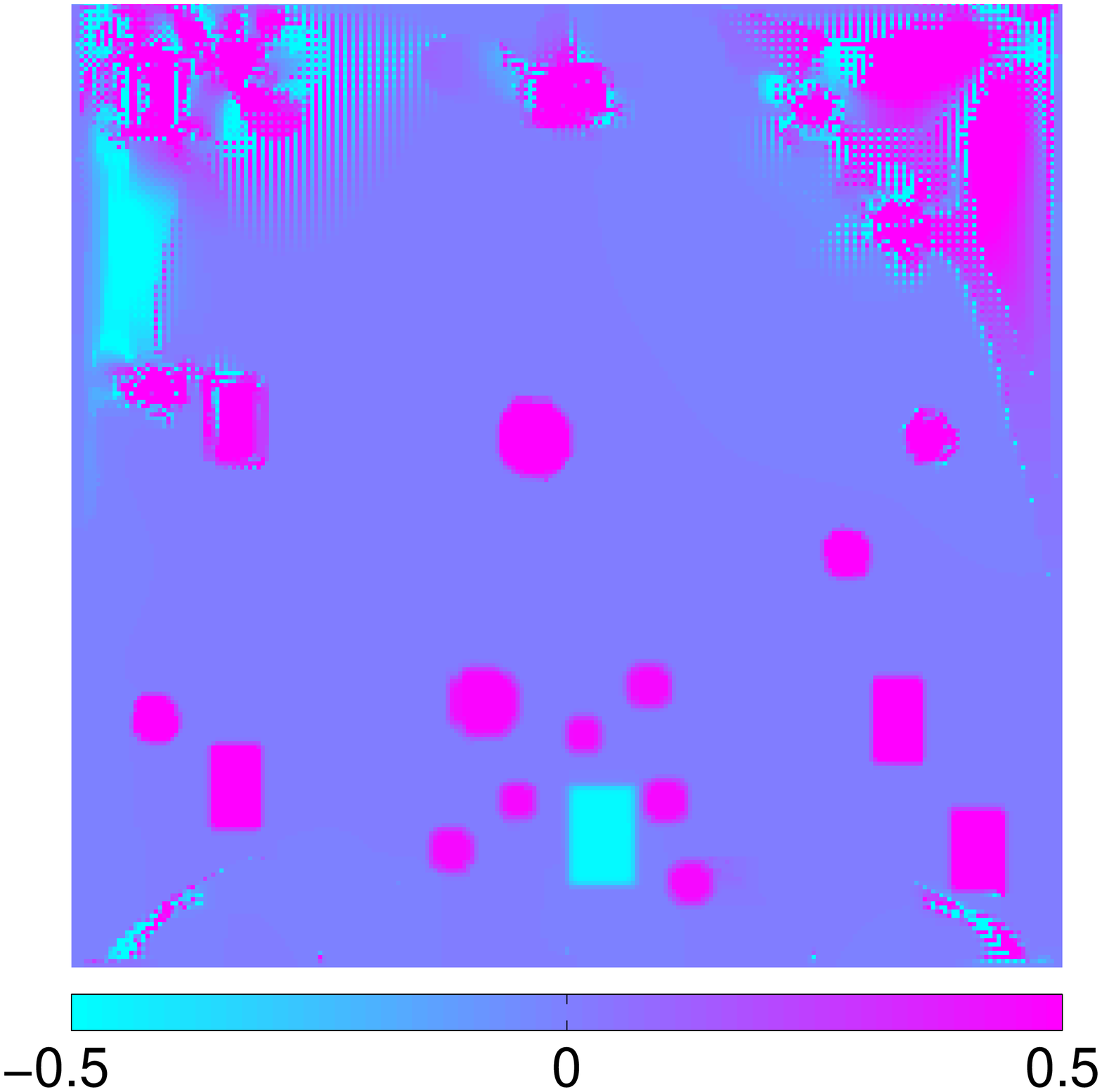}
     \label{ex4lctau}
     }     
     \subfigure[$\zeta$ ($\alpha=4\%$)]{ 
      \includegraphics[width=36.5mm,height=35mm]{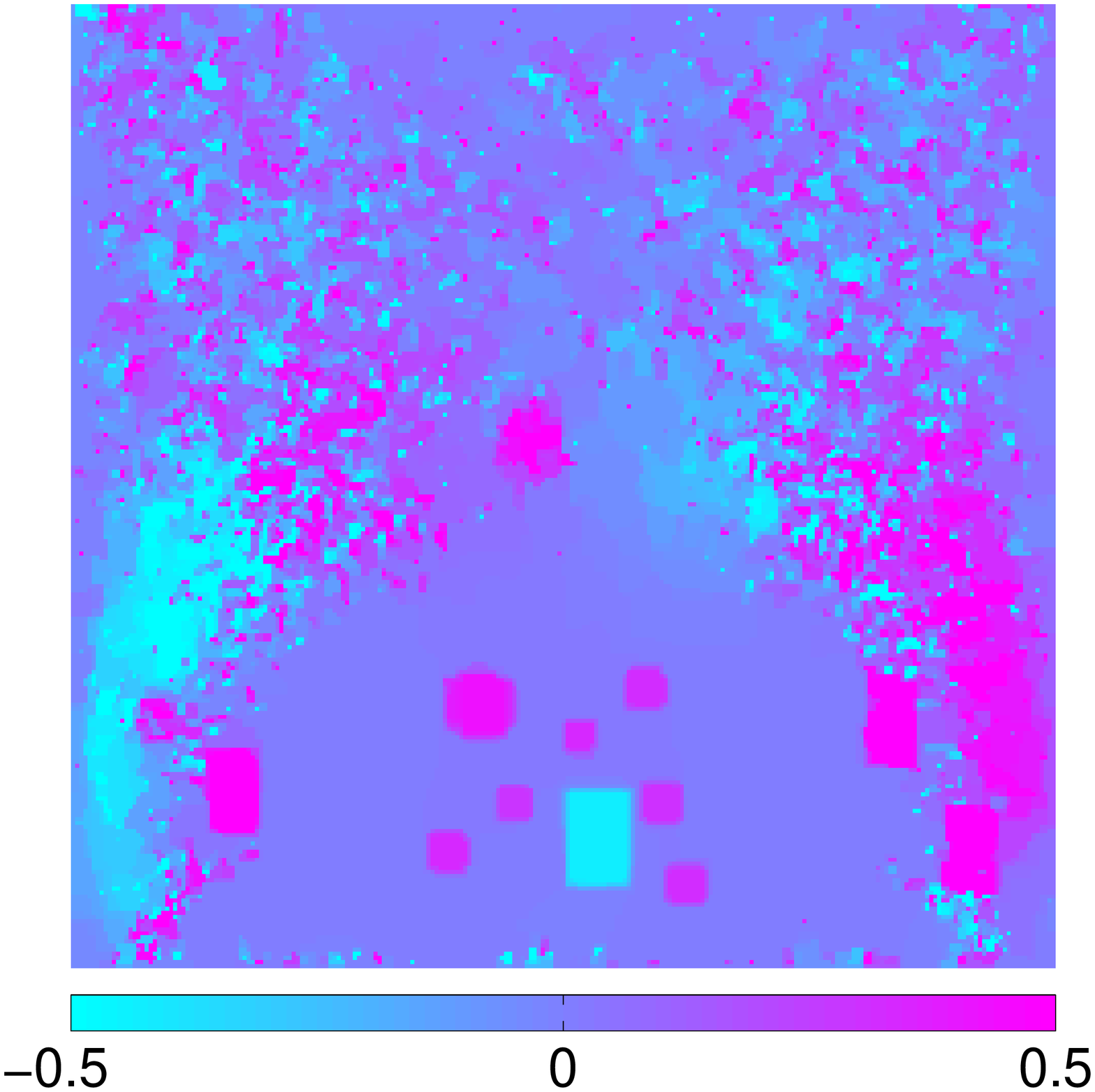} 
     \label{ex4lntau}
     }
      \caption{Simulations on extended domain $X'$. \subref{ex4ltxi}\&\subref{ex4lttau}: true anisotropy $(\xi,\zeta)$. \subref{ex4lcxi}\&\subref{ex4lctau}: reconstructions with noiseless data. \subref{ex4lnxi}\&\subref{ex4lntau}: reconstructions with noisy data($\alpha=4\%$). }
\label{E4large}
  \end{figure}
  
We then use the reconstructions restricted on $X$ to present the desired anisotropy. In the next step, $\beta$ can be recovered on $X$ by using the reconstructed $\tilde\gamma$ in the first step. Figure \ref{E4} displays the numerical results with noiseless data and noisy data($\alpha=1\%,4\%$). A $l_1$ regularization using the split Bregman iteration method is used for both the anisotropic and isotropic components in this simulation. The relative $L^2$ errors in the reconstructions are $\mathcal{E}^C_{\xi}=9.4\%$, $\mathcal{E}^C_{\zeta}=27.6\%$, $\mathcal{E}^C_{\beta}=7.2\%$; $\mathcal{E}^N_{\xi}=9.6\%$, $\mathcal{E}^N_{\zeta}=28.1\%$, $\mathcal{E}^N_{\beta}=7.6\%$ when $\alpha=1\%$; $\mathcal{E}^N_{\xi}=15.8\%$, $\mathcal{E}^N_{\zeta}=38.3\%$ $\mathcal{E}^N_{\beta}=13.7\%$ when $\alpha=4\%$.

\begin{figure}[htp]
  \centering
	\subfigure[$\xi$ ($\alpha=0\%$)]{ 
			\includegraphics[width=37mm,height=35mm]{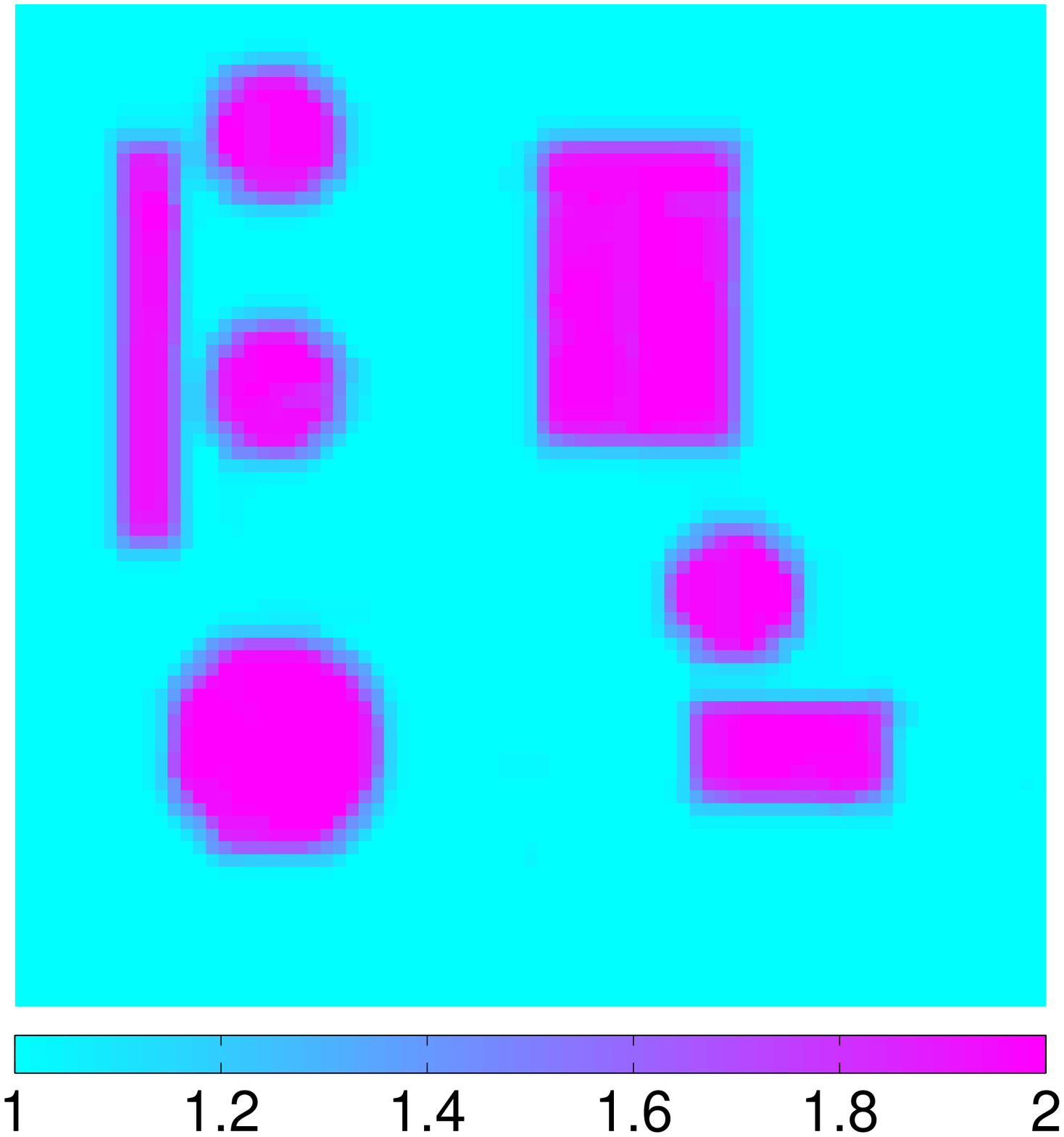}
			\label{ex4cxi}
			}
  \subfigure[$\xi$ ($\alpha=1\%$)]{ 
			\includegraphics[width=37mm,height=35mm]{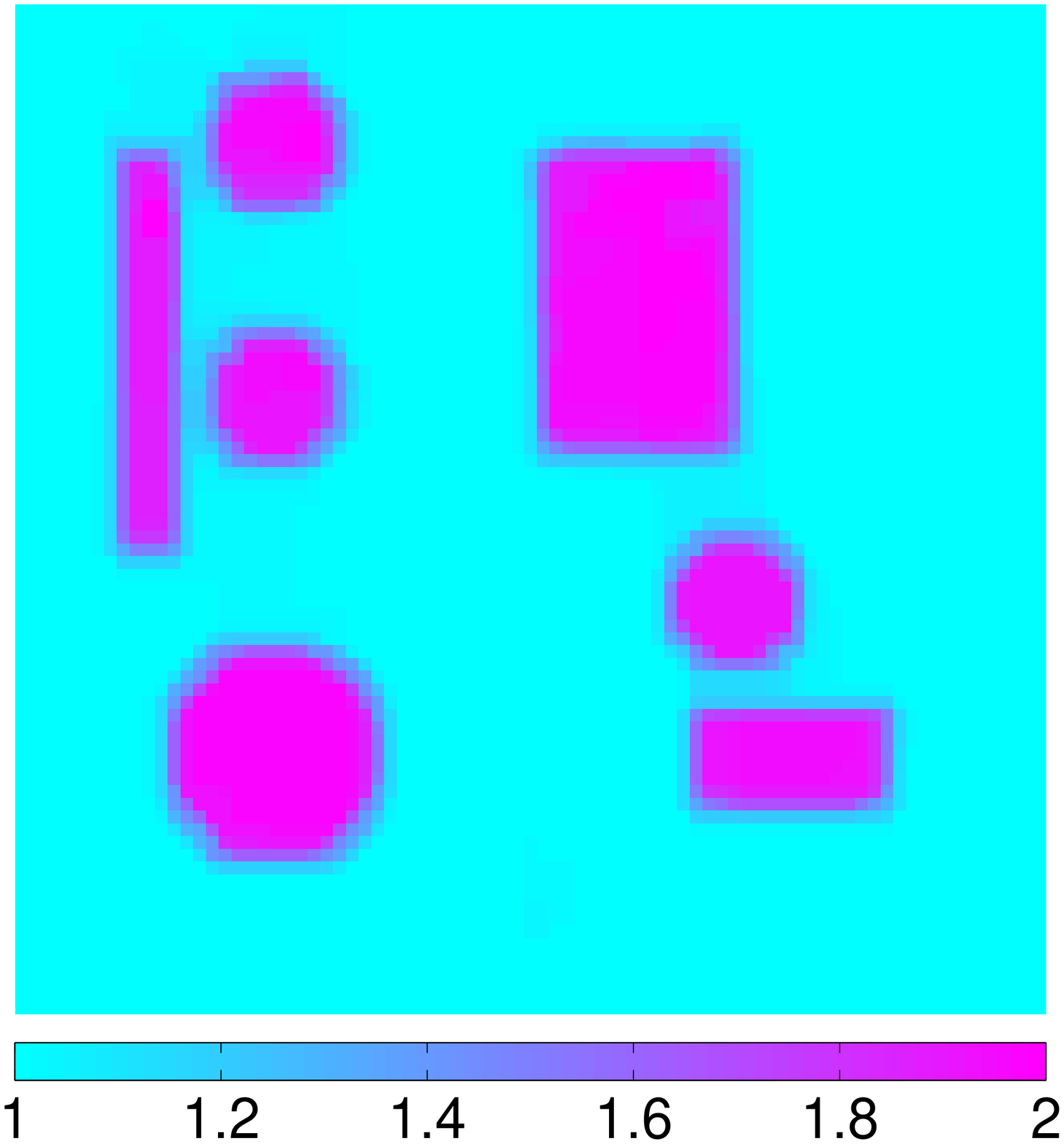}
      \label{ex4nxi}
      }
  \subfigure[$\xi$ ($\alpha=4\%$)]{ 
      \includegraphics[width=37mm,height=35mm]{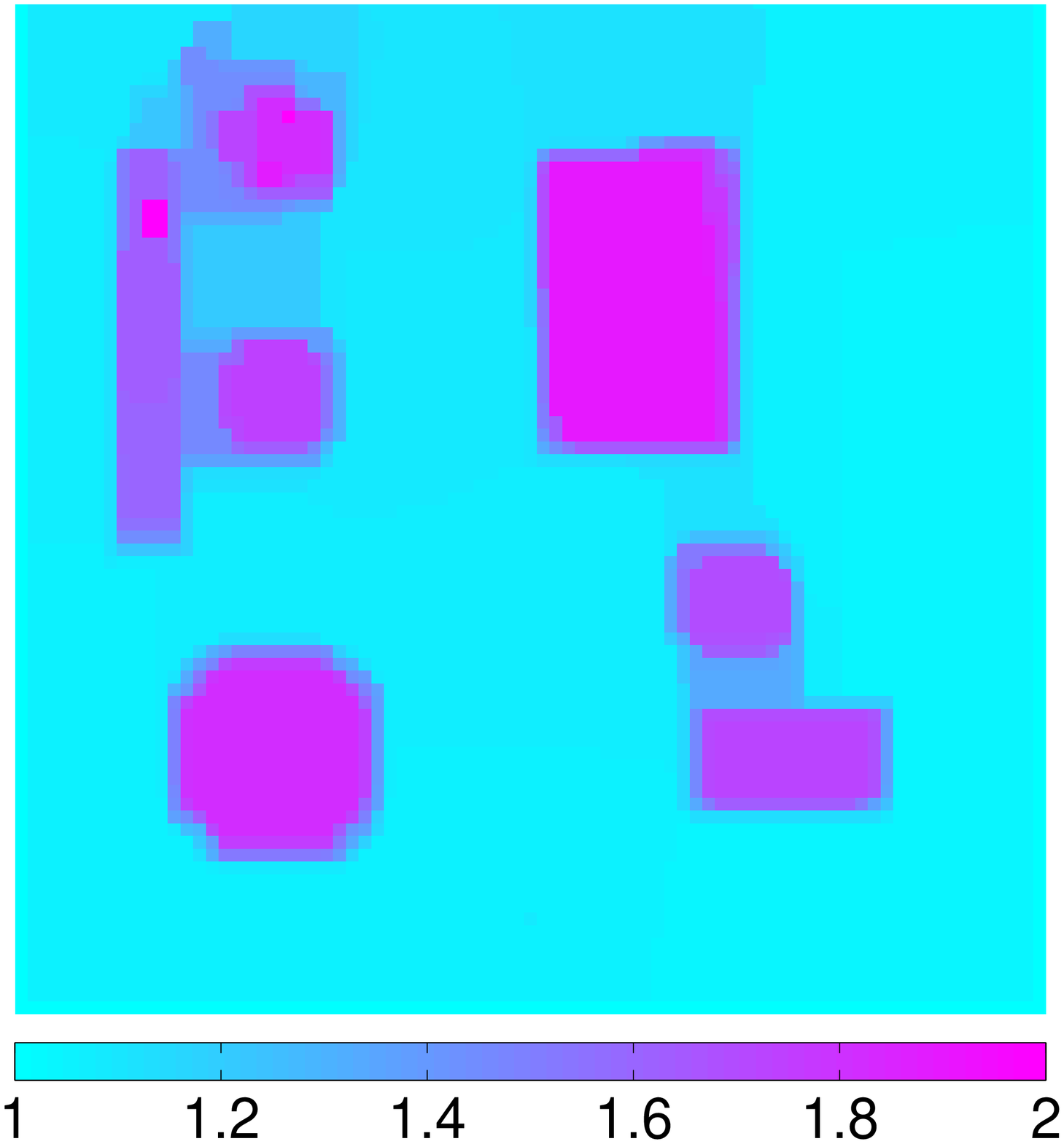}
      \label{ex4nxi4}
      }
  \subfigure[$\xi$ at $\{y=-0.5\}$]{ 
      \includegraphics[trim=10mm 5mm 10mm 0mm,clip,width=35mm,height=38mm]{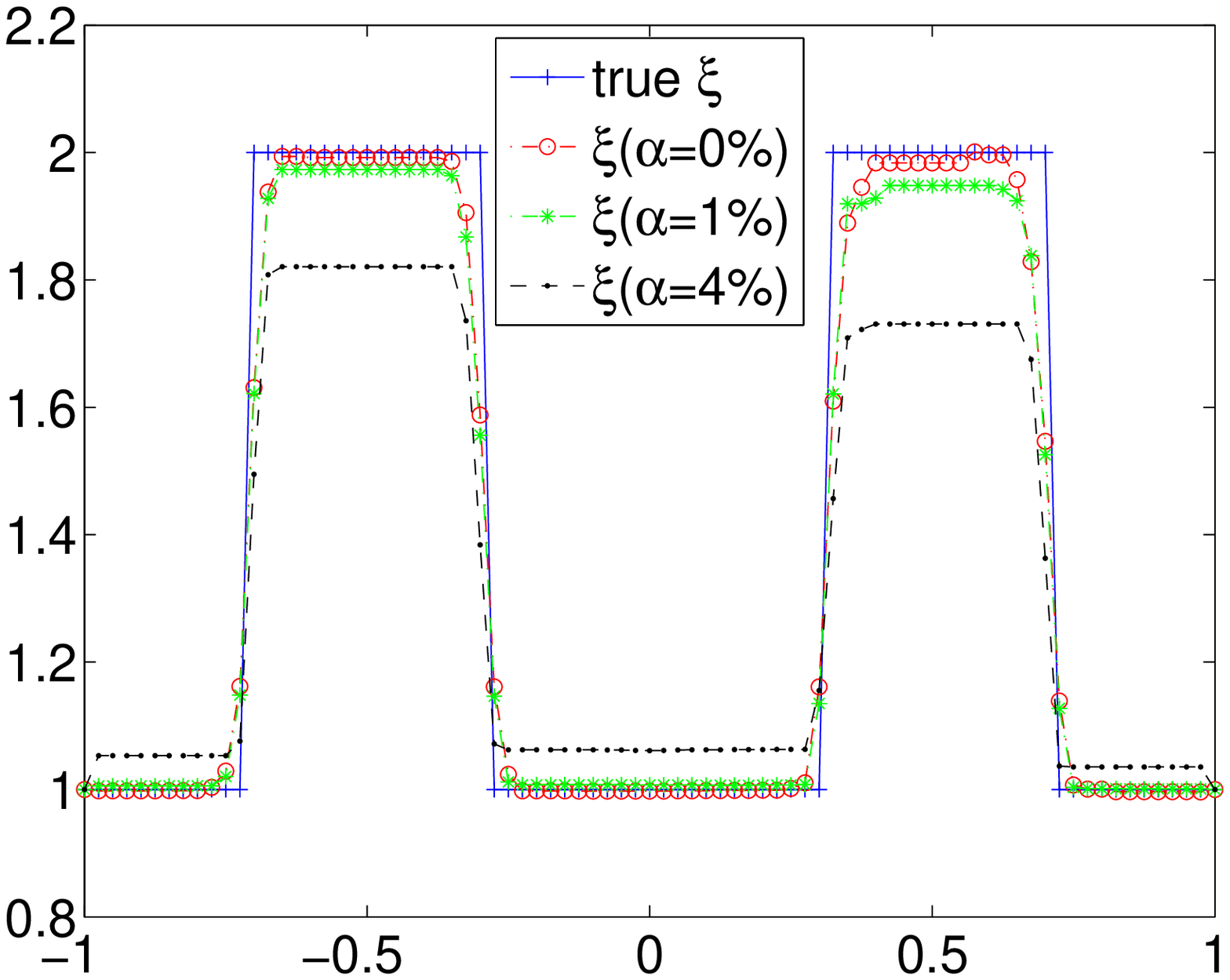} 
      \label{ex4rxi}
      }    
  \subfigure[$\zeta$ ($\alpha=0\%$)]{ 
      \includegraphics[width=37mm,height=35mm]{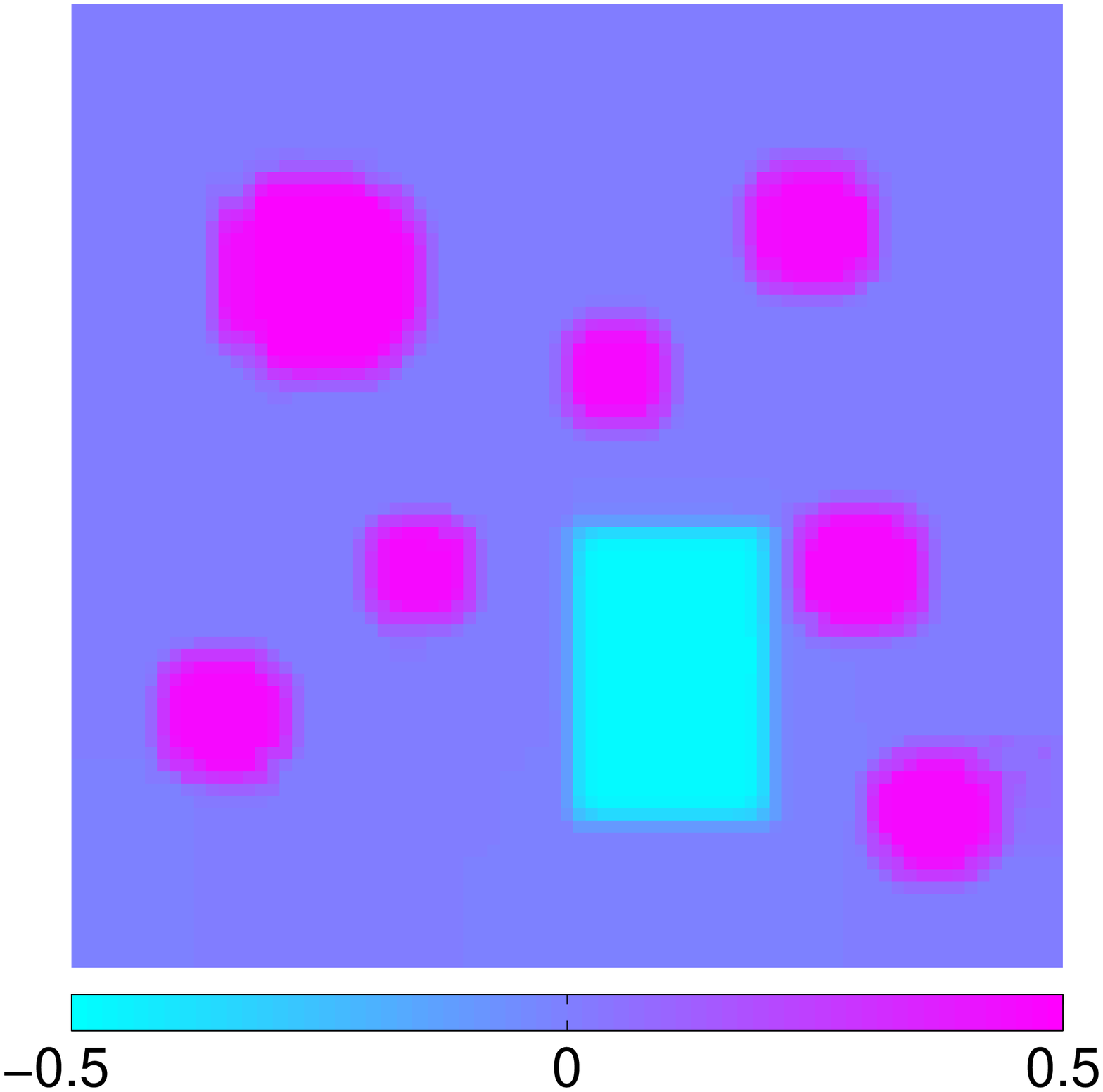}
      \label{ex4ctau}
      }
  \subfigure[$\zeta$ ($\alpha=1\%$)]{ 
      \includegraphics[width=37mm,height=35mm]{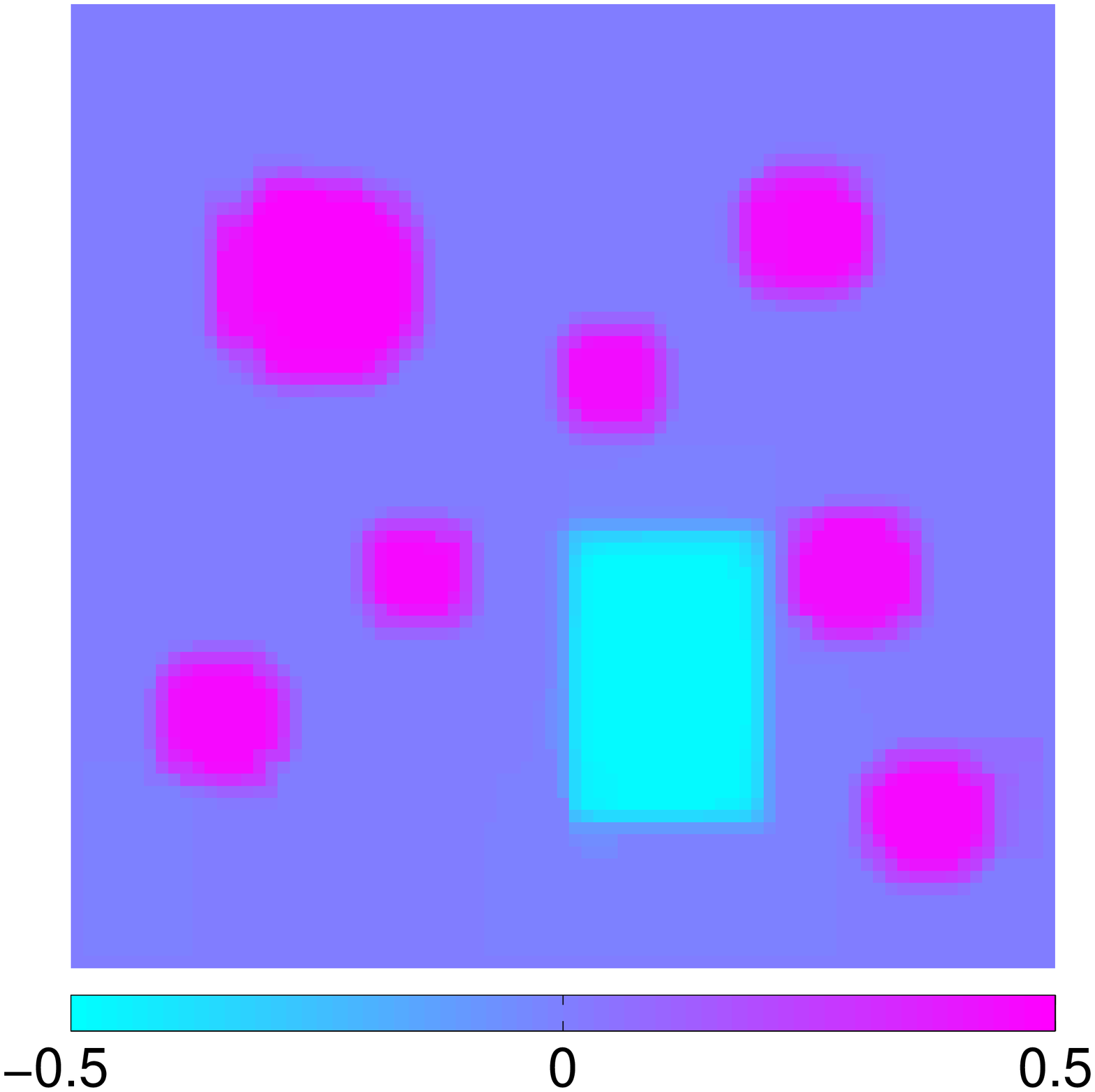}
      \label{ex4ntau}
      }
  \subfigure[$\zeta$ ($\alpha=4\%$)]{ 
      \includegraphics[width=37mm,height=35mm]{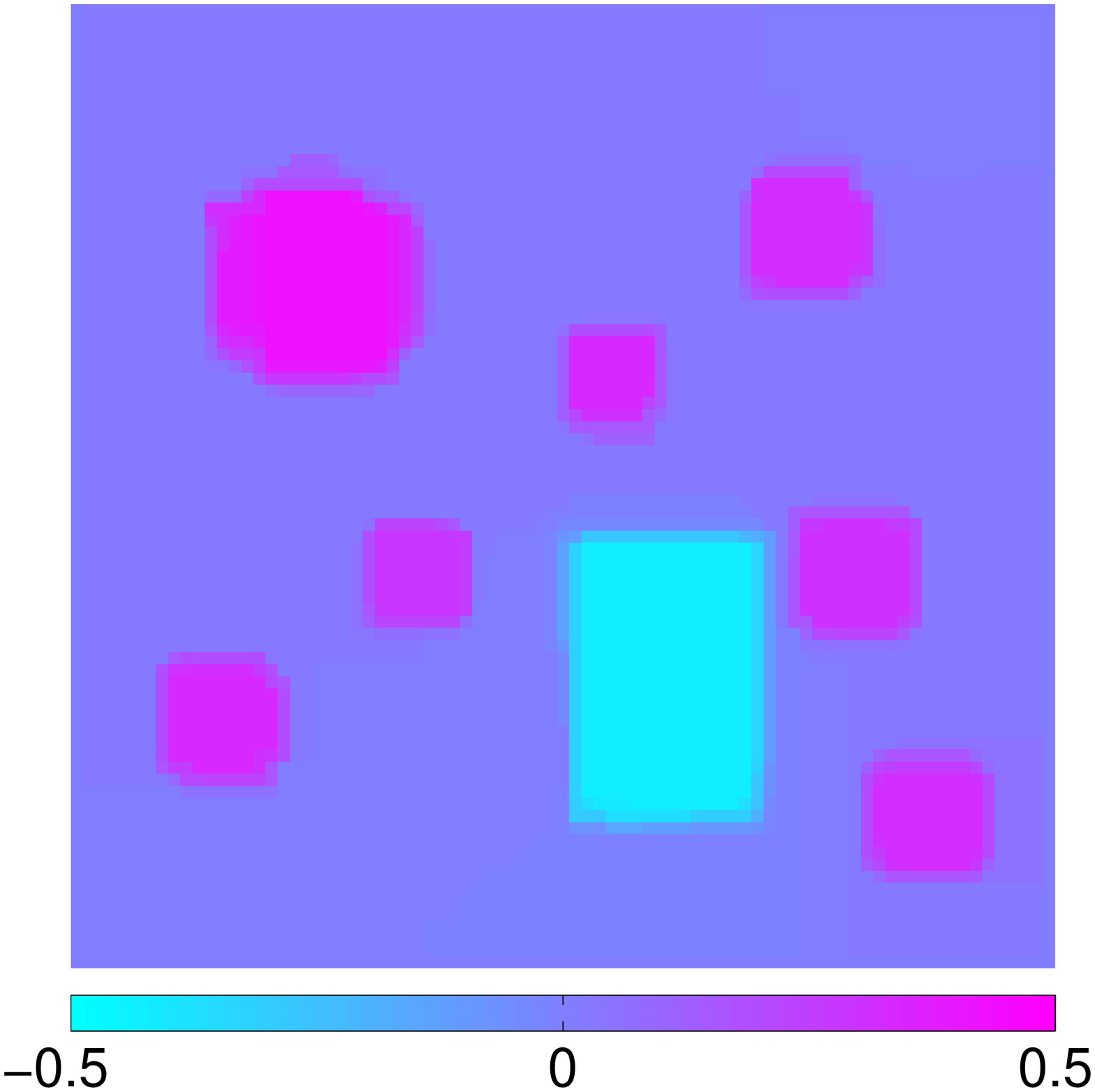}
      \label{ex4ntau4}
      }
  \subfigure[$\zeta$ at $\{y=-0.5\}$]{ 
      \includegraphics[trim=10mm 5mm 10mm 0mm,clip,width=35mm,height=38mm]{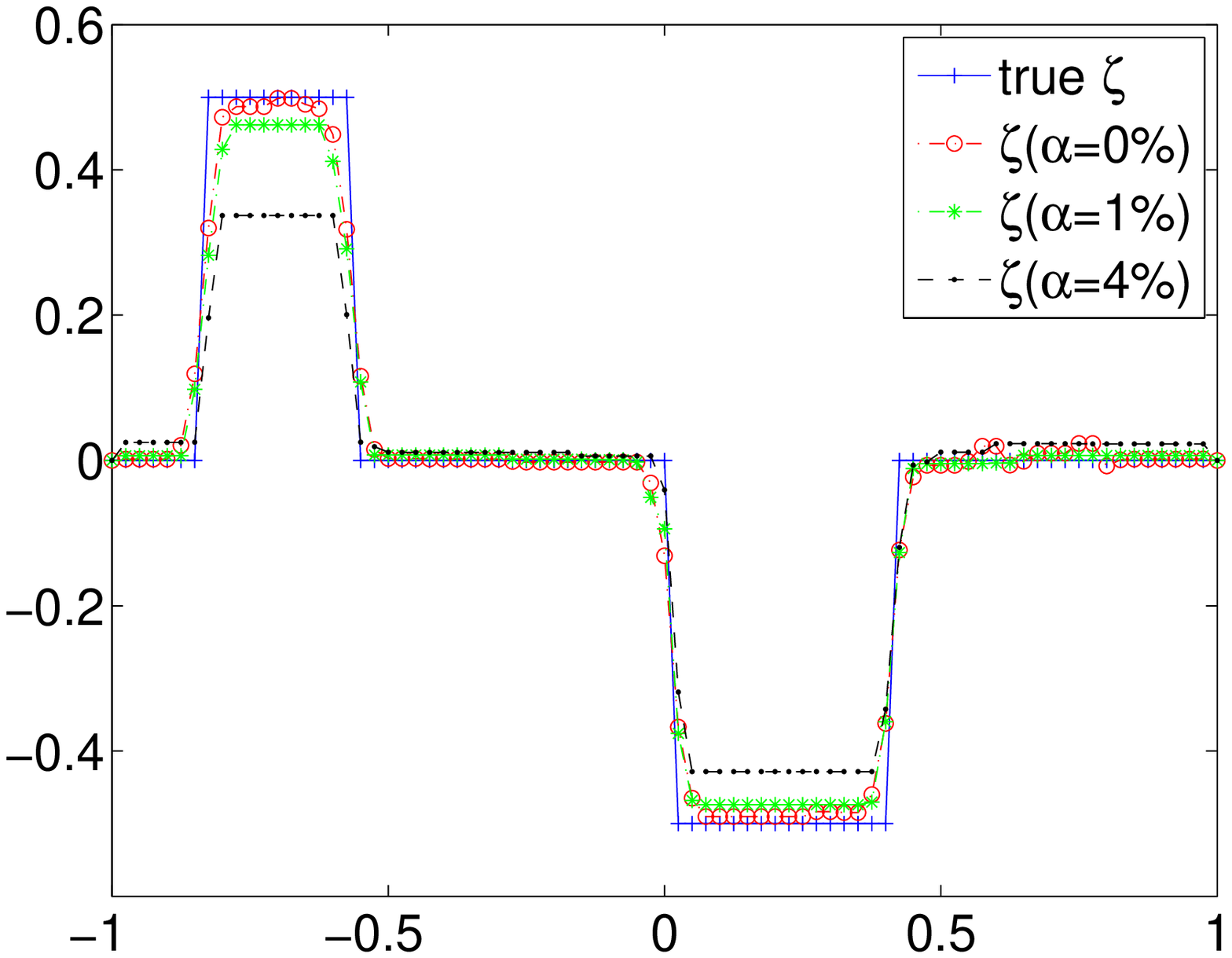} 
      \label{ex4rtau}
      }
  \subfigure[$|\gamma|^{\frac{1}{2}}$ ($\alpha=0\%$)]{
      \includegraphics[width=37mm,height=35mm]{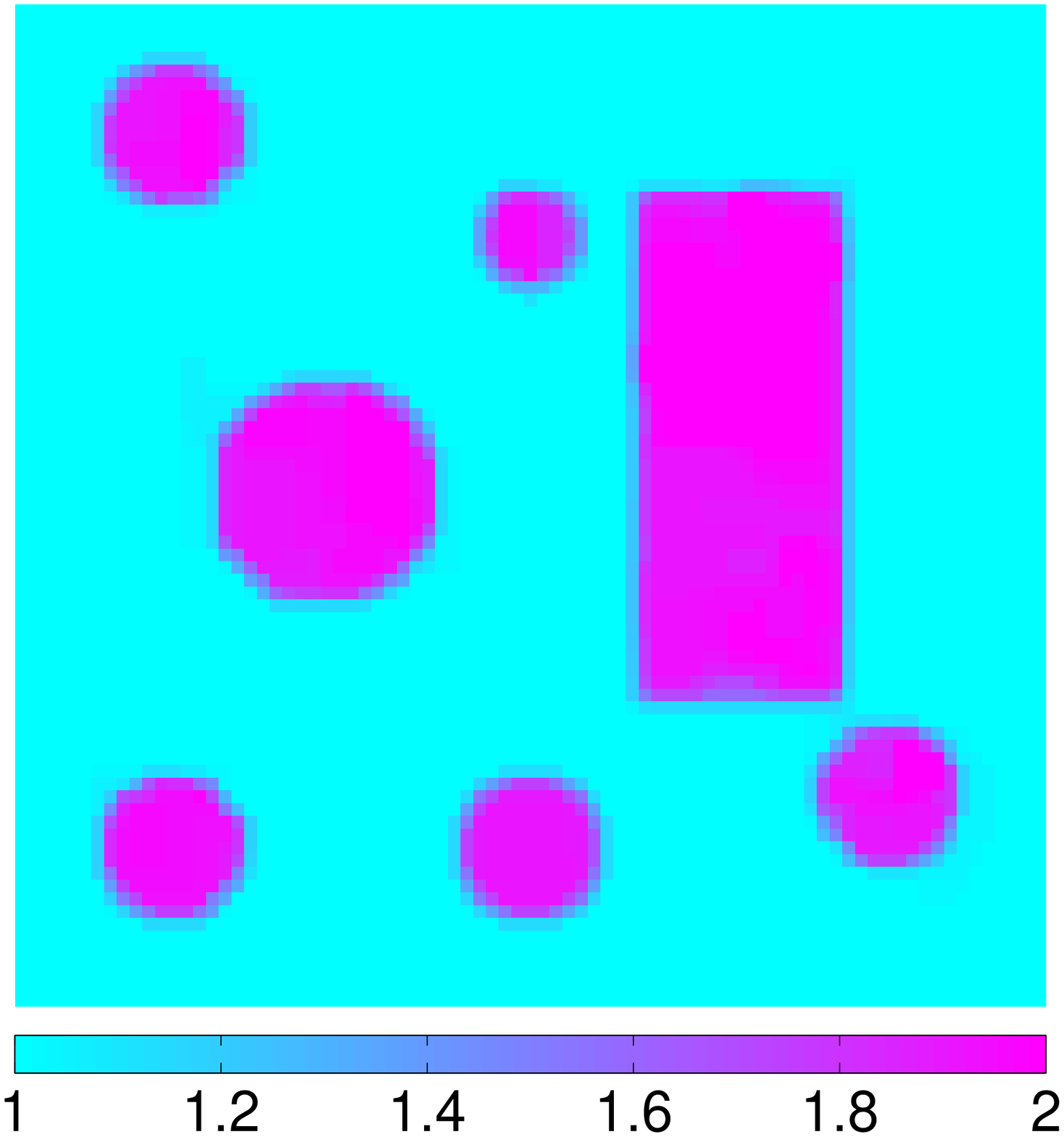}
      \label{ex4cbeta}
      }
  \subfigure[$|\gamma|^{\frac{1}{2}}$ ($\alpha=1\%$)]{
      \includegraphics[width=37mm,height=35mm]{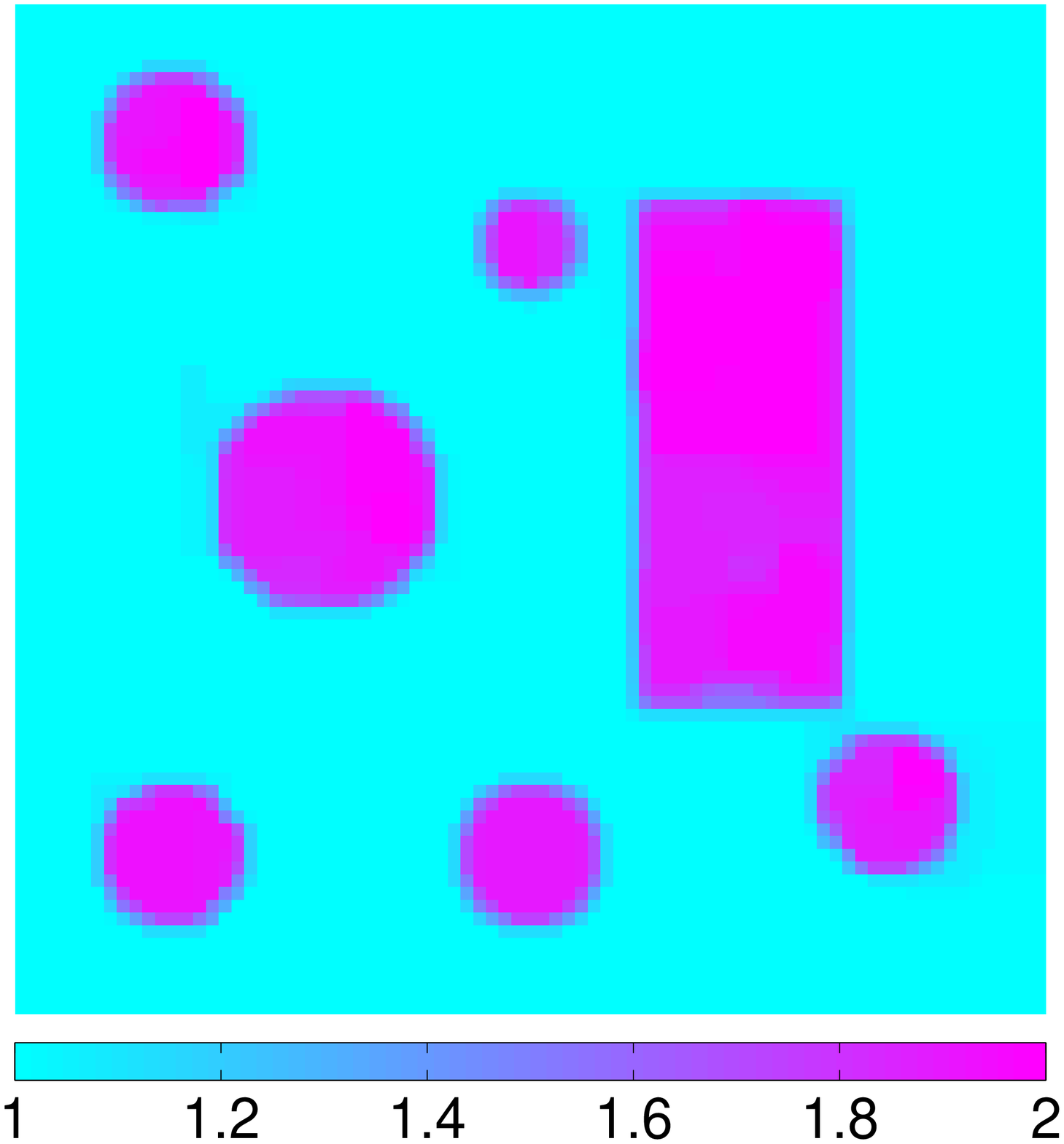}
      \label{ex4nbeta}
      }
  \subfigure[$|\gamma|^{\frac{1}{2}}$ ($\alpha=4\%$)]{
      \includegraphics[width=37mm,height=35mm]{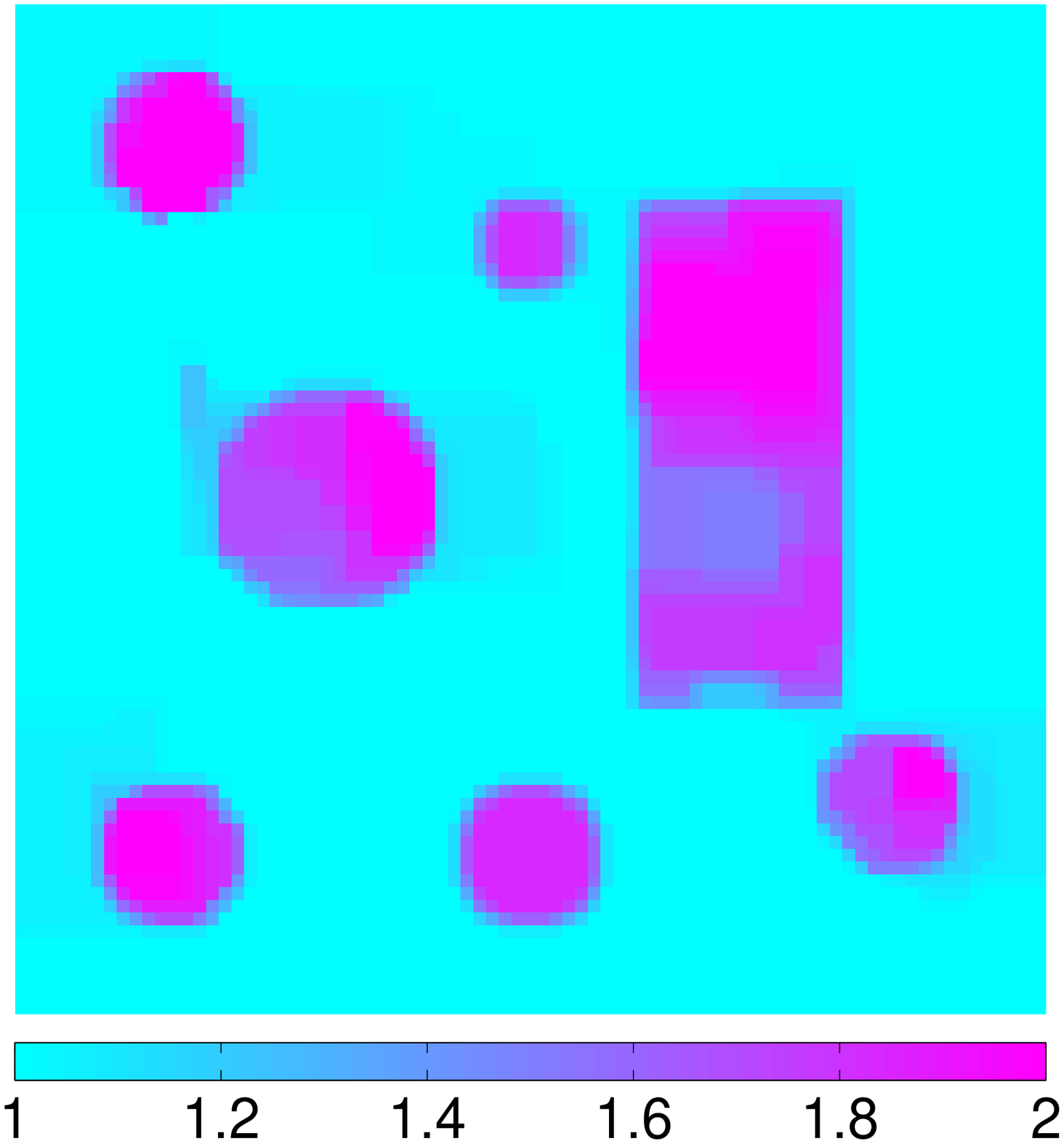}
      \label{ex4nbeta4}
      }
  \subfigure[$|\gamma|^{\frac{1}{2}}$ at $\{y=-0.5\}$]{ 
      \includegraphics[trim=15mm 5mm 14mm 0mm,clip,width=35mm,height=38mm]{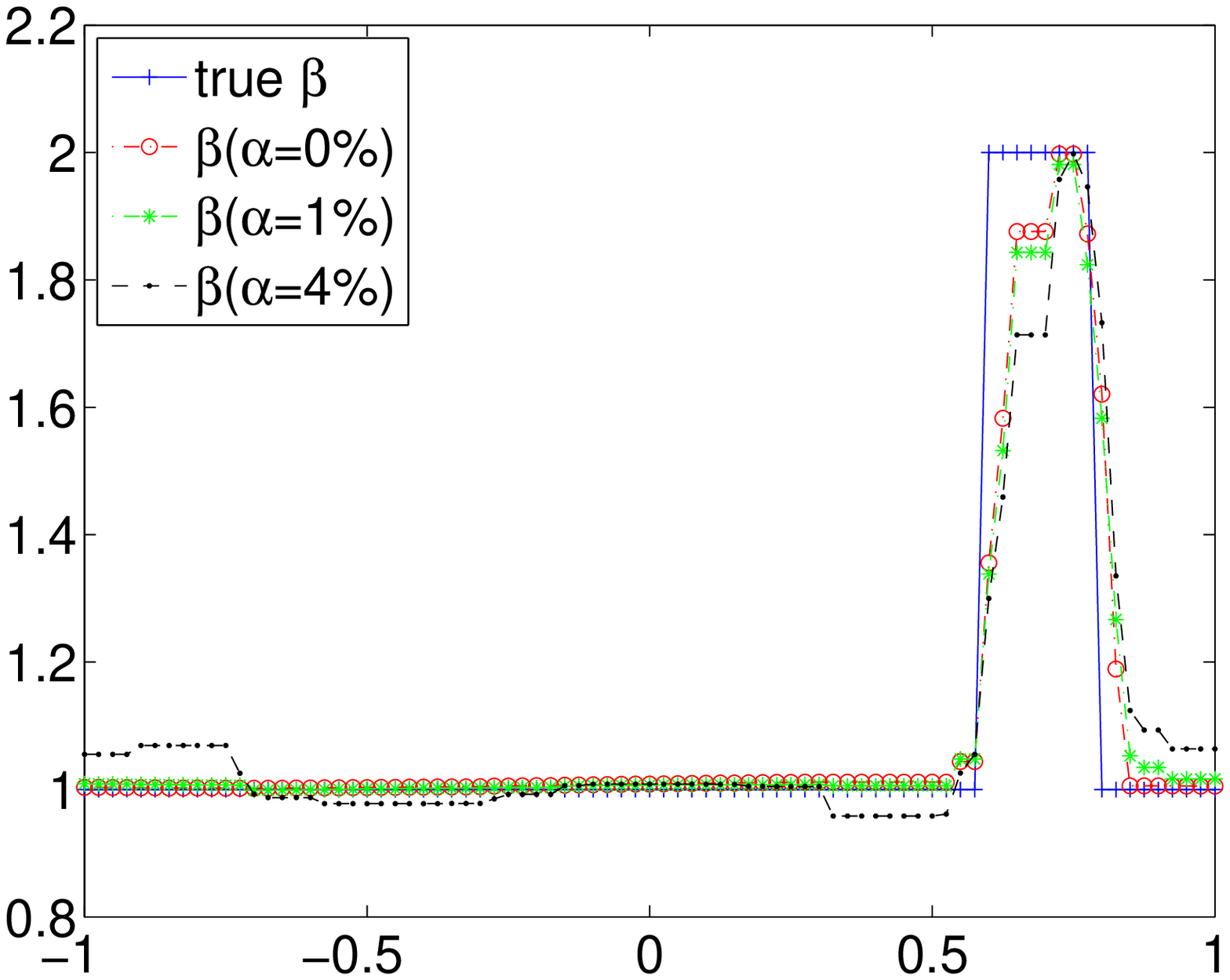} 
      \label{ex4rbeta}
      }
      \caption{Experiment 4. \subref{ex4cxi}\&\subref{ex4ctau}\&\subref{ex4cbeta}: reconstructions with noiseless data. \subref{ex4nxi}\&\subref{ex4ntau}\&\subref{ex4nbeta}: reconstructions with noisy data ($\alpha=1\%$). \subref{ex4nxi4}\&\subref{ex4ntau4}\&\subref{ex4nbeta4}: reconstructions with noisy data ($\alpha=4\%$). \subref{ex4rxi}\&\subref{ex4rtau}\&\subref{ex4rbeta}: cross sections along $y=-0.5$.}
\label{E4}
\end{figure}
  
\paragraph{Experiment 5.} In this experiment, we repeat Experiment 4 on the extended domain $X'$, replacing homogeneous Dirichlet boundary conditions on the left, top and right edges, by homogeneous Neumann conditions on the three other edges. The same (Dirichlet) boundary conditions are used on the bottom edge of the domain,
\begin{align}
\left\{\begin{array}{ll}
u(\mathbf{x})=(2\pi\cdot0.2^2)^{-\frac{1}{2}}\exp\{-\frac{1}{2\cdot0.2^2}(x+x_i)^2\}, \enskip &\mathbf{x}\in \partial X'_B\\
\frac{\partial u}{\partial n}(\mathbf{x})=0, \enskip &\mathbf{x}\in \partial X'\setminus\partial X'_B
\end{array}\right. \quad 1\leq i\leq 5
\end{align}
where $\{\mathbf{x}_i\}_{1\leq i\leq 5}=\{-2.8,-1.5,0,1.5,2.8\}$. As in the last experiment, we first apply the reconstruction algorithm of $\tilde\gamma$ on $X'$ and present its restriction on $X$. Then $\beta$ can be recovered on $X$ by using the reconstructed $\tilde\gamma$. Figure \ref{E5Neumann} displays the numerical results with noiseless data and noisy data ($\alpha=1\%,4\%$). An $l_1$ regularization procedure is again used in this simulation. The relative $L^2$ errors in the reconstructions are $\mathcal{E}^C_{\xi}=9.4\%$, $\mathcal{E}^C_{\zeta}=26.9\%$, $\mathcal{E}^C_{\beta}=6.8\%$; $\mathcal{E}^N_{\xi}=9.5\%$, $\mathcal{E}^N_{\zeta}=28.7\%$, $\mathcal{E}^N_{\beta}=7.7\%$ when $\alpha=1\%$; $\mathcal{E}^N_{\xi}=14.3\%$, $\mathcal{E}^N_{\zeta}=52.1\%$, $\mathcal{E}^N_{\beta}=13.5\%$ when $\alpha=4\%$.
  
\begin{figure}[htp]
  \centering
    \subfigure[$\xi$ ($\alpha=0\%$)]{ 
     \includegraphics[width=37mm,height=35mm]{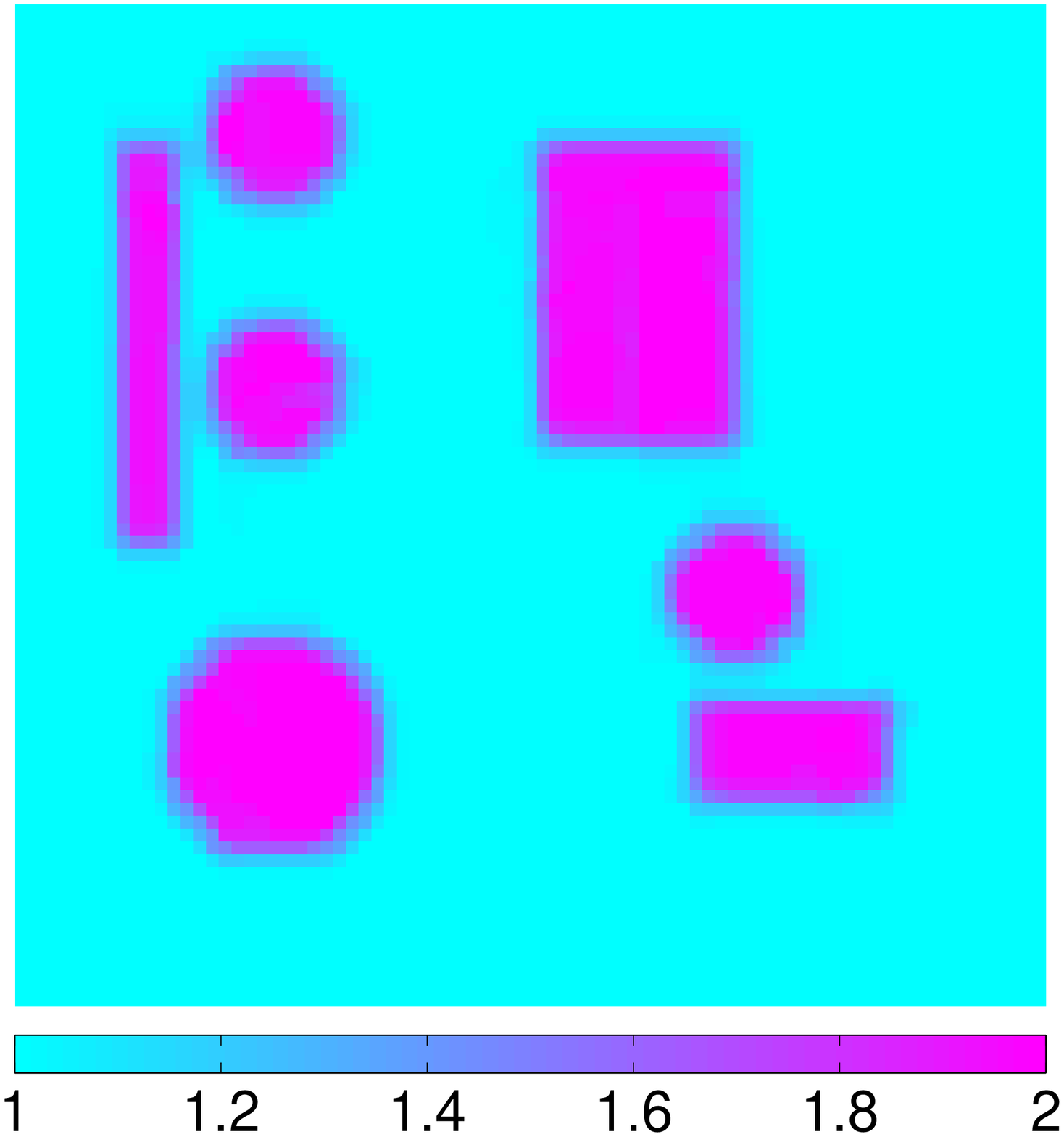}
     \label{ex5cxi}
     }
     \subfigure[$\xi$ ($\alpha=1\%$)]{ 
    \includegraphics[width=37mm,height=35mm]{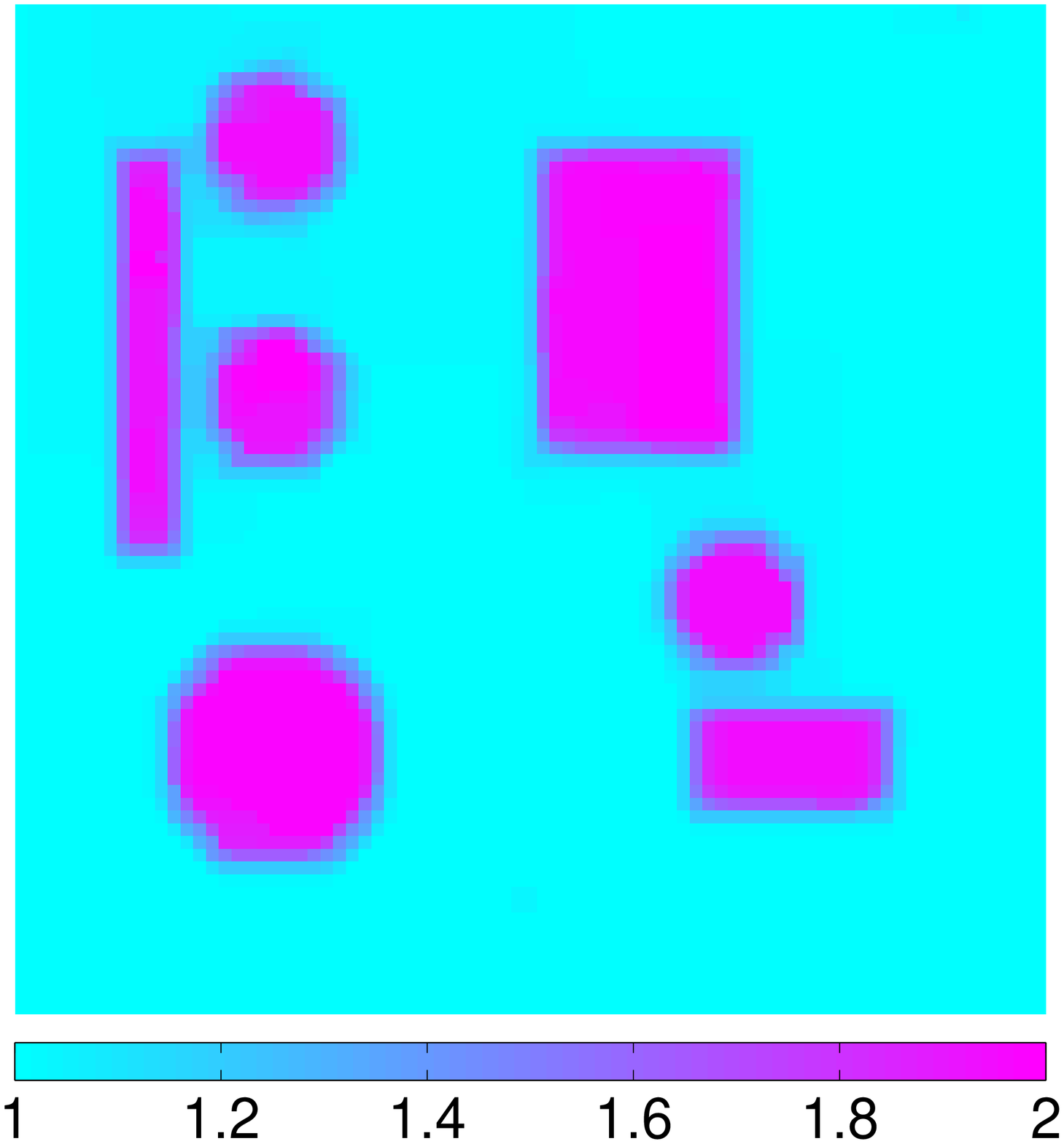}
    \label{ex5nxi}
    }
    \subfigure[$\xi$ ($\alpha=4\%$)]{ 
    \includegraphics[width=37mm,height=35mm]{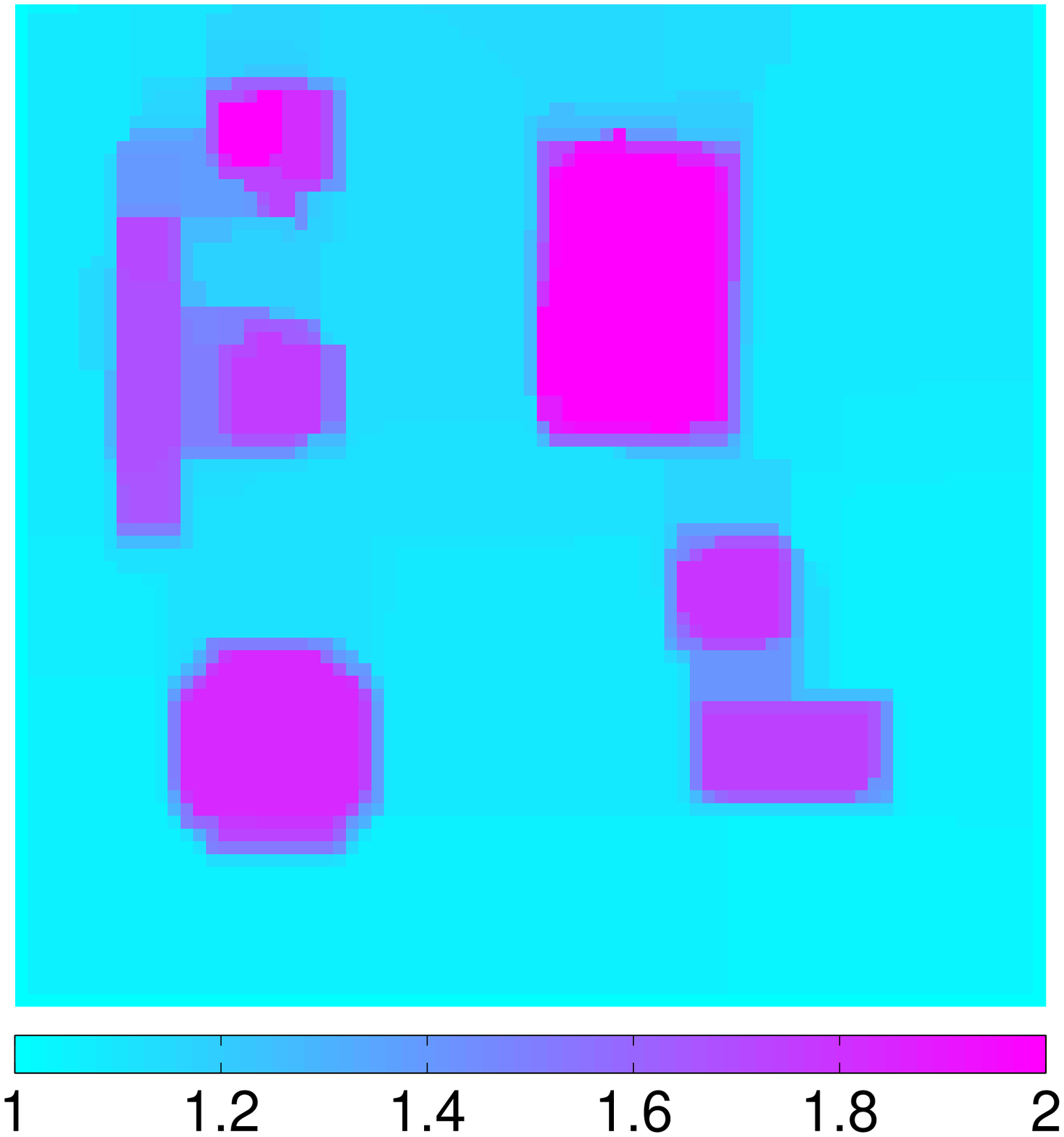}
    \label{ex5nxi4}
    }
    \subfigure[$\xi$ at $\{y=-0.5\}$]{ 
     \includegraphics[trim=10mm 5mm 10mm 0mm,clip,width=35mm,height=38mm]{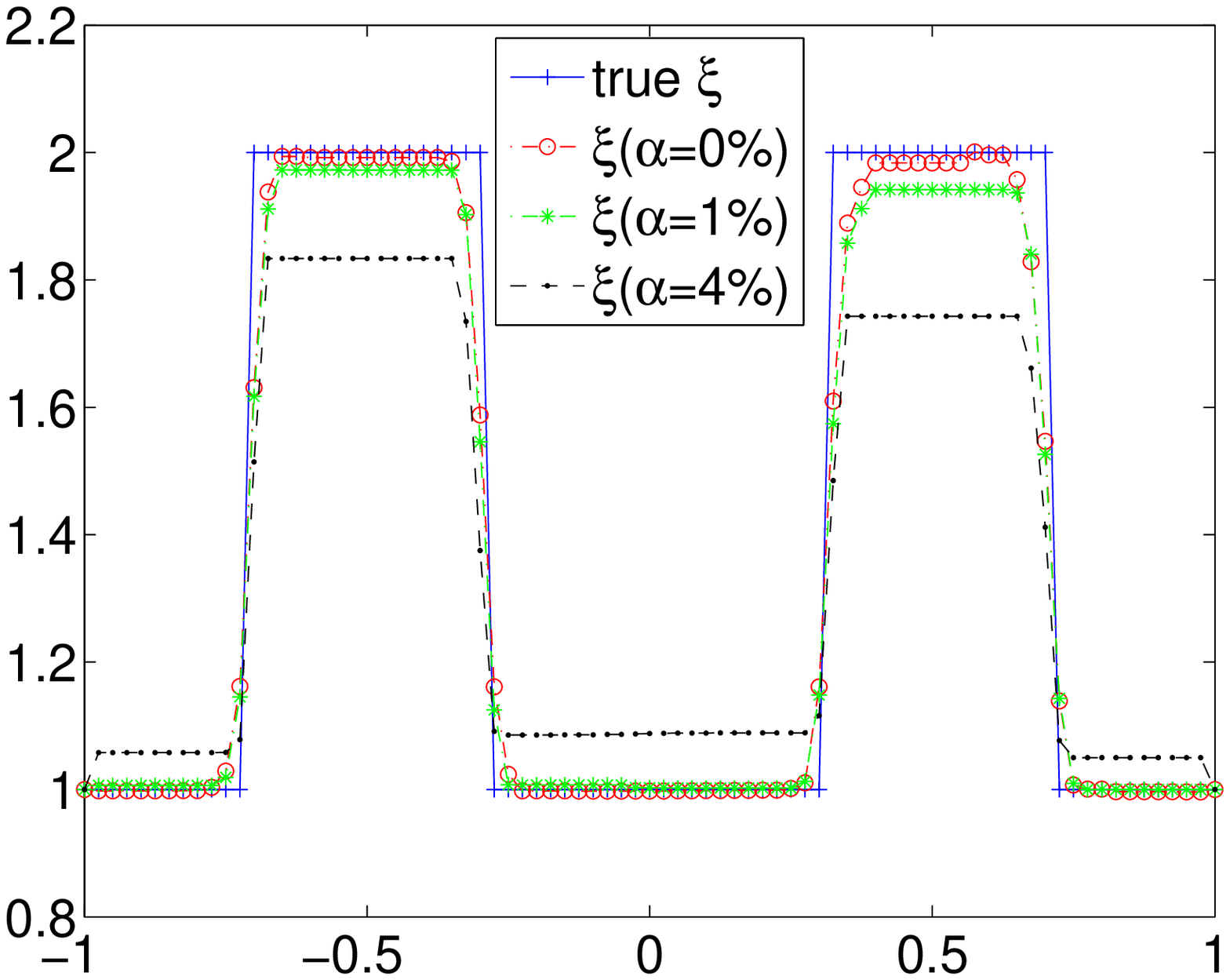} 
      \label{ex5rxi}
     }
      \subfigure[$\zeta$ ($\alpha=0\%$)]{ 
     \includegraphics[width=37mm,height=35mm]{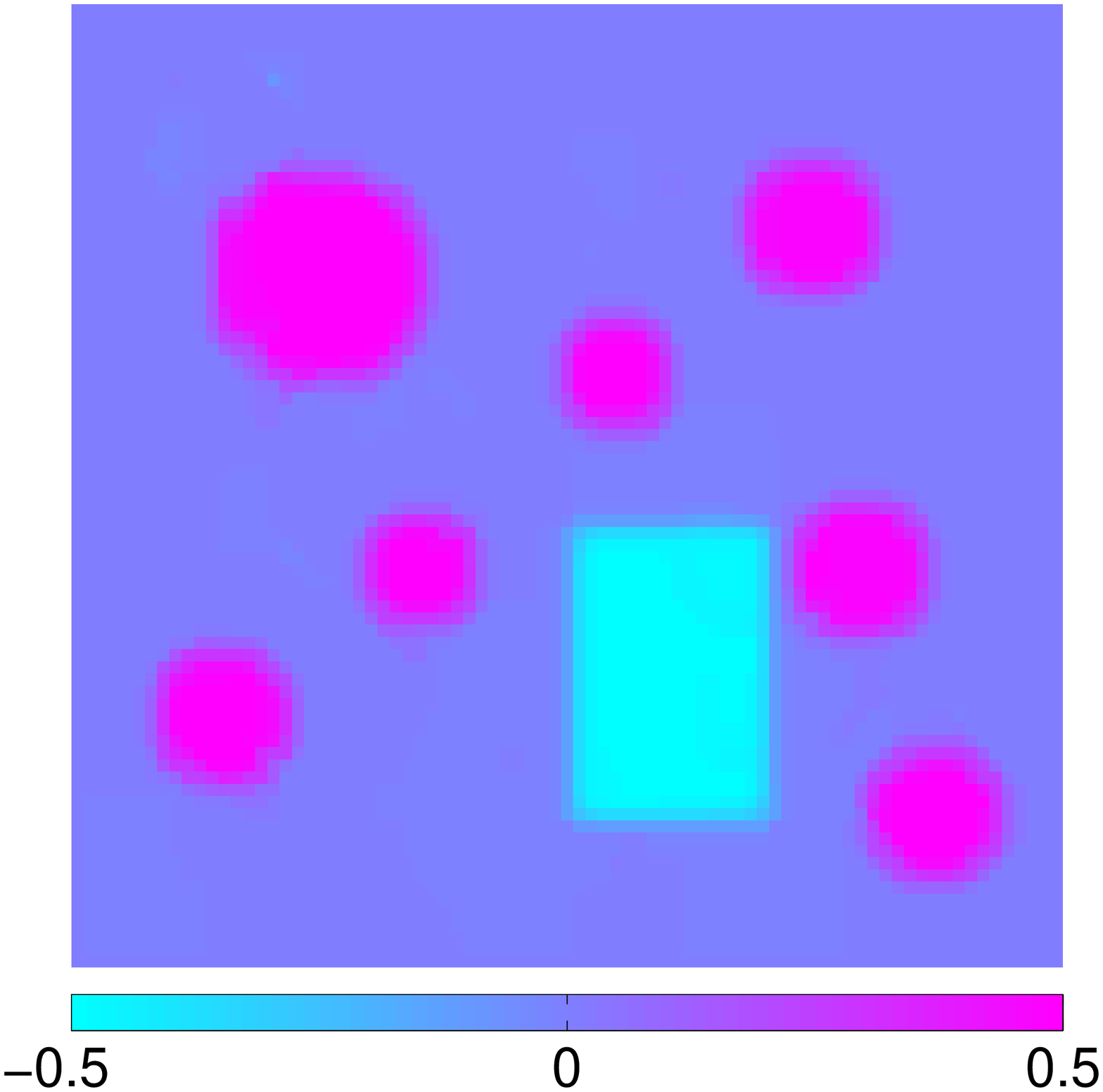}
     \label{ex5ctau}
     }
     \subfigure[$\zeta$ ($\alpha=1\%$)]{ 
     \includegraphics[width=37mm,height=35mm]{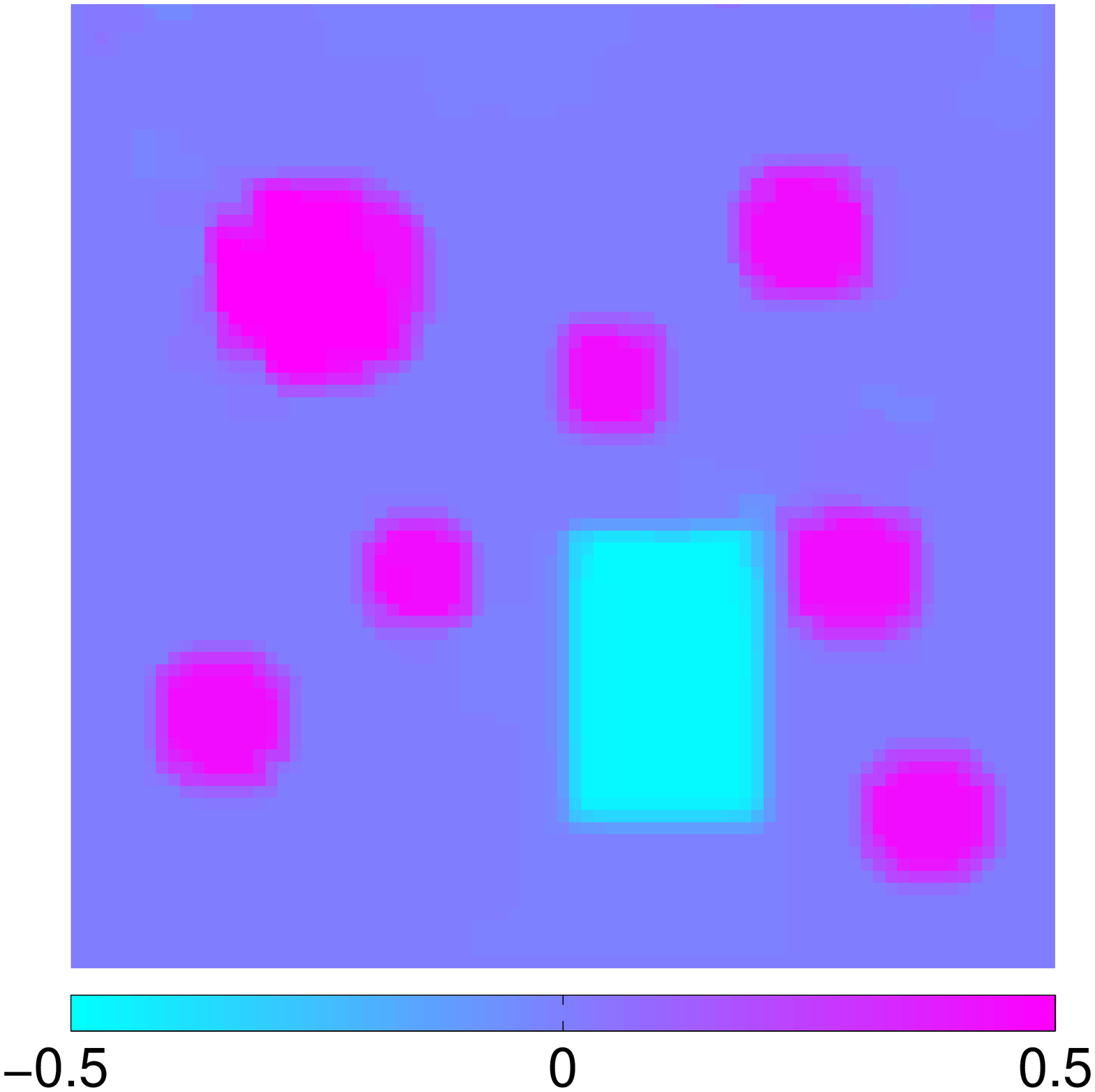}
     \label{ex5ntau}
     }
      \subfigure[$\zeta$ ($\alpha=4\%$)]{ 
     \includegraphics[width=37mm,height=35mm]{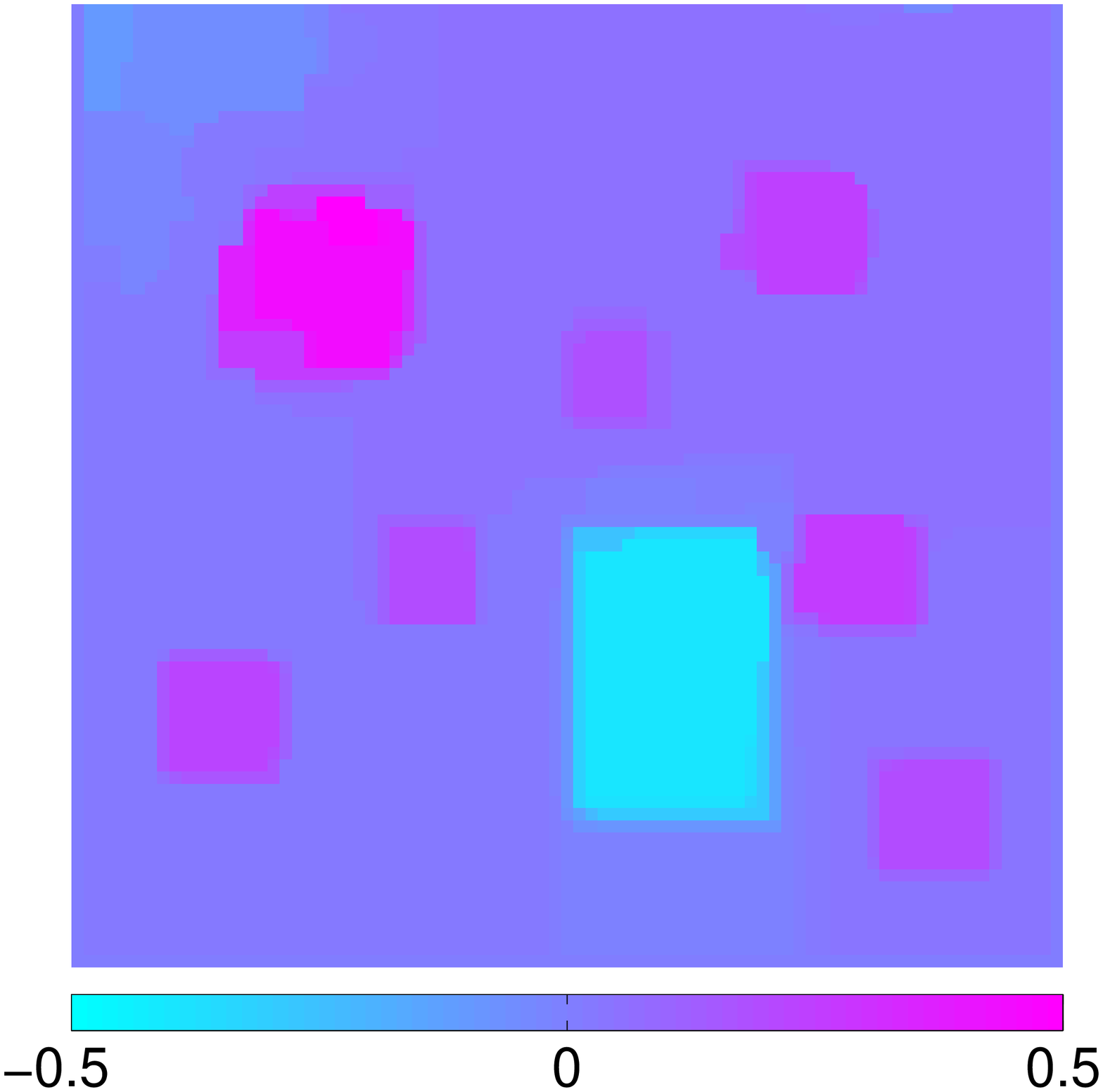}
     \label{ex5ntau4}
     }
    \subfigure[$\zeta$ at $\{y=-0.5\}$]{ 
     \includegraphics[trim=10mm 5mm 10mm 0mm,clip,width=35mm,height=38mm]{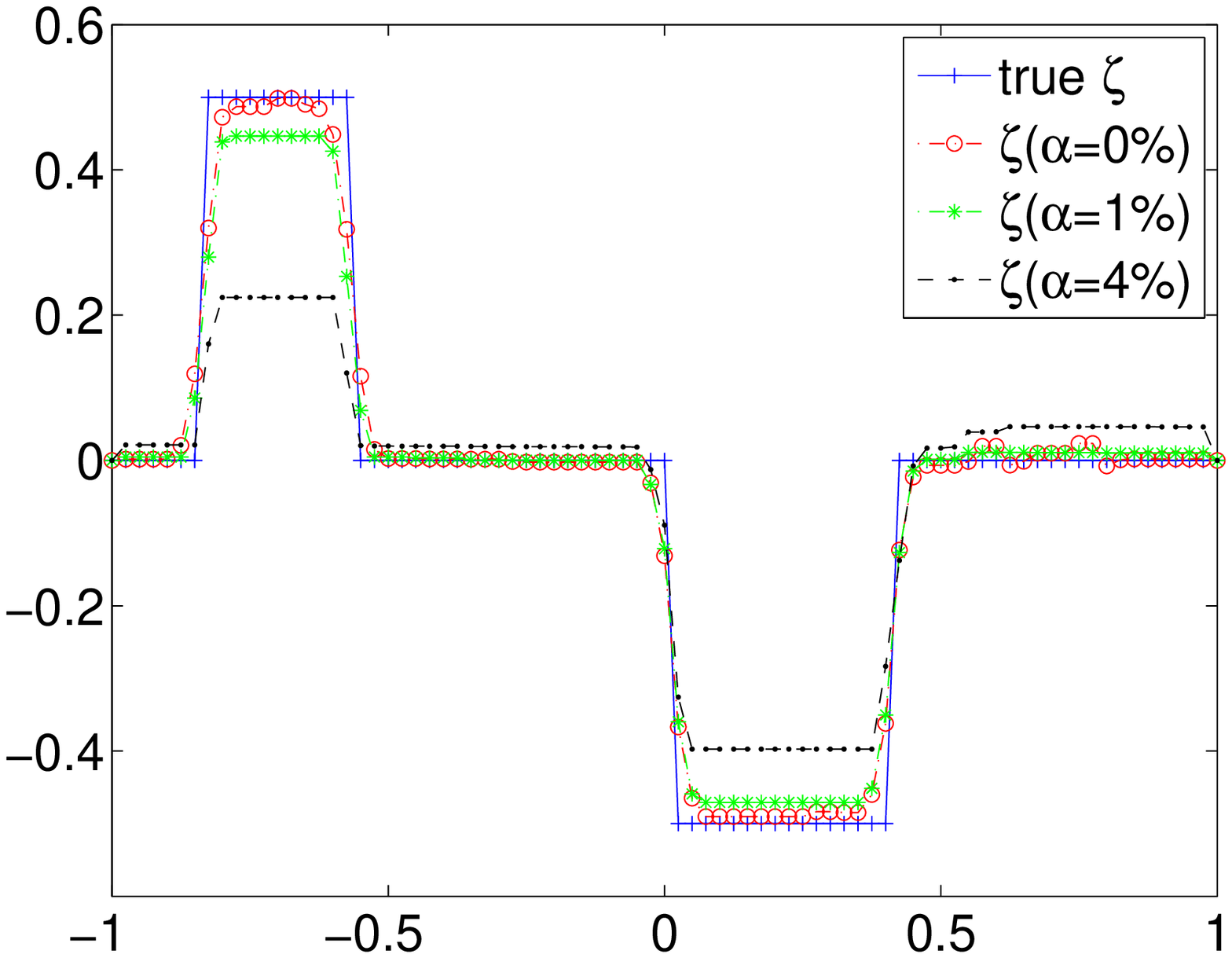} 
     \label{ex5rtau}
     }
     \subfigure[$|\gamma|^{\frac{1}{2}}$ ($\alpha=0\%$)]{
     \includegraphics[width=37mm,height=35mm]{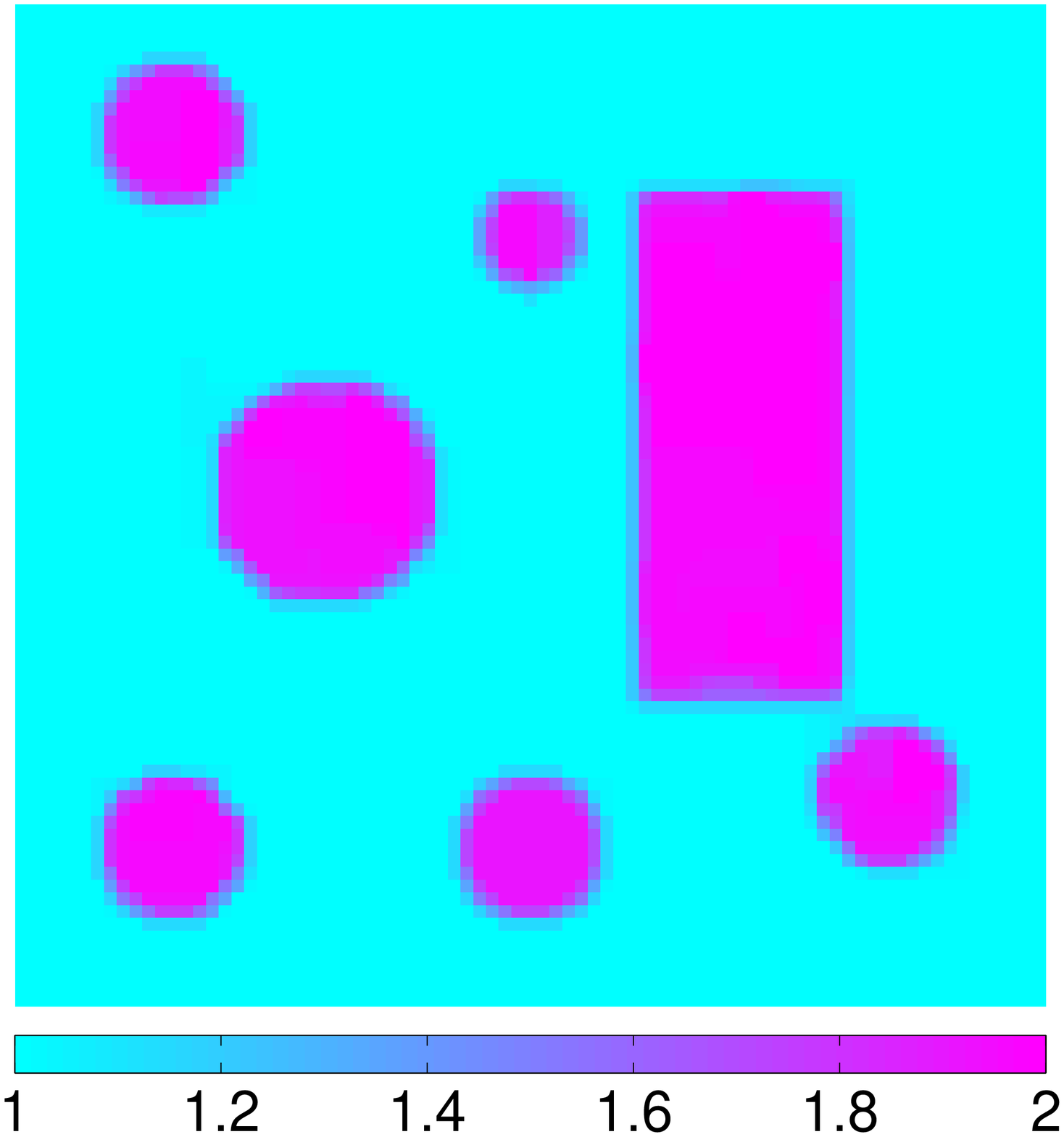}
     \label{ex5cbeta}
     }
    \subfigure[$|\gamma|^{\frac{1}{2}}$ ($\alpha=1\%$)]{
     \includegraphics[width=37mm,height=35mm]{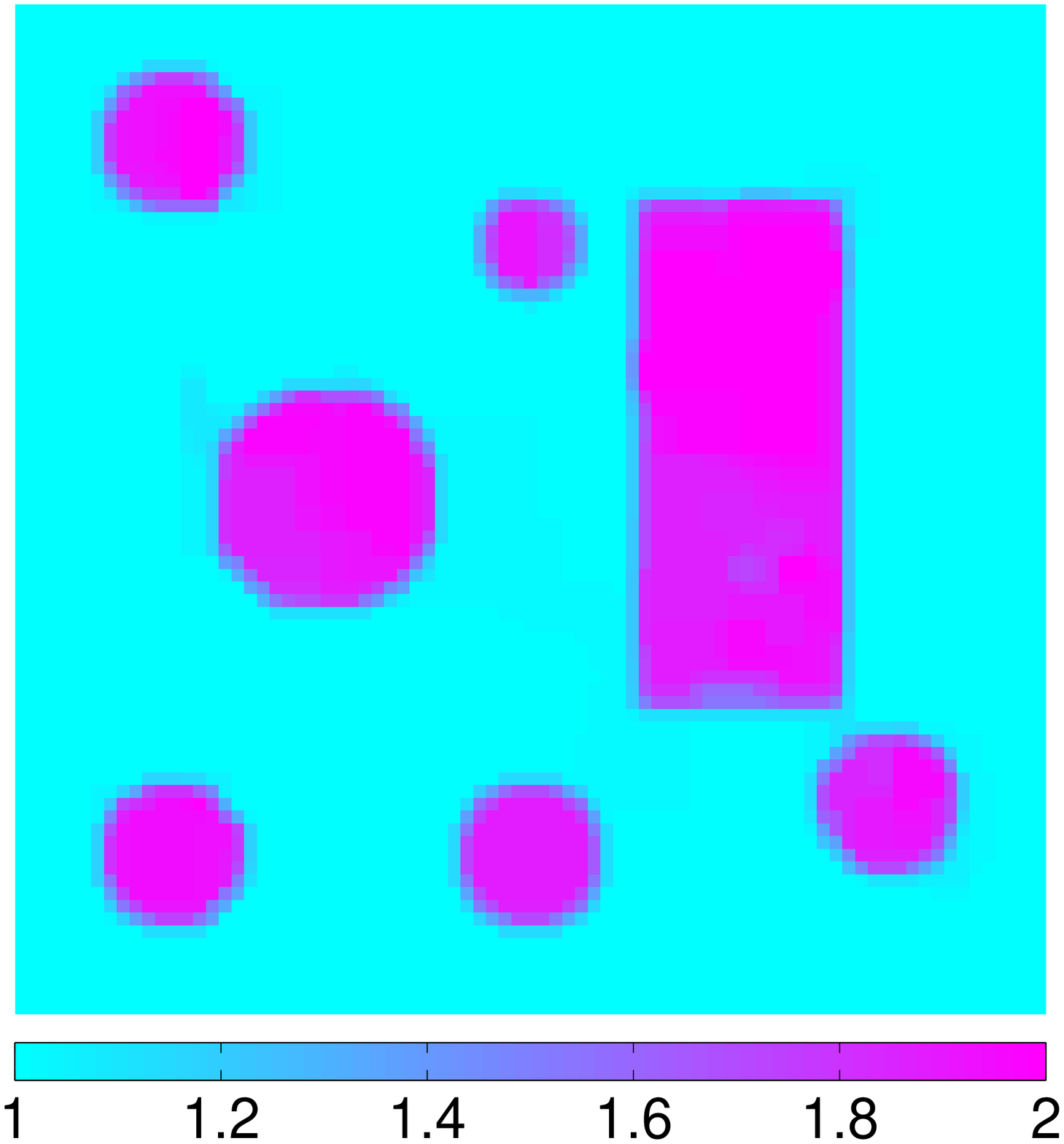}
     \label{ex5nbeta}
     }
     \subfigure[$|\gamma|^{\frac{1}{2}}$ ($\alpha=4\%$)]{
     \includegraphics[width=37mm,height=35mm]{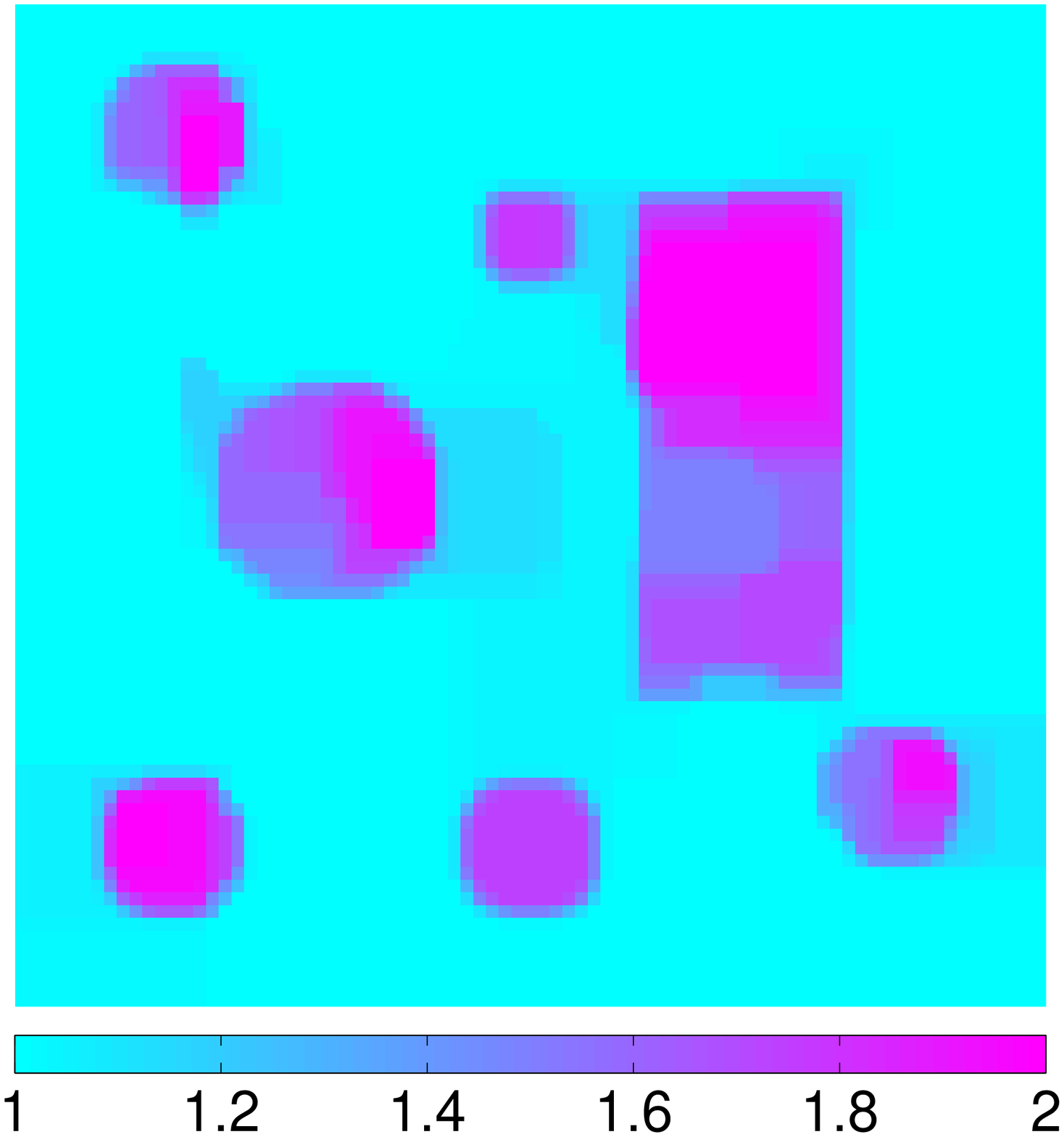}
     \label{ex5nbeta4}
     }
    \subfigure[$|\gamma|^{\frac{1}{2}}$ at $\{y=-0.5\}$]{ 
     \includegraphics[trim=15mm 5mm 14mm 0mm,clip,width=35mm,height=38mm]{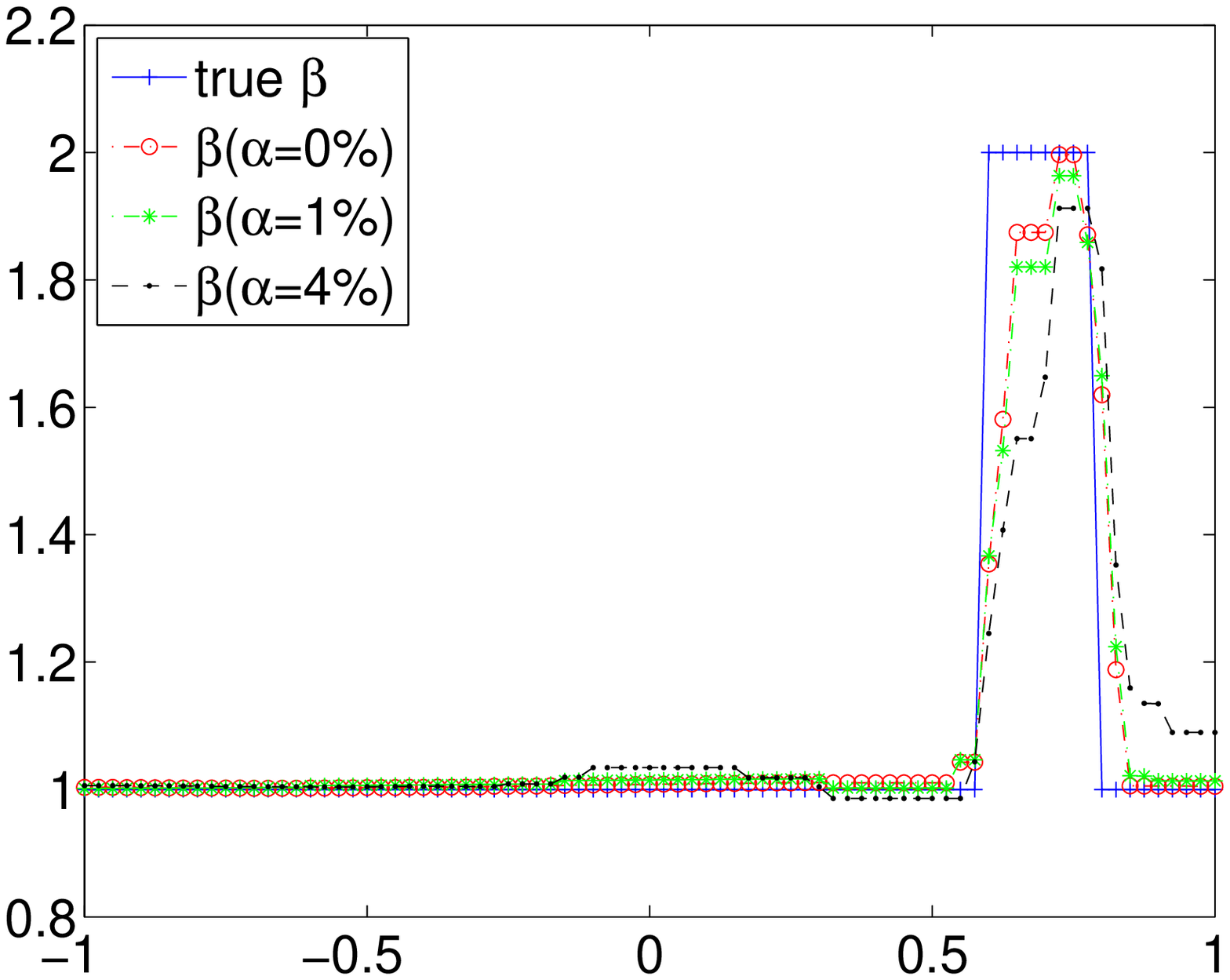} 
     \label{ex5rbeta}
     }
     \caption{Experiment 5. \subref{ex5cxi}\&\subref{ex5ctau}\&\subref{ex5cbeta}: reconstructions with noiseless data. \subref{ex5nxi}\&\subref{ex5ntau}\&\subref{ex5nbeta}: reconstructions with noisy data ($\alpha=1\%$). \subref{ex5nxi4}\&\subref{ex5ntau4}\&\subref{ex5nbeta4}: reconstructions with noisy data ($\alpha=4\%$). \subref{ex5rxi}\&\subref{ex5rtau}\&\subref{ex5rbeta}: cross sections along $y=-0.5$.}
\label{E5Neumann}
\end{figure}

\section{Conclusion} \label{sec:conclu}
This work presents an explicit reconstruction procedure for an anisotropic conductivity tensor $\gamma=(\gamma_{ij})_{1\leq i,j\leq 2}$ from knowledge of current densities of the form $H=\gamma\nabla u$.
 

As explained in Theorem \ref{stability}, these reconstruction algorithms, displaying local reconstruction formulas with Lipschitz stability (with the loss of one derivative from the measurements to the reconstructed quantities) for the anisotropic part of $\gamma$ and Lipschitz stability (with no loss of derivatives) for $\det\gamma$, rely heavily on the ability to construct families of solutions of the conductivity equation with linearly independent gradients (i.e. conditions $A$ and $B$ in Lemma \ref{cond cgo}). As the experimenter pilots these solutions from the boundary, it is then necessary to find appropriate boundary conditions ensuring the linear independence criterion. These linear independence conditions can be directly estimated from the available internal functionals $\{H_j\}_j$ and additional measurements could then be considered if necessary. This method was used in Experiments 4 and 5.

We first prove in Lemma \ref{cond cgo} that, if one can control the entire boundary, then boundary conditions close to traces of Complex Geometrical Optics solutions will generate solutions satisfying conditions $A$ and $B$ throughout the domain. In fact, these conditions can be verified numerically for quite a large class of boundary conditions, such as for instance traces of well-chosen polynomials, and Experiments 1-3 in the numerics section illustrate the success of the method on full reconstruction of both smooth and discontinuous coefficients, as well as its robustness to noise. 

On the other hand, when one has control over only part of the boundary, there will inherently be a breakdown in the reconstruction near the part of the boundary that is not controlled, as homogeneous boundary conditions there will automatically violate the linear independence criterion. On the controlled part of the boundary, using solutions generated with peaked Gaussian profiles at various positions yields satisfactory reconstructions up to a certain depth. As seen numerically on Experiments 4 and 5, the region where reconstructions are stable can be improved by increasing the spacing between the Gaussian profiles.

\section*{Acknowledgment} 
This work was partially funded by AFOSR Grant NSSEFFFA9550-
10-1-0194 and NSF Grant DMS-1108608. FM acknowledges partial support from NSF grant DMS-1025372.

\end{document}